\documentclass[11pt]{amsart}
\usepackage{amsmath}
\usepackage{amsxtra}
\usepackage{amscd}
\usepackage{amsthm}
\usepackage{amsfonts}
\usepackage{amssymb}
\usepackage{eucal}

\newtheorem{thm}{Theorem}[section]
\newtheorem{conj}{Conjecture}
\newtheorem{cor}[thm]{Corollary}
\newtheorem{lem}[thm]{Lemma}
\newtheorem{prop}[thm]{Proposition}

\theoremstyle{remark}
\newtheorem{remark}[thm]{Remark}

\theoremstyle{definition}

\numberwithin{equation}{section}

\newcommand{\Ref}[1]{{$($\ref{#1}$)$}}
\newcommand{\bean}{\begin{eqnarray}}
\newcommand{\eean}{\end{eqnarray}}
\newcommand{\be}{\begin{displaymath}}
\newcommand{\ee}{\end{displaymath}}
\newcommand{\bea}{\begin{eqnarray*}}   
\newcommand{\eea}{\end{eqnarray*}}

\newcommand{\thmref}[1]{Theorem~\ref{#1}}

\newcommand{\secref}[1]{Section~\ref{#1}}
\newcommand{\lemref}[1]{Lemma~\ref{#1}}
\newcommand{\propref}[1]{Proposition~\ref{#1}}

\newcommand{\nc}{\newcommand}

\nc{\bs}{\boldsymbol}

\nc{\arr}{\rightarrow}
\nc{\larr}{\longrightarrow}
\nc{\la}{\lambda}
\nc{\al}{\alpha}
\nc{\La}{\Lambda}
\nc{\ep}{\epsilon}
\nc{\si}{\sigma}
\nc{\om}{\omega}
\nc{\ga}{\gamma}

\nc{\Z}{{\mathbb Z}}
\nc{\C}{{\mathbb C}}

\nc{\g}{{\mathfrak g}}
\nc{\h}{{\mathfrak h}}
\nc{\n}{{\mathfrak n}}
\nc{\N}{\widehat{\n}}
\nc{\G}{\widehat{\g}}

\nc{\ot}{\otimes}
\nc{\wt}{\widetilde}
\nc{\wh}{\widehat}
\nc{\bi}{\bibitem}

\nc{\on}{\operatorname}
\nc{\Res}{\on{Res}}
\nc{\res}{\on{res}}

\nc{\slt}{\widehat{{\mathfrak s}{\mathfrak l}}_2}
\nc{\sw}{{\mathfrak s}{\mathfrak l}}
\nc{\gw}{{\mathfrak g}{\mathfrak l}}
\nc{\glt}{\widehat{{\mathfrak g}{\mathfrak l}}_2}
\nc{\uqgg}{U_q \g}
\nc{\mc}{\mathcal}
\nc{\uqbb}{U_q \wt{\g}}
\nc{\uqslt}{U_q \widehat{\sw}_2}
\nc{\uh}{U \wt{{\mathfrak h}}}
\nc{\uhh}{U \wh{{\mathfrak h}}}
\nc{\asln}{\widehat{\sw}_N}
\nc{\agln}{\widehat{\mathfrak g\mathfrak l}_N}
\nc{\aglm}{\widehat{\mathfrak g\mathfrak l}_M}
\nc{\asli}{\widehat{\sw}_infty}
\nc{\agli}{\widehat{\mathfrak g\mathfrak l}_\infty}
\nc{\sln}{{\sw}_N} 
\nc{\gln}{{\mathfrak g\mathfrak l}_N}
\nc{\glm}{{\mathfrak g\mathfrak l}_M}
\nc{\sli}{{\sw}_\infty}
\nc{\gli}{{\mathfrak g}{\mathfrak l}_\infty}
\nc{\SLi}{{SL}^-_\infty}
\nc{\SLin}{{SL}^{(N),-}_\infty}
\nc{\wSLi}{\wt{SL}^-_\infty}
\nc{\wSLin}{\wt{SL}^{(N),-}_\infty}
\nc{\GLl}{\widehat{GL}^-_l}
\nc{\GLln}{\widehat{GL}^{(N),-}_l}
\nc{\SLl}{\widehat{SL}^-_l}
\nc{\wGLl}{\wt{GL}^-_l}
\nc{\wGLln}{\wt{GL}^{(N),-}_l}
\nc{\wSLl}{\wt{SL}^-_l}

\nc{\asll}{\widehat{\mathfrak s\mathfrak l}_l}
\nc{\ur}{U_\ep^{\rm res}}
\nc{\uqr}{U_q^{\rm res}}
\nc{\urs}{U_{\ep^*}^{\rm res}}
\nc{\dur}{\dot{U}_\ep^{\rm res}}
\nc{\duqr}{\dot{U}_q^{\rm res}}
\nc{\durs}{\dot{U}_{\ep^*}^{\rm res}}
\nc{\fu}{U^{\on{fin}}_\ep}

\nc{\mb}{\mathbf}
\nc{\Rep}{\on{Rep}}

\nc{\M}{\mathcal M}
\nc{\gll}{{\mathfrak g}{\mathfrak l}_l}
\nc{\ol}{\overline}
\nc{\sll}{{\mathfrak s}{\mathfrak l}_l}
\nc{\zz}{{\mathfrak z}}
\nc{\agll}{\widehat{{\mathfrak g}{\mathfrak l}}_l}
\nc{\W}{\mc W}
\def\bin[#1;#2]{\left[\begin{matrix}{\displaystyle
#1}\\{\displaystyle #2}\end{matrix}\right]}

\def\qbin2[#1;#2;#3]{\left[\begin{matrix}{\displaystyle
#1;\displaystyle #2}\\{\displaystyle #3}\end{matrix}\right]}

\nc{\ds}{\displaystyle}
\nc{\Zl}{\Z_{(l)}}

\begin{document}

\title{The Hopf algebra $\Rep\, U_q \agli $}

\author {Edward Frenkel}

\author{Evgeny Mukhin}

\address{Department of Mathematics, University of California, Berkeley, CA
94720, USA}

\date{March 2001. Revised August 2002.}

\begin{abstract}
We define the Hopf algebra structure on the Grothendieck group of
finite-dimensional polynomial representations of $U_q \agln$ in the
limit $N \to \infty$. The resulting Hopf algebra $\Rep\, U_q \agli$ is
a tensor product of its Hopf subalgebras $\Rep_a U_q \agli$,
$a\in\C^\times/q^{2\Z}$. When $q$ is generic (resp., $q^2$ is a
primitive root of unity of order $l$), we construct an isomorphism
between the Hopf algebra $\Rep_a U_q \agli$ and the algebra of regular
functions on the prounipotent proalgebraic group $\wSLi$ (resp.,
$\wGLl$).  When $q$ is a root of unity, this isomorphism identifies
the Hopf subalgebra of $\Rep_a U_q \agli$ spanned by the modules
obtained by pullback with respect to the Frobenius homomorphism with
the algebra generated by the coefficients of the determinant of an
element of $\wGLl$ considered as an $l \times l$ matrix over the
Taylor series. This gives us an explicit formula for the Frobenius
pullbacks of the fundamental representations. In addition, we
construct a natural action of the Hall algebra associated to the
infinite linear quiver (resp., the cyclic quiver with $l$ vertices) on
$\Rep_a U_q \agli$ and describe the span of tensor products of
evaluation representations taken at fixed points as a module over this
Hall algebra.
\end{abstract}

\maketitle
\tableofcontents

\section{Introduction}
\subsection{Background}

Let $\Rep\, GL_N$ be the Grothendieck group of finite-dimensional
polynomial representations of the group $GL_N(\C)$ (i.e., those that
arise in tensor powers of the defining $N$--dimensional representation
$V_{\omega_1}$). Then $\Rep\, GL_N$ is isomorphic to the algebra of
polynomials over $\Z$ in commuting variables $t_i$, where $t_i$ is
the class of $\bigwedge^i V_{\omega_1}$ ($i=1,\ldots,N$).  We have a
natural injection of algebras $\Rep\, GL_N \to \Rep\, GL_{N+1}$
sending $t_i$ to $t_i$ ($i=1,\ldots,N$). Let us set $\Rep\, GL_\infty
= \underset{\longrightarrow}\lim \; \Rep\, GL_N$. Then $\Rep\,
GL_\infty$ carries the structure of a Hopf algebra. The multiplication
comes from the operation of tensor product on representations of
$GL_N$, while the comultiplication is defined as follows. Given a
splitting $N=N_1 + N_2$, we have a natural group homomorphism
$GL_{N_1} \times GL_{N_2} \to GL_N$, and therefore a map $\Rep\, GL_N
\to \Rep\, GL_{N_1} \otimes \Rep\, GL_{N_2}$. One can show that for
large $N,N_1$ and $N_2$ these maps stabilize and give rise to a
well-defined homomorphism $\Delta: \Rep \,GL_\infty \to \Rep
\,GL_\infty \otimes \Rep \,GL_\infty$; this gives us the
comultiplication structure.

A well-known incarnation of the Hopf algebra $\Rep\, GL_\infty$ is the
algebra $\Lambda_\infty$ of symmetric polynomials in infinitely many
variables, in which the basis of irreducible representations
corresponds to the basis of Schur polynomials. The connection with
symmetric functions is due to the existence of the character
homomorphism $\Rep\, GL_N \to \Z[x_1,\ldots,x_N]^{S_N}$. The Hopf
algebra structure on $\Lambda_\infty$ may be easily described in terms
of the operations of multiplication and restriction of functions. The
fact that these operations make $\Lambda_\infty$ a Hopf algebra is
very useful in the theory of symmetric functions (see \cite{Mac}).

Let $q$ be a nonzero complex number.  The goal of the present paper is
to study a ``$q$--analogue'' of the Hopf algebra $\Rep\, GL_\infty$
which, in fact, depends only on the multiplicative order of
$q\in\C^\times$. It is obtained by replacing the groups $GL_N$ with
the quantized universal enveloping algebras $U_q\agln$, where $\agln$
is the affine Lie algebra of $\gln$. As in the case of $GL_N$, the
Grothendieck ring of $U_q\agln$ may be identified with the ring of
polynomials over $\Z$ in commuting variables corresponding to the
evaluation representations $V_{\omega_i}(b)$, where $i=1,\ldots,N$ and
$b \in \C^\times$. This algebra decomposes into a tensor product of
its subalgebras $\Rep_a U_q\agln, a \in \C^\times$, generated by
$V_{\omega_i}(b)$, where $b \in a q^{2\Z}$. The representations
corresponding to different factors do not interact, so without loss of
generality we may concentrate on $\Rep_a U_q\agln$, which is the ring
of polynomials in commuting variables $t_{i,n} =
V_{\omega_i}(aq^{2n})$, where $i=1,\ldots,N$ and $n \in \Z$. We have
natural injections $\Rep_a U_q\agln \to \Rep_a U_q \wh\gw_{N+1}$, and
so we define $\Rep_a U_q\agli$ as $\underset{\longrightarrow}\lim \;
\Rep_a U_q\agln$.

Then $\Rep_a U_q\agli$ is naturally a Hopf algebra, with the operations
of multiplication and comultiplication defined in a similar way to the
case of $GL_\infty$. This Hopf algebra and its versions when $q$ is
specialized to be a root of unity are the central objects of this
paper.

The ring $\Rep_a U_q\agln$ comes equipped with an injective
$q$--character homomorphism $\Rep_a U_q\agln \to
\Z[\La_{i,aq^{2n}}]_{i=1,\ldots,N}^{n \in \Z}$, defined as in
\cite{FR}, \cite{FM2} (see also \cite{FM1}). This homomorphism plays
the role of the ordinary character homomorphism from $\Rep\, GL_N$ to
the ring of symmetric polynomials. Thus, $\Rep_a U_q\agli$ may be
thought of as a ``$q$--analogue'' of the Hopf algebra $\Lambda_\infty$
of symmetric polynomials in infinitely many variables. In fact, since
$\Rep\, GL_\infty$ may be identified with $\Rep\, \gw_\infty$, which
in turn coincides with $\Rep_a U \wh{\gw}_\infty$, the Hopf algebra
$\Lambda_\infty$ appears in the special case $q=1$ of $\Rep_a
U_q\agli$. We expect that many of the well-known structures on the
ring of infinite symmetric polynomials have interesting $q$--analogues
in $\Rep_a U_q\agli$. In this paper we make the first steps towards
understanding these structures.

Our interest in the structure of $\Rep_a U_q\agln$ comes in particular
from the fact that $\Rep_a U_q\agln \otimes \C$ may be thought of as a
classical limit of the two-parameter deformed $\W$--algebra
$\W_{q,t}(\gw_N)$ defined in \cite{FF,AKOS} (for more on this
interpretation, see \cite{FR:simple,FR}). In that sense, $\Rep_a
U_q\agli \otimes \C$ is a special case of the deformed $\W$--algebra
$\W_{q,t}(\gw_\infty)$. Much of the structure on $\Rep_a U_q\agli$
that we discuss in this paper may be generalized to the deformed
$\W$--algebras. Similar structures also exist for quantum affine
algebras of classical types. We plan to discuss these topics
elsewhere.

\subsection{Quantum affine algebras and affine Hecke algebras}

The ring of symmetric polynomials $\Lambda_\infty$ also has another
realization, as the direct sum $\oplus_{n\geq 0} \Rep\, S_n$, where
$\Rep S_n$ is the Grothendieck group of the symmetric group on $n$
letters. The two realizations are equivalent via the Schur duality
functor. From this point of view, a $q$--analogue of $\Lambda_\infty$
may be constructed by taking $\oplus_{n\geq 0} \Rep\, \wh{H}_{q^2,n}$,
where $\wh{H}_{q^2,n}$ is the affine Hecke algebra of $GL_n$. The Hopf
algebra structure on the latter has been introduced by J. Bernstein
and A. Zelevinsky \cite{BZ} and recently described explicitly by
I. Grojnowski \cite{Gr}. This description is also implicit in the
earlier work by S. Ariki \cite{A} (for related works, see \cite{LNT}
and references therein).

When $q$ is not a root of unity, the Hopf algebras $\Rep\, U_q\agli$
and $\oplus_{n\geq 0} \Rep\, \wh{H}_{q^2,n}$ are isomorphic via the
affine $q$--Schur duality functor introduced by M. Jimbo \cite{J} and
studied by V. Chari--A. Pressley \cite{CP:schur}, and
V. Ginzburg--N. Reshetikhin--E. Vasserot \cite{GRV}. Using this
functor, most of our results for generic $q$ on the structure of
$\Rep\, U_q\agli$ may be derived from those of \cite{A,Gr,LNT}.

However, when $q$ becomes a root of unity $\ep$, the Hopf algebra
$\oplus_{n\geq 0} \Rep\, \wh{H}_{\ep^2,n}$ is different from $\Rep\,
\ur \agli$, where $\ur \agln$ is the restricted specialization of $U_q
\agln$ at $q = \ep$. This can be seen, for example, from the fact that
irreducible representations of affine Hecke algebras at roots of unity
are parameterized only by those multisegments (in the language of
\cite{BZ}) which are aperiodic. In fact, the restricted specialization
$\ur \agln$ contains a subalgebra $\fu \agln$, which we call the small
quantum affine algebra. The direct limit of Grothendieck groups $\Rep\,
\fu \agli$ may be identified with the quotient of $\Rep\, \ur \agli$
by the ideal generated by the pullbacks of the Frobenius
homomorphism, and carries a natural Hopf algebra structure. It is this
quotient that is isomorphic to $\oplus_{n\geq 0} \Rep\,
\wh{H}_{\ep^2,n}$.

In this paper we do not use affine Hecke algebras. One of the
advantages of the Hopf algebra $\Rep U_q\agli$ over $\oplus_{n\geq 0}
\Rep\, \wh{H}_{q^2,n}$ is in the existence of the $q$--character
homomorphism. This homomorphism makes the interpretation of $\Rep
U_q\agln$ as an analogue of the Hopf algebra of symmetric polynomials
especially transparent.

\subsection{Main results}

Our first result is an identification of each of the three Hopf
algebras: $\Rep_a U_q \agli$, $\Rep_a \ur \agli$, and $\Rep_a \fu
\agli$ with well-known Hopf algebras.

Let $\wSLi$ be the group of all matrices $M = (M_{i,j})_{i,j \in \Z}$
with integer entries, such that $M_{i,i} = 1$ for all $i$, and
$M_{i,j} = 0$ for all $i<j$. It contains the subgroup $\wGLl$ of those
matrices which satisfy: $M_{i+l,j+l} = M_{i,j}$. Finally, $\wGLl$
contains a normal subgroup $\wSLl$ defined by the equations $a_i=0$
($i>0$), where $a_i$ is the $i$th coefficient of the determinant of a
matrix from $\wGLl$ considered as a subgroup of the formal loop group
of $GL_l$ (see \secref{new sect}).  Thus, we have a sequence of
embeddings
\begin{equation}    \label{3 groups}
\wSLl \hookrightarrow \wGLl \hookrightarrow \wSLi.
\end{equation}
All three groups are groups of $\Z$--points of prounipotent
proalgebraic groups over $\Z$, which we we will denote by the same
symbols by abuse of notation. The ring $\Z[\wSLi]$ of regular
functions on $\wSLi$ is just the ring of polynomials in the variables
$M_{i,j}, i>j$. The ring $\Z[\wGLl]$ is its quotient by the relations
$M_{i+l,j+l} = M_{i,j}$, and $\Z[\wSLl]$ is the quotient of
$\Z[\wGLl]$ by the relations $a_i=0$ ($i>0$). Each of these rings has a
natural structure of Hopf algebra. In addition, the polynomial algebra
$\Z[a_i]_{i>0}$ is a Hopf subalgebra of $\Z[\wGLl]$, which we denote
by $\Z[\wt{\mc Z}^-_l]$.

\bigskip

\noindent {\bf Theorem.} {\it
Let $\ep$ be a nonzero complex number, such that $\ep^2$ has order
$l$. Then we have the following commutative diagram of Hopf algebras
$$
\begin{CD}
\Z[\wSLi] @>{\sim}>> \Rep_a U_q\agli \\
@VVV @VVV \\
\Z[\wGLl] @>{\sim}>> \Rep_a \ur \agli \\
@VVV @VVV \\
\Z[\wSLl] @>{\sim}>> \Rep_a \fu \agli \\
\end{CD}
$$
where the left vertical maps are induced by the embeddings \eqref{3
groups}, the upper right vertical map is induced by specialization of
representations to $q=\ep$, and the lower right vertical map is the
quotient of $\Rep_a U_q\agli$ by the augmentation ideal of the
subalgebra spanned by the Frobenius pullbacks. Moreover, the 
subalgebra spanned by the Frobenius pullbacks is identified with
$\Z[\wt{\mc Z}^-_l] \subset \Z[\wGLl]$. 
}

\bigskip

The top two isomorphisms in the above diagram are given by the formula
$M_{i,j} \mapsto t_{i-j,j}$, where $M_{i,j}$ is the function whose
value on a matrix $M$ in $\wSLi$ (resp., $\wGLl$) is equal to the
$(i,j)$ entry of $M$, and $t_{i,j}$ is the class of evaluation module
$V_{\om_i}(aq^{2n})$. We also obtain that the Frobenius pullback of
the $i$th fundamental representation is equal to the polynomial $a_i$
of the fundamental representations $t_{k,m}$.

The complexifications of the algebras of functions appearing in the
left column of the above diagram carry natural actions of the Lie
algebras $\sw^-_\infty, \wh{\gw}^-_l$, and $\wh{\sw}^-_l$,
respectively. The corresponding actions on the complexified
Grothendieck groups appearing in the right column of the above diagram
may be described quite explicitly in representation theoretic terms.

At the level of integral forms we obtain an action of the Hall
algebras ${\mc H}_\infty$ and ${\mc H}_l$ (associated to the infinite
linear quiver and the cyclic quiver with $l$ vertices), and the
integral form $U_\Z \wh{\sw}^-_l$ of $U \wh{\sw}^-_l$, respectively.
In fact, we prove that $\Rep_a U_q\agli$ (resp., $\Rep_a \ur \agli$)
and ${\mc H}_\infty$ (resp., ${\mc H}_l$) are restricted dual Hopf
algebras.

This duality, at the level of vector spaces after complexification,
has been obtained earlier by M. Varagnolo and E. Vasserot \cite{VV} by
a different method. More precisely, they have identified $\Rep_a
U_q\agli \otimes \C$ (resp., $\Rep_a \ur \agli \otimes \C$) with the
restricted dual of ${\mc H}_\infty \otimes \C$ (resp., ${\mc H}_l
\otimes \C$) by declaring the basis of standard modules in the former
to be dual to the PBW basis in the latter. Then, using the equivalence
between the geometric descriptions of representations of quantum
affine algebras \cite{GV,V} and the corresponding Hall algebras
\cite{Lu}, they showed that the basis of irreducible modules in the
former is dual to the canonical basis in the latter.

In contrast to \cite{VV}, we obtain all of our results by algebraic
methods. In fact, most of our proofs are rather elementary. The main
technical tool is the theory of $q$--characters. We hope that our
methods may be extended so as to enable us to identify the basis of
irreducible modules with the dual canonical basis in a purely
algebraic way.
 
\subsection{Contents} The paper is organized as follows.

In Section \ref{algebras for generic q} we start with the definition
of different realizations of quantized enveloping 
algebras when $q$ is not a root
of unity.  In Section \ref{isomorphisms section} we describe the
isomorphisms between these realizations. We try to fill the gaps in
the literature by collecting all relevant results (some unpublished)
with compatible normalizations.  The main result of this section is
Theorem \ref{cube theorem}.

In Section \ref{rep generic section} we discuss the polynomial
representations of $U_q\agln$ and introduce our main technical tool,
the $q$--characters. In this section (and in Sections \ref{algebras at
roots}, \ref{rep at roots}, \ref{frob section}) we extend to the case
of $U_q\agln$ some result that so far have only been established for
$U_q\asln$.

In Section \ref{ind limit def generic section} we define the inductive
limit $\Rep\,U_q\agli$ and a Hopf algebra structure on it. The main
results of this section are Corollary \ref{irrep form basis} and
Theorem \ref{hopf str}.  In Section \ref{ind lim structure generic} we
study the Hopf algebra $\Rep\,U_q\agli$. The main results are stated
in Theorems \ref{repa to uminus}, \ref{action of sl} and Proposition
\ref{many gl to agl lemma}.

In Section \ref{algebras at roots} we describe the quantum groups at
roots of unity. In Section \ref{rep at roots} we describe the
polynomial representations of $\ur\agln$ and the corresponding
$\ep$--characters. Sections \ref{algebras at roots} and \ref{rep at
roots} are the root of unity versions of Sections \ref{algebras for
generic q} and \ref{rep generic section}.

In Section \ref{frob section} we study the Frobenius homomorphism and
the modules obtained by the pullback with respect to the Frobenius
homomorphism.

In Section \ref{hall section} we present a self-contained analysis of
the Hall algebras associated to the infinite linear quiver and the
cyclic quivers when the deformation parameter is equal to $1$. We also
define the algebras of functions on the groups $\wGLl, \wSLl$ and
$\wt{\mc Z}_l^-$, and various decompositions and dualities.

In Section \ref{ind str root} we define and study the Hopf algebra
$\Rep\,\ur\agli$. This section is the root of unity version of
Sections \ref{ind limit def generic section} and \ref{ind lim
structure generic}.  The main results are given by Lemma \ref{can
lem}, Theorems \ref{repa to uminus root}, \ref{hall action}, \ref{fr
ort sl} and \ref{sl ort fr}, and Propositions \ref{many gl to agl root
lemma} and \ref{many gl to agl root factor lemma}.

Finally we give a summary of the general picture in Section
\ref{summary section}, see Theorem \ref{summary thm}.

\subsection{Acknowledgments} We thank V.~Chari, I.~Cherednik, M.~Jimbo,
A. Okounkov, V.~Tarasov, and especially N.~Reshetikhin for useful
discussions. The first named author thanks M.-F. Vigneras for
inspiring discussions of the connection between representations of
quantum affine algebras and $p$--adic groups in non-defining
characteristic.

This research was supported in part by E.~Frenkel's fellowship from
the Packard Foundation and NSF grant DMS-0070874.

\section{Quantum algebras for generic $q$}\label{algebras for generic
  q}

In this section we collect information about different realizations of
the quantized enveloping
algebras $U_q\sln$, $U_q\gln$, $U_q\asln$, $U_q\agln$ and relations among
them. Though most of this material is known to experts, the results in
the literature are not complete, scattered in many places and have
different normalizations.

In this paper we deal exclusively with finite-dimensional
representations of quantum affine algebras, all of which have level
zero. The central element of the quantum affine algebra acts as the
identity on any such representation. To simplify our formulas, we will
impose the additional relation that this central element is equal to
$1$ (in all realizations).

In Sections \ref{algebras for generic q}, \ref{rep generic section},
\ref{ind limit def generic section}, \ref{ind lim structure generic},
we consider quantized enveloping algebras over the field of complex
numbers, and $q$ is a nonzero complex number which is not a root of
unity. These algebras may instead be considered over the field of
rational functions $\C(q)$. This point of view is adopted in the
remaining sections of the paper.

\subsection{The algebra $U_q\asln$}\label{asln} 

In this section we recall two realizations of $U_q\asln$. 

Let $C=(C_{ij})_{i,j=1,\dots,N-1}$ be the Cartan matrix and
$\{\al_i\}_{i=1,\dots, N-1}$ the set of simple roots of $\sln$. Let
$C=(C_{ij})_{i,j=0,\dots,N-1}$ be the Cartan matrix and
$\{\al_i\}_{i=0,\dots, N-1}$ the set of simple roots of $\asln$.  Set
\be [n]_q = \frac{q^n - q^{-n}}{q - q^{-1}},\qquad
[n]_q!=[1]_q[2]_q\dots[n]_q.  \ee 

The quantized enveloping algebra of $\asln$ 
in the {\em Drinfeld-Jimbo realization}
\cite{Dr1,J} is an associative algebra, denoted
by $U_q^{\on{DJ}} \asln$, with generators $X_i^{{}\pm{}}$, $K_i^{{}\pm
1}$ ($i=0,\dots,N-1$), and relations:
\begin{align*}
&K_0K_1\dots K_{N-1}=1, \qquad \qquad\qquad K_iK_j = K_jK_i,&\\
&K_iX_j^{{}\pm{}}K_i^{-1} = q^{{}\pm C_{ij}}X_j^{{}\pm},
\qquad\qquad\;\;\;[X_i^+ , 
X_j^-] = \delta_{ij}\frac{K_i - K_i^{-1}}{q -q^{-1}},&\\
&X_i X_j = X_j X_i \qquad\qquad\qquad\qquad \quad(|i-j|\neq 1 \mod N),&\\
&[X_i^\pm,[X^\pm_i,X^\pm_j]_{q^{-1}}]_q=0\qquad \qquad\;\;\;(|i-j|=1 \mod N, \;\; N>2),&\\
&[X^\pm_i,[X^\pm_i,[X^\pm_i,X^\pm_j]_{q^{-2}}]]_{q^2}=0 \qquad\, (
|i-j|=1 \mod N,\;\;N=2).&
\end{align*}
Here we used the following notation for commutators:
\be
[A,B]_a=AB-aBA,\qquad [A,B]=[A,B]_1.
\ee

Note that the first relation identifies the central element of
$U_q\asln$ with $1$.

Denote the subalgebra of $U_q\asln$ generated by $K_i^{\pm 1}, X_i^+$
(respectively, $K_i^{\pm 1}, X_i^-$), $i=0,\dots,N-1$, by $U_q {\mathfrak
b}_+$ (resp., $U_q {\mathfrak b}_-$).

The algebra $U_q^{\on{DJ}} \asln$ has the structure of a Hopf
algebra. The comultiplication $\Delta$ on the generators is given by
the formulas 
\begin{align}\label{sln comult}
\Delta(K_i) &= K_i \otimes K_i,\notag \\ \Delta(X^+_i)
&= X^+_i \otimes 1 + K_i \otimes X^+_i,\\ \Delta(X^-_i)&= X^-_i
\otimes K_i^{-1} + 1 \otimes X^-_i.\notag  
\end{align}

\medskip

Now consider the Drinfeld ``new''
realization of quantized enveloping algebra of $\asln$. 
This is the algebra, denoted by $U^{\rm new}_q\asln$,
with generators $x_{i,n}^{\pm}$, $k_i^{\pm 1}$ ($i=1,\dots,N-1$,
$n\in\Z$), $h_{i,n}$ ($i\in 1,\dots,N-1$, $n\in \Z\backslash 0$) and
the following relations:
\begin{align*}
k_ik_j = k_jk_i,\quad  k_ih_{j,n}& =h_{j,n}k_i,
\quad h_{i,n}h_{j,m}=h_{j,m} h_{i,n},\\
k_ix^\pm_{j,n}k_i^{-1} = q^{\pm C_{ij}}x_{j,n}^{\pm},&\quad
[h_{i,n} , x_{j,m}^{\pm}] = \pm \frac{1}{n} [n C_{ij}]_q x_{j,n+m}^{\pm},\\
[x_{i,n}^+, x_{j,m}^-] &=\delta_{ij}
\frac{\phi_{i,n+m}^+ - \phi_{i,n+m}^-}{q - q^{-1}},\\
[x_{i,n+1}^{\pm},x_{j,m}^{\pm}]_{q^{\pm C_{ij}}}&=
-[x_{j,m+1}^{\pm}x_{i,n}^{\pm}]_{q^{\pm C_{ij}}},\\
[x^\pm_{i,n} , [x^\pm_{j,k},x^\pm_{i,m}]_q]_q&+[x^\pm_{i,m} ,
[x^\pm_{j,k},x^\pm_{i,n}]_q]_q = 0 \qquad (|i-j|=1).
\end{align*}
Here $\phi_{i,n}^{\pm}$'s are determined by the formula
\be\label{series}
\Phi_i^\pm(u) := \sum_{n=0}^{\infty}\phi_{i,\pm n}^{\pm}u^{\pm n} =
k_i^{\pm 1} \exp\left(\pm(q-q^{-1})\sum_{m=1}^{\infty}h_{i,\pm m}
u^{\pm m}\right).
\ee

Here we again set the central element equal to $1$.

The algebras $U_q^{\on{DJ}} \asln$ and $U^{\rm new}_q\asln$ are
isomorphic, under the isomorphism $\on {B}: U^{\rm new}_q\asln\to
U_q^{\on{DJ}} \asln$ described in Section \ref{new-DJ}.
We will therefore use a single notation $U_q \asln$ for both of them.

Lusztig has introduced in \cite{L}, \S 37.1.3 four families of algebra
automorphisms of $U_q\asln$, 
for which we will use the notation ${\mc B}_i^{(\pm,
j)}$, where $j=1,2$, $i=0,1,\dots, N-1$ (in \cite{L} they were denoted by
$T'_{i,e}$ and $T''_{i,e}, e=\pm 1$). For $N>2$ they are given in the
Drinfeld-Jimbo realization by the
formulas \begin{align}\label{braid1} {\mc B}_i^{(\pm, 1)}
(X_i^+)&=-K_i^{\pm1}X_i^-,\qquad {\mc B}_i^{(\pm ,1)}
(X_i^-)=-X_i^+K_i^{\mp1},\notag \\ {\mc B}_i^{(\pm ,1)}
(X_j^+)&=X_j^+X_i^+-q^{\pm 1}X_i^+X_j^+ \qquad (j=i-1,i+1),\\ {\mc
B}_i^{(\pm, 1)} (X_j^-)&=X_i^-X_j^--q^{\mp 1}X_j^-X_i^-\qquad
(j=i-1,i+1),\notag\\ {\mc B}_i^{(\pm, 1)} (X_j^+)=X_j^+&,\qquad {\mc
B}_i^{(\pm, 1)} (X_j^+)=X_j^+ \qquad (i\neq j,j-1,j+1) \notag \end{align} and
\begin{align}\label{braid2} {\mc B}_i^{(\pm, 2)} (X_i^+)&=-X_i^-K_i^{\mp
1},\qquad {\mc B}_i^{(\pm, 2)} (X_i^-)=-K_i^{\pm 1}X_i^+,\notag\\ {\mc
B}_i^{(\pm, 2)} (X_j^+)&=X_i^+X_j^+-q^{\pm 1}X_j^+X_i^+ \qquad
(j=i-1,i+1),\\ {\mc B}_i^{(\pm, 2)} (X_j^-)&=X_j^-X_i^--q^{\mp
1}X_i^-X_j^-\qquad (j=i-1,i+1),\notag\\ {\mc B}_i^{(\pm, 2)}
(X_j^+)=X_j^+&,\qquad {\mc B}_i^{(\pm, 2)} (X_j^+)=X_j^+ \qquad (i\neq
j,j-1,j+1),\notag \end{align} where we consider indices modulo
$N$. The automorphisms ${\mc B}^{(\pm, j)}_i$ are called the braid
group automorphisms.

Also let $\mc D$ be the automorphism of $U_q\asln$ defined by
\bean\label{mc T}
\mc D(X_i^\pm)=X_{i+1}^\pm,\qquad \mc D(K_i^{\pm1})=K_{i+1}^{\pm1}
\qquad (i=0,\dots, N-1),
\eean
where we again consider indices modulo $N$.

The algebra $U_q\sln$ is the Hopf subalgebra of 
$U^{\rm DJ}_q\asln$ generated by
$K_i^{\pm 1}$, $X_i^\pm$ ($i=1,\dots,N-1$). The algebra $U_q\sln$ is
preserved by $\mc B_i^{(\pm, j)}$ ($j=1,2$, $i=1,\dots,N-1$).

\subsection{Universal $R$--matrix}    \label{rmatr}

Let $(\cdot,\cdot)$ be the restriction of the invariant bilinear form
on $\sln$ to the Cartan subalgebra, 
normalized so that the square of any root is equal to
$2$. Denote by $P$ the weight lattice of $\sw_N$. Then $(\cdot,\cdot)$
induces a pairing $P \times P \to \frac{1}{N} \Z$.

A vector $w$ in a $U_q\sln$--module $W$ is called a vector of weight
$\bs \la \in P$, if $K_i \cdot w = q^{(\bs \la,\al_i)} w$,
$i=1,\dots,N-1$.  A representation $W$ of $U_q\sln$ is said to be of
type 1 if it is the direct sum of its weight spaces $W = \oplus_{\bs
\la} W_{\bs\la}$, where $W_{\bs \la} = \{w\in W | K_i \cdot w =
q^{(\bs\lambda,\al_i)}w\}$.

A representation $V$ of $U_q\asln$ is called of type 1 if $V$ is of
type 1 as a representation of $U_q\sln$. Every finite-dimensional
irreducible representation of $U_q\asln$ can be obtained from a type 1
representation by twisting with an automorphism of $U_q\asln$. In this
paper we consider only finite-dimensional type 1 representations.

Fix an $N$th root of $q$.  Given a pair $V, W$ of $U_q\sln$--modules
of type 1, define an operator $q^{H\otimes H}$ on $V \otimes W$, which
acts by multiplication by $q^{(\bs \la,\bs \mu)}\in q^{\frac{1}{N}\Z}$
on the subspace $V_{\bs \la} \otimes W_{\bs \mu}$. Informally, one can
think that $K_i=q^{H_i}$ and $H\otimes H$ is the canonical element of
the bilinear form on the Cartan subalgebra of $\sw_N$.

If $V$ is an $U_q \asln$--module of type 1, and $\pi: U_q \asln \to
\on{End} V$ is the corresponding action map, then $(\pi \otimes
\on{Id}) (q^{H \otimes H})$ and $(\on{Id} \otimes \pi)(q^{H\otimes
H})$ are well-defined elements of $U_q \asln$ extended by $K_i^{\pm
1/N}$ and $q^{1/N}$.

The algebra $U_q\asln$ is a quasitriangular Hopf algebra. More
precisely, a completed tensor product $$U_q {\mathfrak b}_+
\;\widehat{\otimes}\; U_q {\mathfrak b}_- \subset U_q\asln
\;\widehat{\otimes}\; U_q\asln$$ contains a unique invertible element
$\wt{\mathcal R}$, such that ${\mathcal R} = \wt{\mathcal R}
q^{-H\otimes H}$ satisfies the following identities: 
\bea \Delta^{\on{opp}} =
{\mathcal R} \Delta {\mathcal R}^{-1}, \qquad (\Delta \otimes \on{id})
{\mathcal R} = {\mathcal R}^{13} {\mathcal R}^{23}, \qquad (\on{id}
\otimes \Delta) {\mathcal R} = {\mathcal R}^{13} {\mathcal R}^{12},
\eea 
These identities imply that $\mc R$ satisfies the Yang-Baxter
relation. The element
${\mathcal R}$ is called the universal $R$--matrix of $U_q \asln$. See
\cite{Da,EFK} and Section \ref{isomorphisms section} below for more
details.

\subsection{The algebra $U_q\agln$}\label{agln}

In this section we give two realizations of the quantized enveloping
algebra of $\agln$.

The first realization is obtained by modifying the definition of
$U_q^{\on{new}} \asln$. Namely, we define the algebra
$U_q^{\on{new}}\agln$ with generators $x_{i,n}^{\pm}$
($i=1,\dots,N-1$, $n\in\Z$), $t_i^{\pm 1}$, $g_{i,n}$ ($i=1,\dots,N$,
$n\in \Z\backslash 0$) satisfying the following relations: 
\begin{align*} 
t_it_j
= t_jt_i,\quad t_ig_{j,n}&=g_{j,n}t_i,\quad
g_{i,n}g_{j,m}=g_{j,m}g_{i,n}, \\ t_ix^\pm_{j,n}t_i^{-1}&=
q^{\pm(\delta_{i,j}-\delta_{i-1,j})}x_{j,n}^{\pm},\\ \;[g_{i,n} ,
x_{j,m}^{\pm}]&=0 \qquad (j\neq i,\;j\neq i+1),\\ \; [g_{i,n} ,
x_{i,m}^{\pm}] = \pm \frac{q^{-ni}}{n} [n]_q x_{i,n+m}^{\pm}, &\qquad \;
[g_{i,n} , x_{i-1,m}^{\pm}] = \mp \frac{q^{n(1-i)}}{n} [n]_q
x_{i-1,n+m}^{\pm},\\
\;[x_{i,n}^+,x_{j,m}^-]&=\delta_{ij} \frac{ \phi_{i,n+m}^+ -
\phi_{i,n+m}^-} {q - q^{-1}}, 
\end{align*}
and the last two relations of $U_q^{\rm new}
\asln$. Here we set 
\bean &&\Theta_i^\pm(u):=
\sum_{n=0}^{\infty}\theta_{i,\pm n}^{\pm}u^{\pm n} = t_i^{\pm 1}
\exp\left(\pm(q-q^{-1})\sum_{m=1}^{\infty}g_{i,\pm m} u^{\pm
m}\right),\notag \\ \label{phi theta}&& \Phi_i^\pm(u) : =
\sum_{n=0}^{\infty}\phi_{i,\pm n}^{\pm}u^{\pm n} =
\Theta_i^\pm(uq^{i-1})\Theta_{i+1}^\pm(uq^{i+1})^{-1}.  
\eean

\medskip

Now we define the ``$R$--matrix'' version of algebra $U_q\agln$
following \cite{DF}. This is an associative algebra, denoted by
$U_q^{\on{R}} \agln$ generated by the coefficients $L^\pm_{ij}[k]$ of
the series $L^\pm_{ij}(z)=\sum_{k\geq 0} L^\pm_{ij}[\pm k]z^{\pm
k}$ $(i,j=1,\dots,N)$, satisfying the relations 
\begin{align*}
L^+_{ji}[0]&=L_{ij}^-[0]=0 \qquad\qquad\; (1\leq i<j\leq N),\\
L_{ii}^-[0]L_{ii}^+[0]&=L_{ii}^+[0]L_{ii}^-[0]=1\qquad(1\leq i \leq
N),\\
R(z/w)L_1^{\pm}(z)L_2^{\pm}(w)&=L_2^{\pm}(w)L_1^{\pm}(z)R(z/w),\\
R(z/w)L_1^{+}(z)L_2^{-}(w)&=L_2^{-}(w)L_1^{+}(z)R(z/w),
\end{align*}
where
\bea \lefteqn{R(z)=(q^{-1}z-q)\sum_{i} E_{ii}\ot E_{ii}+}\\
&&+(z-1)\sum_{i\neq j} E_{ii}\ot E_{jj} +(q^{-1}-q)\sum_{i<j}E_{ij}\ot
E_{ji}+z(q^{-1}-q)\sum_{i<j}E_{ji}\ot E_{ij}.  
\eea
Here
$E_{ij}$ is the $N \times N$ matrix with the $(i,j)$ entry equal to
$1$ and all other entries equal to $0$.

Define the Hopf algebra structure on $U_q^{\on{R}} \agln$ by the
formula \be \Delta(L_{ij}^\pm(z))=\sum_{k=1}^N L_{ik}^\pm(z)\ot
L_{kj}^\pm(z).  \ee 

The algebras $U_q^{\on{R}} \agln$ and $U^{\rm new}_q\agln$ are
isomorphic under the isomorphism $\on{DF}: U^{\rm new}_q\agln\to
U_q^{\on{R}} \agln$ described in Section \ref{R-new}.
We will therefore use a single notation $U_q \agln$ for both of them.

\subsection{The algebra $U_q\gln$}\label{gln}
In this section we recall two realizations of quantized enveloping
algebra of $\gln$. 

The algebra $U_q^{\on{DJ}} \gln$ (see \cite{J,DF}) has generators
$X_j^{\pm}$, $T_i^{\pm1}$ ($j=1,\dots,N-1$, $i=1,\dots,N$), and
relations 
\begin{align*} T_iT_j&=T_iT_j,\\
T_iX_j^\pm T_i^{-1}=
q^{\pm(\delta_{i,j}-\delta_{i-1,j})}X_j^\pm,&\qquad [X^+_i,X^{-}_j]=
\delta_{i,j}\frac{T_iT_{i+1}^{-1}-T_{i+1}T_i^{-1}}{q-q^{-1}}, \\
\;[X^{\pm}_i,X^{\pm}_j]&=0\qquad(|i-j|>1),\\
\;[X_i^\pm,[X^\pm_i,X^{\pm}_j]_{q^{-1}}]_{q}&=0 \qquad(|i-j|=1).
\end{align*} 

The comultiplication in $U_q^{\on{DJ}} \gln$ is given by the formula
\begin{align*}
\Delta(T_i) &= T_i \otimes T_i,\notag \\
\Delta(X^+_i) &= X^+_i \otimes 1 + T_iT_{i+1}^{-1}\otimes X^+_i,\\
\Delta(X^-_i)&= X^-_i \otimes T_{i+1}T_i^{-1} + 1 \otimes X^-_i.\notag
\end{align*}

The algebra $U^{\rm DJ}_q\sln$ is the Hopf subalgebra of $U^{\rm
DJ}_q\gln$ generated by $X_i^\pm$, $K_i^{\pm 1}:=(T_i/T_{i+1})^{\pm
1}$ ($i=1,\dots,N-1$).

\medskip

Now we introduce the $R$--matrix realization of the same algebra. 
Let $U_q^{\on{R}}\gln$ be the associative algebra generated by
a$L^\pm_{ij}$, $i,j=1,\dots,N$, satisfying the relations 
\begin{align*}
L^+_{ji}&=L_{ij}^-=0 \qquad\;\;\;\;(1\leq i<j\leq N),\\
L_{ii}^-L_{ii}^+&=L_{ii}^+L_{ii}^-=1\qquad (1\leq i \leq N),\\
RL_1^{\pm}L_2^{\pm}&=L_2^{\pm}L_1^{\pm}R, \qquad
RL_1^{+}L_2^{-}=L_2^{-}L_1^{+}R, 
\end{align*}
where 
\be R=-R(0)=q\sum_{i}
E_{ii}\ot E_{ii}+\sum_{i\neq j} E_{ii}\ot E_{jj}
+(q-q^{-1})\sum_{i<j}E_{ij}\ot E_{ji}. 
\ee

The comultiplication is given by the formula
\be
\Delta(L_{ij}^\pm)=\sum_{k=1}^N L_{ik}^\pm\ot L_{kj}^\pm.
\ee

The algebras $U_q^{\on{DJ}} \gln$ and $U^{\on{R}}_q\gln$ are
isomorphic. The isomorphism $I: U^{\on{DJ}}_q\gln\to
U_q^{\on{R}} \gln$  is described in Section \ref{DJ-R}.
Therefore we will use a single notation $U_q \gln$ for both of them.

\subsection{Evaluation homomorphism} \label{evaluation section}
In this section we define two
evaluation homomorphisms. The relation between them is described in
Section \ref{first cube section}.

Jimbo \cite{J} introduced an evaluation homomorphism
$ev_a:U_q^{\rm DJ}\asln\to U_q^{\rm DJ}\gln$, which
is given on generators by the formulas
\begin{align}\label{ev sln}
ev_a(K_i)&=T_iT_{i+1}^{-1},\qquad ev_a(X_i^\pm)=X_i^\pm\qquad
(i=1,\dots,N-1),\notag\\ 
ev_a(X_0^+)&=
(-q)^N aq^{-1} T_1T_N
[X_{N-1}^-,[X_{N-2}^-,\dots,[X_2^-,X^-_1]_{q^{-1}}\dots
]_{q^{-1}}]_{q^{-1}}, \\ ev_a(X_0^-)&=
a^{-1}q^{-1}T_1^{-1} T_N^{-1}
[X_{N-1}^+,X_{N-2}^+,\dots,[X_{2}^+,X_{1}^+]_{q^{-1}}\dots
]_{q^{-1}}]_{q^{-1}},\notag
\end{align}
(this homomorphism coincides with 
$ev_{(-q)^N}$ in the notation of Section 12.2C in \cite{CP}). 
We remark that in \cite{J} another family of evaluation
homomorphisms was also introduced. However, while the above
homomorphisms $ev_a$ are compatible with our choice of the embedding
$U_q^{\rm DJ} \widehat{\sw}_{N-1} \hookrightarrow U_q^{\rm DJ}\asln$
(namely, ``left upper corner'' embeddings, see Section \ref{first cube
section}), the homomorphisms from the other family are not (they are
compatible with the ``lower right corner'' embeddings).

Also note that we can rewrite the image of $X_0^\pm$ using the braid
group automorphisms \bean\label{ev braid} ev_a(X_0^+)&=&(-q)^N
aq^{-1} T_1T_N {\mc B}_1^{(+,2)} {\mc B}_2^{(+,2)}\dots {\mc
B}_{N-2}^{(+,2)}(X_{N-1}^-),\\
ev_a(X_0^-)&=&a^{-1}q^{-1}T_1^{-1}T_N^{-1}{\mc B}_1^{(-,1)}
{\mc B}_2^{(-,1)}\dots {\mc B}_{N-2}^{(-,1)}(X_{N-1}^+).  \eean

We have a family of Hopf algebra automorphisms $\tau_a: U_q^{\rm DJ}
\asln\to U_q^{\rm DJ} \asln$ depending on a parameter $a\in\C^\times$,
defined by the formulas
$$
\tau_a(X_i^\pm) = X^\pm_i \quad (i\neq 0), \qquad \tau_a(X_0^\pm) =
a^{\pm 1} X_0^\pm, \qquad \tau_a(K_i^{\pm 1}) = K_i^{\pm 1}.
$$
In the new realization they are given by the formulas
\begin{equation}    \label{taua}
\tau_a(x_{i,n}^{\pm})=a^nx_{i,n}^{\pm}, \qquad
\tau_a(h_{i,n}^{\pm})=a^nh_{i,n}^{\pm},\qquad \tau_a(k_i^{\pm1})=k_i^{\pm1}.
\end{equation}
We have:
$$
ev_a \circ \tau_b = ev_{ab}.
$$

\medskip

Next we define the evaluation homomorphism in the $R$--matrix
realization.

\begin{lem}    \label{ev homo}
For any nonzero complex number $a$, there exists a surjective
homomorphism $ev_a:U_q^{\rm R}\agln\to U_q^{\rm R}\gln$ defined on the
generators by the formula \bean\label{gl evaluation}
ev_a(L^\pm_{ij}(z))=\frac{L^\pm_{ij}-(z/a)^{\pm 1}L^\mp_{ij}}{1-(z/a)^{\pm 1}}.
\eean
\end{lem}

\begin{proof}
The verification is straightforward, using the identities \begin{align*}
R(z)&=-R+zR_{21}^{-1},\\ R&=R_{21}^{-1}+(q-q^{-1})P, \end{align*} where
$R_{21}=PRP$ and $P:(\C^n)^{\otimes 2}\to (\C^n)^{\otimes2}$ is the
permutation of the factors.
\end{proof}

We have a family of Hopf algebra automorphisms $\tau_a:U^{\rm
R}_q\agln\to U^{\rm R}_q\agln$ depending on a parameter
$a\in\C^\times$, defined by the formula \bean\label{twist}
\tau_a(L^\pm_{ij}(z))=L^\pm_{ij}(a^{-1}z).\eean Clearly, $ev_a \circ \tau_b =
ev_{ab}$.

In Section \ref{isomorphisms section} we show that the first
evaluation map coincides with the restriction of the second evaluation
to $U_q\asln$.
We remark that the evaluation maps are not homomorphisms of Hopf algebras. 

\subsection{Some obvious maps}\label{inclusions section}

We have the following obvious embeddings of algebras:
\bea
U_q^{\rm new}\widehat{\sw}_{N-1}&\to &U^{\rm new}_q\asln,\qquad
U_q^{\rm new }\widehat{\gw}_{N-1}\to U^{\rm new}_q\agln,\\
U_q^{\rm R}\widehat{\gw}_{N-1}&\to& U^{\rm R}_q\agln,\qquad\;\;
U_q^{\rm R}{\gw}_{N-1}\to U^{\rm R}_q\gln,\\
U_q^{\rm DJ}{\sw}_{N-1}&\to& U^{\rm DJ}_q\sln,\qquad
U_q^{\rm DJ}{\gw}_{N-1}\to U^{\rm DJ}_q\gln,
\eea
given by the ``upper left corner'' rule.
For example, the first embedding maps
\bea
x^{\pm}_{i,n}\mapsto x^{\pm}_{i,n},\qquad
k^{\pm 1}_i\mapsto k_i^{\pm 1},\qquad
h_{i,n}\mapsto h_{i,n},
\eea
where $i=1,\dots,N-2$. 

Note that there also exist ``lower right corner'' embeddings of
algebras, given by similar formulas. However, the ``upper left
corner'' embeddings are compatible with our choice of isomorphisms
between different realizations (see Section \ref{first cube section}),
and the ``lower right corner'' ones are not.

We also have obvious embeddings of Hopf algebras
\bea
U^{\rm DJ}_q\sln &\to& U^{\rm DJ}_q\asln,\qquad
U^{\rm R}_q\gln \to  U^{\rm R}_q\agln,
\eea
where the second map is given by $L^\pm_{ij}\mapsto L^\pm_{ij}[0]$.

Finally, we have a Hopf algebra embedding \be U^{\rm new}_q\asln\to
U^{\rm new}_q\agln, \ee which maps \bea x^{\pm}_{i,n}\mapsto
x^{\pm}_{i,n},\qquad k^{\pm 1}_i\mapsto (t_i t_{i+1}^{-1})^{\pm
1},\qquad h_{i,n}\mapsto q^{n(i-1)}g_{i,n}-q^{n(i+1)}g_{i+1,n} \eea
(see Section \ref{first cube section} for
the compatibility with the comultiplication).

Let $\mc C$ be the algebra generated by the coefficients of the series
$\prod_{i=1}^N\Theta^\pm_i(z)$. All of these coefficients
belong to the center of $U^{\rm new}_q\agln$, and so ${\mc C}$ is a
central subalgebra of $U^{\rm new}_q\agln$.

We have an embedding of algebras \bean\label{det decom} U^{\rm
new}_q\asln\otimes {\mathcal C}\to U^{\rm new}_q\agln. \eean If we
extend the algebra $U^{\rm new}_q\asln$ by $q^{\pm 1/N}$ and $k_i^{\pm
1/N}$ and the algebra $\mc C$ by $q^{\pm 1/N}$ and
$(\prod_{i=1}^Nt_i)^{\pm 1/N}$ (we use the same notation for the
extended algebras here), then we have the inverse embedding
\bean\label{det decom2} U^{\rm new}_q\agln\to U^{\rm
new}_q\asln\otimes {\mathcal C}. \eean So the algebras $U^{\rm
new}_q\agln$ and $U^{\rm new}_q\asln\otimes {\mathcal C}$ are
``almost'' isomorphic.

We will show in Section \ref{isomorphisms section} that the power
series $\prod_{i=1}^N\Theta^\pm_i(z)$ is a group-like element with
respect to the comultiplication, and hence \Ref{det decom} and
\Ref{det decom2} are Hopf algebra embeddings.

\section{Isomorphisms between different
realizations}\label{isomorphisms section} In this section we describe
the isomorphisms between different realizations of $U_q\asln$,
$U_q\agln$ and $U_q\gln$ defined in the previous section.

\subsection{Isomorphism between $U_q^{\rm new}\asln$ and $U_q^{\rm
DJ}\asln$}\label{new-DJ}
The isomorphism \bea\label{Beck} \on{B}: U^{\rm
new}_q\asln\overset{\sim}\to U^{\rm DJ}_q\asln \eea is essentially due
to Beck \cite{B} (see also\cite{KT,LSS}). We will follow the
normalization of \cite{BCP}. 

Let $\delta=\sum_{i=0}^{N-1}\al_i$ be the
imaginary root of $\asln$. For 
a root $\al$ of $\sln$ and $k\in\Z$, let 
$E_{k\delta+\alpha},F_{k\delta+\alpha}\in U^{\rm DJ}_q\asln$ be the
root vectors given in formula (1.7) of \cite{BCP}. 

Note that $E_{\al_i}=X^+_i$, $F_{\al_i}=X^-_i$, and that the braid
group automorphisms $T_i$ used in \cite{BCP} coincide with the
automorphisms ${\mc B}_i^{(-,2)}$ introduced in formula
\eqref{braid1}.

\begin{lem}\label{B lemma} {\em (Lemma 1.5 in \cite{BCP})}. The following
formulas define 
an isomorphism of algebras $\on{B}: U^{\rm
new}_q\asln\overset{\sim}\to U^{\rm DJ}_q\asln$:
\bea
\;\;\on B(x_{i,r}^+)=\left\{\begin{matrix} (-1)^iE_{r\delta+\alpha_i}&(r\geq
0),\\ (-1)^{i+1}F_{-r\delta-\alpha_i}K_i^{-1} &(r<0),
\end{matrix}\right. \quad 
\on B(x_{i,r}^-)=\left\{\begin{matrix} 
(-1)^{i+1}K_iE_{r\delta-\alpha_i}&(r> 0),\\ 
(-1)^iF_{-r\delta+\alpha_i} &(r\leq 0).
\end{matrix}\right.
\eea
\end{lem}
Note that we chose the function $o(i)$ in Lemma 1.5 of \cite{BCP} to be
\be
o(i)=(-1)^i.
\ee

We need the following explicit formula for $x_{1,n}^\pm$ which
follows from Theorem 1 in \cite{BCP}:

\begin{lem}\label{w1}
We have 
\bea
x_{1,n}^\pm&=&(-1)^n{\mc B}^{\mp n}_{\om_1}(X^\pm_1),\\
{\mc B}_{\om_1}&=&\mc B_0^{(-,2)}\mc B_{N-1}^{(-,2)}\mc B_{N-2}^{(-,2)}\dots\mc
B_2^{(-2)}\mc D,
\eea
where $\mc D$ is given by \Ref{mc T}.
\end{lem}

Now we compute the inverse isomorphism.

\begin{lem}\label{inverse B}
The following formulas define
an isomorphism of algebras 

\noindent
$\on{B}^{-1}: U^{\rm
DJ}_q\asln\overset{\sim}\to U^{\rm new}_q\asln$, such that $\on B^{-1}\circ
\on B=\on{Id}$:
\begin{align*}
\on B^{-1} (X_i^\pm)=x_{i,0}^\pm,\qquad \on B^{-1}(K_i^{\pm 1})=k_i^{\pm 1}
\qquad(i=1,\dots,N-1),\\
\on B^{-1}(X_0^+)=(-q)^N[x_{N-1,0}^-
,[x_{N-2,0}^-,\dots,[x_{2,0}^-,x_{1,1}^-]_{q^{-1}}\dots ]_{q^{-1}}]_{q^{-1}}K_0,\\ 
\on B^{-1}(X_0^-)=(-q)^{-N}K_0^{-1}[[\dots[x_{1,-1}^+,x_{2,0}^+]_q,\dots,x_{N-2,0}^+]_q,x_{N-1,0}^+]_q.
\end{align*}
\end{lem}
\begin{proof}
From Lemma \ref{w1} we obtain
\bea
x_{1,1}^-&=&K_1[\dots[[X^+_0,X^+_{N-1}]_{q^{-1}},X_{N-2}^+]_{q^{-1}},\dots,X_2^+]_{q^{-1}},
\\ 
x_{1,-1}^+&=&[X_2^-,[X_3^-,\dots,[X_{N-1}^-,X_0^-]_q\dots]_q]_q K_1^{-1}.
\eea
Now the formulas for $X_0^\pm$ are checked by substituting this for
$x_{1,1}^\pm$ and expanding the commutators of Drinfeld-Jimbo
generators. The formulas for other generators are obvious.
\end{proof}

\subsection{The first fundamental representation of $U_q\asln$}
We describe the explicit action of generators of $U_q\asln$ in the
representation obtained by pullback with respect to $ev_z$ of the
first fundamental representation of $U_q\sln$.

The following lemma is proved by direct verification.
\begin{lem}\label{DJ fund}
For any $z\in\C^\times$, there is an $N$--dimensional vector
representation $V(z)$, of $U^{\rm DJ}_q\asln$ with basis
$\{v_0,\dots,v_{N-1}\}$ 
and action $\pi_z$ given by 
\begin{align*} 
\pi_z(X_i^+)&=E_{i-1,i},\qquad \pi_z(X_i^-)=E_{i,i-1},\\
\pi_z(K_i)=qE_{i-1,i-1}&+q^{-1}E_{i,i}+\sum_{j;\;j\neq i,i-1}E_{j,j}
\qquad (i=1,\dots, N-1),\\
\pi_z(X_0^+)&= (-q)^{N}z
E_{N,1}, \qquad 
\pi_z(X_0^-) = (-q)^{-N}z^{-1} E_{1,N}.
\end{align*}
Moreover, this representation is isomorphic to $ev_z^*(V)$, where
$V$ is the defining $N$--dimensional representation of $\sln$. 
\end{lem}

Now we compute the action of the new generators on $V(z)$.
\begin{lem}\label{new fund}
The action of $U^{\rm new}_q\asln$ on $V(z)$ is given by the formulas
\bea \pi_z(x_{i,k}^+)&=&
(q^iz)^kE_{i-1,i},\\ \pi_z(x_{i,k}^-)&=& (q^iz)^kE_{i,i-1},\\
\pi_z(h_{i,k})&=&\frac{[k]_q}{k}(q^iz)^k(q^{-k}E_{i-1,i-1}-q^kE_{ii}).
\eea 
\end{lem}
\begin{proof}
It is straightforward to check that the formulas above define a
representation of $U^{\rm new}_q\asln$ on $V(z)$ 
(cf. Proposition 3.2.B in \cite{Ko}).
Using Lemma \ref{inverse B}, we obtain then that the action of
$U^{\rm DJ}_q\asln$ is as in Lemma \ref{DJ fund}.
\end{proof}

\subsection{Isomorphism between $U^{\rm new}_q\agln$ and $U^{\rm
R}_q\agln$}\label{R-new}

The isomorphism \bea\label{DF} \on{DF}: U^{\rm new}_q\agln
\overset{\sim}\to U^{\rm R}_q\agln \eea is obtained by J. Ding and
I. Frenkel in \cite{DF}. Here
we describe the inverse map (in our normalization).

Recall the representation $V(z)$ of $U_q^{\rm new}\asln$ 
with action $\pi_z$ given by Lemma \ref{new fund}.
Introduce the $N\times N$ matrices \be \wt{L}^+(z)=(\pi_z\otimes
\on{id}){\mathcal R}^{-1},\qquad
\wt{L}^-(z)=(\on{id}\otimes\pi_z){\mathcal R}, \ee with coefficients
in $U^{\rm new}_q\asln$. 
More precisely, the entries of $\wt{L}^\pm(z)$
are power series in $z^{\pm 1}$, whose coefficients are matrices whose
entries belong to an
extension of $U^{\rm new}_q\asln$ by $k_i^{\pm 1/N}$ and $q^{\pm
1/N}$. The matrices $\wt{L}^\pm(z)$ are called the $L$--operators of
the $N$--dimensional 
representation of $U_q \asln$. The Yang-Baxter relation on the
universal $R$--matrix implies that $L^\pm(z)$ satisfy the
relations of $U^{\rm R}_q\agln$.

Note that in the above formula we used the inverse matrix ${\mathcal
R}^{-1}$ instead of ${\mathcal R}$. This choice allows us to construct
the isomorphism which is compatible with the embeddings $U^{\rm
R}_q{\gw}_{N-1}\to U^{\rm R}_q\agln$ as above (that is using
the ``left upper corner'' rule).

Using the embedding $U^{\rm new}_q\asln\to U^{\rm new}_q\agln$, we
obtain two $N\times N$ matrices $\tilde{L}^\pm(z)$ with entries in
$U^{\rm new}_q\agln$. Now set \be L^\pm(z)=\tilde L^\pm(z)\;
\exp\left(
(q-q^{-1})\sum_{n>0}\frac{q^{-Nn}z^{n}}{[N]_{q^n}}\sum_{i=1}^Ng_{i,-n}
\right)\prod_{i=1}^Nt_i^{1/N}. \ee 
The coefficients of $L^\pm(z)$ a
priori belong to an extension of $U^{\rm new}_q\agln$ by the elements
$t_i^{1/N}$, but according to the formulas below, they actually
belong to $U^{\rm new}_q\agln$.
Note that $L^\pm(z)$ differs from $\wt{L}^\pm(z)$ only by a series, whose
coefficients are central elements of $U_q^{\rm new}\agln$ extended by
$t_i^{1/N}$.  

Define a homomorphism $(\on{DF})^{-1}:\;U^{\rm R}_q\agln\to
U^{\rm new}_q\agln$ sending the generators $L^\pm_{ij}[n]$ of $U^{\rm
R}_q\agln$ to the corresponding entries of the above matrices
$L^\pm(z)$, which are elements of $U^{\rm new}_q\agln$. By
construction, $(\on{DF})^{-1}$ is a well-defined algebra homomorphism.

\begin{lem}
The homomorphism $(\on{DF})^{-1}$ is inverse to the isomorphism from
\cite{DF}.
\end{lem}
\begin{proof}
We use the 
triangular decomposition of the universal $R$-matrix of $U^{\rm
new}_q\asln$ obtained in \cite{Da}. We have $${\mathcal R}={\mathcal
R}^+{\mathcal R}^0 {\mathcal R}^-q^{-H\otimes H},$$ where $H\otimes H$
is defined in \secref{rmatr}. Explicitly we have \bea {\mathcal
R}^+=\prod_{\alpha\in\Delta_+}\prod_{n\geq 0}
\exp_q((q^{-1}-q)E_{\alpha+n\delta}\otimes F_{\alpha+n\delta}),\\
{\mathcal R}^-=\prod_{\alpha\in\Delta_+}\prod_{n<0}
\exp_q((q^{-1}-q)E_{\alpha+n\delta}\otimes F_{\alpha+n\delta}), \eea
where $\Delta_+$ is the set of positive roots of ${\mathfrak
{sl}}_N$, and \be \exp_q(x)=\sum_{m\geq 0}\frac{x^m}{(m)_{q^2}!},\qquad
(m)_{q}!=\prod_{i=1}^m\frac{1-q^m}{1-q} \ee (see Definition 4 in
\cite{Da} for the description of the order of these factors). The explicit
formula for ${\mathcal R}^0$ is given in Theorem 2 of \cite{Da} (see
also formula (3.7) in \cite{FR}).

Therefore we find (compare with formula (3.22) of \cite{DF}): 
\begin{align*}
&L^+(z)=
\left(\begin{matrix}
1 & 0 & 0&\dots &0&0\\
x_1^+(zq)_{< 0}&1 & 0 &\dots & 0&0\\
*&x_2^+(zq^2)_{< 0}&1 & \dots & 0& 0\\
\dots & \dots & \dots & \dots &\dots & \dots \\
* & * & * & \dots  & 1 & 0 \\
\tilde X_0^-z+\dots & * & * & \dots & x_{N-1}^+(zq^{N-1})_{< 0} & 1
\end{matrix}\right)\times\\
&\left(\begin{matrix}
\Lambda_1^-(z)& 0 & \dots & 0\\
0 & \Lambda_2^-(z)& \dots & 0\\
\dots & \dots & \dots & \dots\\
0 & 0 & \dots &  \Lambda_N^-(z)
\end{matrix}\right)
\left(
\begin{matrix}
1 & x_1^-(zq)_{\leq 0} &*&\dots &*\\
0 & 1 &x_2^-(zq^2)_{\leq 0} &\dots& *\\
0 & 0 & 1 & \dots & *\\
\dots & \dots & \dots & \dots & x_{N-1}^+(zq^{N-1})_{< 0}  \\
0 & 0 & 0 & \dots & 1
\end{matrix}\right)\notag
\end{align*}
and
\begin{align*}
&L^-(z)=
\left(\begin{matrix}
1 & 0 & 0&\dots &0&0\\
-x_1^+(zq)_{\geq 0}&1 & 0 &\dots & 0& 0\\
*&-x_2^+(zq^2)_{\geq 0}&1 & \dots & 0&0\\
\dots & \dots & \dots & \dots & \dots &\dots \\
* & * & * & \dots &1 &0 \\
* & * & * & \dots &-x_{N-1}^+(zq^{N-1})_{\geq 0}  &1
\end{matrix}\right)\times\\
&\left(\begin{matrix}
\Lambda_1^{+}(z)& 0 & \dots & 0\\
0 & \Lambda_2^{+}(z) & \dots & 0\\
\dots & \dots & \dots & \dots\\
0 & 0 & \dots & \Lambda_N^{+}(z)
\end{matrix}\right)
\left(\begin{matrix}
1 & -x_1^-(zq)_{>0} &*&\dots & -\tilde X_0^+z^{-1}+\dots\\
0 & 1 &-x_2^-(zq^2)_{> 0} &\dots& *\\
0 & 0 & 1 & \dots & *\\
\dots & \dots & \dots & \dots & -x_{N-1}^-(zq^{N-1})_{>0} \\
0 & 0 & 0 & \dots & 1
\end{matrix}\right)
\end{align*}
where 
\begin{align*}
x_i^\pm(z)_{\leq 0}=(q-q^{-1})&\sum_{j\leq 0} x^\pm_{i,j}z^{-j},\qquad
(x_i^\pm)(z)_{< 0}=(q-q^{-1})\sum_{j<0} x^\pm_{i,j}z^{-j},\\
x_i^\pm(z)_{\geq 0}=(q-q^{-1})&\sum_{j\geq 0} x^\pm_{i,j}z^{-j},\qquad
(x_i^\pm)(z)_{> 0}=(q-q^{-1})\sum_{j>0} x^\pm_{i,j}z^{-j},\\
&\tilde X_0^-z+\dots = z(-q)^N(q-q^{-1})t_NX_0^-t_1^{-1}+z^2(\ldots),\\
-&\tilde X_0^+z^{-1}+\dots=
-z^{-1}(-q)^{-N}(q-q^{-1})t_1X_0^+t_N^{-1}+z^{-2}(\ldots),\\
&\Lambda_i^\pm(z)=\frac{1}{\Theta^\pm_i(z^{-1}q)}\prod_{j=1}^{i-1}
\frac{\Theta^\pm_j(z^{-1}q^{-1})}{\Theta^\pm_j(z^{-1}q)}, \\
\end{align*}
(note that we have moved the term $q^{\pm H\otimes H}$ to the middle
factor).
\end{proof}

Introduce the series
$$
S^\pm(z) = \sum_{\sigma \in S_N} (-q)^{-l(\sigma)}
L^\mp_{1\sigma(1)}(z) L^\mp_{2\sigma(2)}(zq^{-2}) \ldots
L^\mp_{N\sigma(N)}(zq^{2-2N}).
$$
According to \cite{J}, the coefficients of these series are central
elements of the algebra $U_q^{\rm R} \agln$. Moreover, because
$S^\pm(z)$ may be interpreted as the $L$--operators corresponding to a
one-dimensional representation of $U_q^{\rm R} \agln$ (called the
quantum determinant), they have ``group--like''
comultiplication:
\begin{equation}    \label{grcom}
\Delta(S^\pm(z)) = S^\pm(z) \otimes S^\pm(z).
\end{equation}

The coefficients of the series \bean \label{La and Theta}
\prod_{i=1}^N\La^\pm_i(zq^{2-2i}) =
\prod_{i=1}^N\Theta_i^\pm(z^{-1}q^{2N-1})^{-1} \eean are
central elements of $U_q^{\on{new}} \agln$. 

The above formulas for $(\on{DF})^{-1}$ imply that $S^\pm(z)$ and
$\prod_{i=1}^N\La^\pm_i(zq^{2-2i})$ act by the same power series in
$z^{\pm 1}$ on any lowest weight vector. Therefore we obtain
\begin{lem} We have
\be (\on{DF})^{-1}(S^\pm(z)) = \prod_{i=1}^N \Theta_i^\pm(z^{-1}q^{2N-1})^{-1}.
\ee In particular, the maps \eqref{det decom} and \Ref{det decom2} are
embeddings of Hopf algebras.
\end{lem}

\subsection{Embedding of $U_q^{\rm DJ}\asln$ into $U_q^{\rm
R}\agln$}\label{DJ-R}
Now we can complete the picture and compute explicitly the embedding
of $U_q^{\rm DJ}\asln$ into $U_q^{\rm R}\agln$ and the isomorphism
between $U_q^{\rm DJ}\gln$ and $U_q^{\rm R}\gln$.

\begin{lem}\label{agl dj to r}
The following formulas 
define an embedding of Hopf algebras 

\noindent $\hat I: U_q^{\rm DJ}\asln\to U_q^{\rm R}\agln$,
\bea
\hat I(X^+_0)&=& (-q)^N ({q^{-1}-q})^{-1} L_{1N}^-[-1]L_{NN}^+[0],\\
\hat I(X^-_0)&=& (-q)^{-N} ({q-q^{-1}})^{-1} L_{NN}^-[0]L_{N1}^+[1] ,\\
\hat I(X_i^+)&=& ({q^{-1}-q})^{-1} L_{i+1,i}^-[0]L_{ii}^+[0],\\
\hat I(X_i^-)&=& (q-q^{-1})^{-1} L_{ii}^-[0]L_{i,i+1}^+[0] \\
\hat I(K_i)&=&L_{ii}^+[0]L_{i+1,i+1}^-[0] \qquad (i=1,\dots,N-1).
\eea
\end{lem}
\begin{proof}
The formulas for $K_i$, $i=1,\dots,N-1$, follow from comparing the
constant coefficients of the diagonal entries in the above formula for
$L^+(z)$. The formulas for $X^{+}_i$ and $X^{-}_i$, $i=1,\dots,N-1$,
follow from comparing the constant coefficients of the $(i,i+1)$
entries in $L^-(z)$ and of the $(i+1,i)$ entries in $L^+(z)$
respectively. Finally the formulas for $X^{+}_0$ and $X^{-}_0$ follow
from comparing the $z^{-1}$ coefficient of the $(1,N)$ entry of
$L^-(z)$ and the $z$ coefficient of the $(N,1)$ entry of $L^+(z)$
respectively.

After that it is straightforward to check on the generators that this
map is compatible with comultiplication.
\end{proof}

Similarly, we obtain the following

\begin{lem}\label{gl dj to r}
The following formulas define an isomorphism of Hopf algebras

\noindent
$I:\;U^{\rm DJ}_q\gln\overset{\sim}\to U^{\rm R}_q\gln$,
\bea
I(X_i^+)&=& ({q^{-1}-q})^{-1} L_{i+1,i}^- L_{ii}^+ \qquad(i=1,\dots,
N-1),\\
I(X_i^-)&=& (q-q^{-1})^{-1} L_{ii}^-  L_{i,i+1}^+\qquad(i=1,\dots,
N-1),\\
I(T_i)&=&L_{ii}^+ \qquad(i=1,\dots, N).
\eea
\end{lem}

\subsection{Summary}\label{first cube section}
The maps described above are summarized in Figure \ref{figure 1}.
\setlength{\unitlength}{1mm}
\begin{figure}
\begin{picture}(40,112)(50,-60)    \label{pic1}

\put(49,46){\vector(-3,-2){21}} 
\put(113,49){\vector(-1,0){50}}
\put(86,27){\vector(-1,0){50}}
\put(21,26){$U_q^{\rm new}\widehat{\gw}_N$}
\put(48,48){$U_q^{\rm new}\asln$}
\put(115,48){$U^{\rm new}_q\widehat{\sw}_{N-1}$}
\put(88,26){$U^{\rm new}_q\widehat{\gw}_{N-1}$}
\put(116,46){\vector(-3,-2){21}} 

\put(25,23){\vector(0,-1){32}} 
\put(119,45){\vector(0,-1){32}}
\put(52,24){\vector(0,-1){11}}
\put(52,44){\line(0,-1){15}}
\put(92,23){\vector(0,-1){32}}

\put(14,5){$\sim \on{DF}$}
\put(44,19){$\sim \on B$}
\put(81,3){$\sim \on{DF}$}
\put(120,28){$\sim \on B$}

\put(49,5){\vector(-3,-2){21}}
\put(86,-13){\vector(-1,0){50}}
\put(21,-14){$U_q^{\rm R}\agln$}
\put(48,8){$U_q^{\rm DJ}\asln$}
\put(115,8){$U^{\rm DJ}_q\widehat{\sw}_{N-1}$}
\put(88,-14){$U_q^{\rm R}\widehat{\gw}_{N-1}$}
\put(116,5){\vector(-3,-2){21}} 

\put(39,-6){$\hat I$} 
\put(105,-6){${\hat I}$}

\put(25,-17){\vector(0,-1){32}} 
\put(119,5){\vector(0,-1){32}}
\put(52,5){\line(0,-1){16}}
\put(52,-15){\vector(0,-1){12}}
\put(92,-17){\vector(0,-1){32}}

\put(19,-35){$ev_a$}
\put(46,-21){$ev_a$}
\put(86,-40){$ev_a$}
\put(121,-13){$ev_a$}

\put(49,-35){\vector(-3,-2){21}}
\put(89,-31){\vector(-1,0){27}}
\put(95,-31){\line(1,0){18}}
\put(86,-53){\vector(-1,0){50}}
\put(21,-54){$U_q^{\rm R}\gln$}
\put(48,-32){$U_q^{\rm DJ}\gln$}
\put(114,-32){$U_q^{\rm DJ}{\gw}_{N-1}$}
\put(88,-54){$U_q^{\rm R}{\gw}_{N-1}$}
\put(116,-35){\vector(-3,-2){21}} 

\put(39,-46){$\sim I$} 
\put(105,-46){$\sim I$}
\end{picture}
\caption{Different realizations of quantum groups.}\label{figure 1}
\end{figure}
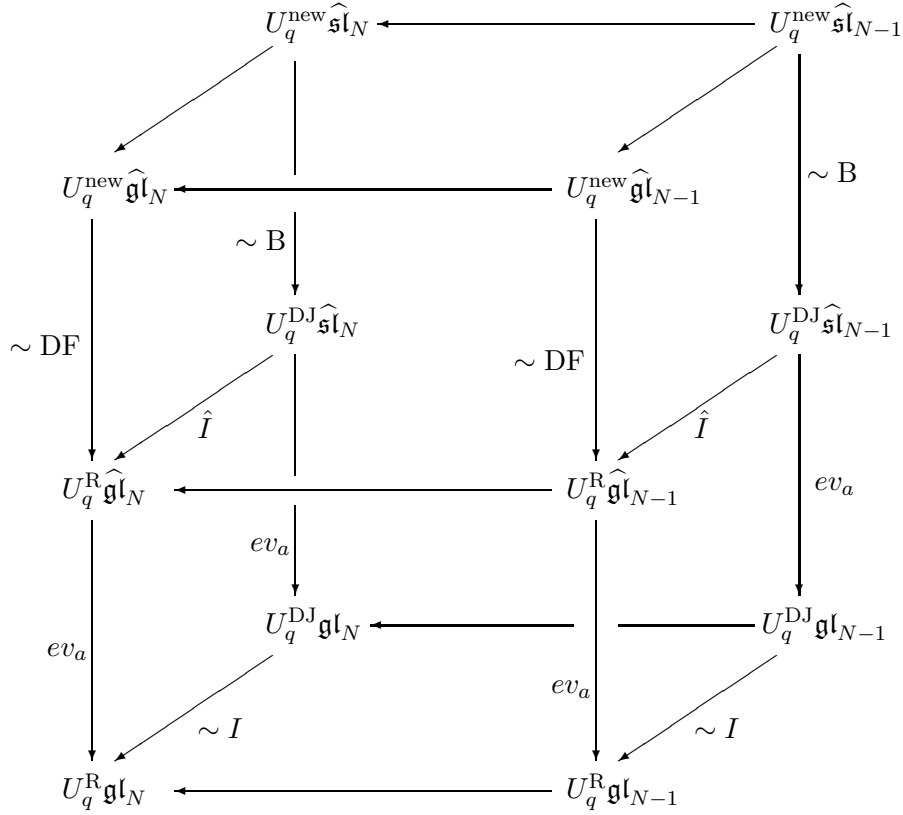
All maps in Figure \ref{figure 1} 
are algebra homomorphisms, moreover the maps in the upper
left and the upper right faces are Hopf algebra homomorphisms.

\begin{thm}\label{cube theorem}
The diagram in Figure \ref{figure 1} is commutative.
\end{thm}
\begin{proof}
The commutativity of all faces except for the back, lower left and
lower right ones
follow from the definitions of the maps.

The commutativity of the lower left face is checked explicitly on
generators. The only nontrivial case is that of $X^\pm_0$.  We use the
relations \bea (RL^+_1L^+_2)_{ik}^{kj}=(L_2^+L_1^+R)_{ik}^{kj},\qquad
(1\leq i<k<j\leq N),\\
(RL^-_1L^-_2)_{ik}^{kj}=(L_2^-L_1^-R)_{ik}^{kj},\qquad (N\leq
i>k>j>0), \eea which hold in $U_q^{\rm R}\gln$. They can be
rewritten in the form \bea \;[ L_{ik}^+,L_{kj}^+ ] &=&
(q^{-1}-q)L_{kk}^+ L_{ij}^+ \qquad (1\leq i<k<j\leq N),\\
\;[L_{ik}^-,L_{kj}^-] &=&(q-q^{-1})L_{kk}^-L_{ij}^- \qquad (N\geq
i>k>j>0).  \eea Now, we compute
\begin{align*}
&I^{-1}(ev_a(\hat I(X^+_0))) \\
&=I^{-1}\left(ev_a
\left(\frac{(-q)^{N+1}}{q-q^{-1}} L_{1N}^-[-1]L_{NN}^+[0]\right)\right)=
I^{-1}\left(\frac{a(-q)^{N}}{q-q^{-1}}L_{1N}^+L_{NN}^+\right) \\
&
=I^{-1}\left(\frac{a(-q)^{N}}{(q-q^{-1})^{N-1}}
L_{N-1,N-1}^-[L_{N-1,N}^+,\dots,L_{33}^-
[L^+_{34},L_{22}^-[L_{23}^+,L_{12}^+]]\dots] L_{NN}^+ \right)\\
&=a(-q)^{N}T_{N-1}^{-1}[T_{N-1}X^-_{N-1},\dots,T_3^{-1}[T_3X_3^-,T_2^{-1}
[T_2X_2^-,T_1X_1^-]]\dots] T_N\\
&\qquad \qquad \quad 
=(-1)^N aq^{N-1}
[X_{N-1},\dots,[X_3^-,[X_2^-,X_1^-]_{q^{-1}}]_{q^{-1}}\dots]_{q^{-1}}
T_1T_N=ev_a(X_0^+).
\end{align*}
The computation for $X_0^-$ is similar. The same argument also proves
the commutativity of the lower right face.

The commutativity of the back face is proved by standard diagram
chasing argument from the commutativity of the other faces.
\end{proof}

\section{Finite-dimensional representations and $q$--characters for
generic $q$}\label{rep generic section}

In this section we recall known facts about finite-dimensional
representations of $U_q\asln$ in the case of generic $q$ and extend
them to the case of $U_q\agln$. We also compute the $q$--characters of
the evaluation modules of $U_q\agln$.

\subsection{Finite-dimensional representations of $U_q \asln$}
\label{fd}

First we recall the classification of finite-dimensional
representations of $U_q\asln$, due to Chari and Pressley. The details
can be found in \cite{CP} and \cite{CP:gen}.

Recall that in this paper we restrict ourselves to type 1
representations of $U_q \asln$ (see the definition in \secref{rmatr}).

Define the elements ${\mc P}_{i,n}\in U^{\rm new}_q\asln$
($i=1,\ldots,N-1$, $n \in \Z$) through the generating
functions \be {\mathcal P}_i^\pm(u)=\sum_{n\geq 0} {\mc P}_{i,\pm n}
u^{\pm n}= \exp\left(\mp \sum_{m>0}\frac{h_{i,\pm m}}{[m]_q}u^{\pm
m})\right).  \ee We have \bean\label{phi p} \Phi_i^\pm(u)=k_i^{\pm
1}{\mathcal P}_i^\pm(uq^{-1}){\mathcal P}_i^\pm(uq)^{-1}.  \eean

Let $V$ be a representation of $U_q\asln$. 
A vector $v\in V$ is called a highest weight vector if
\begin{equation}    \label{hwv}
x_{i,n}^+ \cdot v=0,\quad \quad {\mc P}_{i,n}\cdot v = P_{i,n}v, \quad
\quad k_i \cdot v = q^{\mu_i} v
\end{equation}
($i=1,\dots,N-1$, $n\in\Z$),
where $P_{i,n} \in \C$ and $\mu_i
\in \Z$ (actually, then necessarily $\mu_i \in \Z_{\geq 0}$). The
representation $V$ is called a highest weight representation if
$V=U_q\asln\cdot v$, for some highest weight vector $v$.
Every finite-dimensional irreducible representation $V$ of $U_q\asln$
of type 1 is a highest weight representation and 
\bean\label{hw}
P_i^\pm(u):= \sum_{n\geq 0} P_{i,\pm n} u^{\pm
n}=\prod_{j=1}^{\mu_i}(1-a_j^{\pm 1}u^{\pm 1}). \eean
In particular, the polynomials $P_i^+(u)$ determine the polynomials
$P_i^-(u)$ and the numbers $\mu_i$. In order to simplify notation, we
will denote $P_i^+(u)$ by $P_i(u)$. The polynomials $P_i(u)$ 
($i=1,\ldots,N-1$), are called the Drinfeld polynomials of $V$. For any
$(N-1)$--tuple ${\bs P} = (P_1(u),\ldots,P_{N-1}(u))$ of polynomials
with constant term 1, there exists a unique representation of $U_q
\asln$ with such Drinfeld polynomials.

For $a\in\C^\times$, the Drinfeld polynomials of the pullback
of irreducible $U_q\asln$--module with
Drinfeld polynomials $P_i(u)$ by the automorphism $\tau_a$
(see formula \Ref{taua}) are given by $P_i(au)$ ($i=1,\dots,N-1$).

\subsection{Finite-dimensional representations of
$U_q\agln$}\label{agln rep}

A vector $w$ in a $U_q\gln$--module $W$ is called a vector of weight
$\bs\la=(\la_1,\dots,\la_N) \in \Z^n$, if $$T_i \cdot w = q^{\la_i} w
\qquad (i=1,\dots,N).$$ A representation $W$ of $U_q\gln$ is said to be
of type 1 if it is the direct sum of its weight spaces $W =
\oplus_{\bs\la} W_{\bs\la}$, where $W_{\bs\la} = \{w\in W | T_i \cdot
w = q^{\la_i}w\}$.

Let ${\mc F}_N$ be the category of finite-dimensional representations
$W$ of $U_q\gln$ of type 1, such that all weights ${\bs \la} =
(\la_1,\dots,\la_N)$ occurring in $W$ have nonnegative integral
components, i.e., $\la_i\in\Z_{\geq 0}$ ($i=1,\dots,N$). These
representations are quantum analogues of the polynomial
representations of $GL_n$. The category ${\mc F}_N$ is semisimple,
with simple objects being the irreducible representations with highest
weights $\bs\la=(\la_1,\dots,\la_N)$, where
$\la_1\geq\la_2\geq\dots\geq \la_N\geq 0$.

A representation $V$ of $U_q\agln$ is said to be of type 1 if $V$ is
of type 1 as a representation of $U_q\gln$.

Define elements ${\mc Q}_{i,n}\in U_q^{\rm new}\agln$ ($i=1,\dots,N$,
$n\in\Z$), through the 
generating functions \be {\mathcal Q}_i^\pm(u)=\sum_{n\geq 0} {\mc
Q}_{i,\pm n} u^{\pm n}= \exp \left(\mp \sum_{m>0} \frac{g_{i,\pm
m}}{[m]_q}u^{\pm m})\right). \ee We have 
\bean\label{theta q}
\Theta_i^\pm(u)&=&t_i^{\pm 1}{\mathcal
Q}_i^\pm(uq^{-1}){\mathcal Q}_i^\pm(uq)^{-1}, \\
{\mc P}^\pm_i(u)&=&{\mc Q}^\pm_i(uq^{i-1}) {\mc
Q}^\pm_{i+1}(uq^{i+1})^{-1}.  \label{PandQ}
\eean
Note that the coefficients $C_{n}$ ($n \in \Z$), of the
series
\begin{equation}    \label{Cn}
C^\pm(u): = \sum_{n\geq 0} C_{\pm n} u^{\pm n} = \prod_{i=1}^N {\mc
Q}^\pm_i(u)
\end{equation}
are central elements of $U_q\agln$. Recall that $\mc C$ is the algebra
generated by the coefficients of the products
$\prod_{i=1}^N\Theta_i^\pm(u)$. Then $\mc C$ is also generated by $C_n$
($n\in\Z$) and the element $\prod_{i=1}^N t_i$.

Let $V$ be a finite-dimensional representation of $U_q\agln$ of type
1.  Recall that the elements $\{ t_i,{\mc Q}_{i,n} \}$ commute among
themselves, so every finite-dimensional representation $V$ is a direct
sum of generalized eigenspaces \be V_{(\bs\la,\bs\ga)} = \{ x \in
V_{\bs\la} |\; ({\mc Q}_{i,n} - \ga_{i,n})^p \cdot x = 0\;\; \on{for}\;
\on{some}\; p, \;\; (i=1,\dots,N,\;n \in \Z )\}. \ee The
decomposition of $V$ into a direct sum of subspaces
$V_{(\bs\la,\bs\ga)}$ is a refinement of its weight decomposition.

Given a collection $(\bs \la,\bs\ga)$ of generalized eigenvalues,
we form the generating functions
$$
\Gamma^\pm_i(u) = \sum_{n\geq 0} \ga_{i,\pm n} u^{\pm n}.
$$
We refer to the collections 
$(\bs\la,\bs\Gamma):=(\la_i,\Gamma^\pm_i(u))_{i=1,\dots,N}$
occurring on a given representation $V$ as the common spectra of $\{
t_i,{\mathcal Q}^\pm_i(u) \}$.

\begin{lem}    \label{mult prop}
The generalized eigenvalues of $\{t_i, {\mc Q}^\pm_i(u)\}$ satisfy the
multiplicative property, i.e., the generalized eigenvalues on $V_1
\otimes V_2$ are the products of the generalized eigenvalues on
$V_1$ and $V_2$.
\end{lem}

\begin{proof}
The comultiplication formula $\Delta(t_i) = t_i \otimes t_i$ implies
that the eigenvalues of $t_i$'s satisfy the multiplicative property.

Let us extend $U_q \agln$ by adding $t_i^{\pm 1/N}$ and
$q^{1/N}$. Consider the commutative subalgebra ${\mc J}$ (resp.,
$\wt{\mc J}$) of $U_q \agln$ generated by $k_i$, ${\mc P}_{i,n}$, $n \in
\Z$, $i=1,\ldots,N-1$ (resp., $t_i$, ${\mc Q}_{i,n}$, $n \in
\Z$, $i = 1,\ldots,N$). Then we have an embedding $\wt{\mc
J} \hookrightarrow {\mc J} \otimes {\mc C}$.

Using this embedding, it suffices to establish the
multiplicative property of ${\mc P}^\pm_i(u)$ and $C^\pm(u)$.  The
former follows from Proposition 1 in \cite{FR}, and the latter follows
from the comultiplication formula \eqref{grcom}.
\end{proof}

The maps \Ref{det decom}, \Ref{det decom2} allow us to identify 
irreducible $U_q \agln$--modules with tensor products of the form
$V' \otimes V''$, where $V'$ is an irreducible representation of
$U_q\asln$, on which ${\mc C}$ acts trivially, and $V''$ is a
one-dimensional representation, on which $U_q \asln$ acts trivially.
That enables us to carry over the
results about representations of $U_q \asln$ to the case of $U_q
\agln$. 

Let $V$ be a representation of $U_q\agln$. A vector $v\in V$ is called
a highest weight vector, if
\bean   \label{qi}
x_{j,n}^+ \cdot v=0,\qquad {\mc Q}^\pm_i(u) \cdot v =
Q^\pm_i(u) v \qquad  t_i \cdot v = q^{\la_i} v 
\eean
($i=1,\dots,N,\;j=1,\dots,N-1,\;n\in\Z$).
The module $V$ is called a highest weight module if $V=U_q\agln\cdot v$
for some highest weight vector $v$.

The following lemma follows from existence of maps 
\Ref{det decom}, \Ref{det decom2} and the
description of $U_q\asln$--modules.
\begin{lem}
Any irreducible $U_q
\agln$--module $V$, whose restriction to $U_q \gln$ is in ${\mc F}_N$,
is a highest weight module and the ratios 
of the eigenvalues of $\mc Q_i^\pm(u)$ corresponding to the highest
weight vector 
\be
Q_i^\pm(uq^{i-1}) Q^\pm_{i+1}(uq^{i+1})^{-1},
\ee 
are of the form \Ref{hw}, where $\mu_i =
\la_i-\la_{i+1}$. Moreover, there is a bijection between 
$N$-tuples $Q_i^\pm(u)$
$(i=1,\dots,N)$ with such property and the irreducible
$U_q\agln$--modules whose restriction to $U_q \gln$ is in ${\mc F}_N$.
\end{lem}

\subsection{Polynomial representations of $U_q\agln$}\label{polynomial rep}

In order to avoid the redundancy related to the central characters, we
impose a condition on the possible spectra of $\{t_i,{\mc
Q}_i^\pm(u)\}$ on the finite-dimensional representations of $U_q
\agln$.

Let $\wh{\mathcal F}_N$ be the category of all finite-dimensional
representations $V$ of $U_q\agln$, such that

\begin{itemize}
\item the restriction of $V$ to $U_q\gln$ is in ${\mathcal F}_N$;

\item each generalized common eigenvalue of ${\mathcal Q}^\pm_i(u)$
in $V_{\bs\la}$, where ${\bs \la} = (\la_1,\ldots,\la_N)$, is
a polynomial $\Gamma^\pm_i(u)$ in $u^{\pm 1}$ of degree $\la_i$, so that
the zeroes of the functions $\Gamma^+_i(u)$ and $\Gamma^-_i(u)$
coincide.
\end{itemize}
The last condition means that common spectra
of $\{ t_i,{\mathcal
Q}^\pm_i(u) \}$ on a polynomial representation $V$ of $U_q \agln$
are of the form $(\bs\la,\bs\Gamma)$, where
\begin{equation}    \label{formula for common spectra}
\Gamma_{i}^{\pm}(u)=\prod_{k=1}^{\la_i}(1-a_{ik}^{\pm 1} u^{\pm
1}), \qquad a_{ik}\in\C^\times.
\end{equation}
In particular, the
(generalized) eigenvalues of $t_i$ and ${\mc Q}^-_i(u)$ are determined
by those of ${\mc Q}^+_i(u)$.

We will call the objects of the category $\wh{\mathcal F}_N$ {\em
polynomial representations of $U_q\agln$}. In this paper we 
restrict ourselves to polynomial representations of $U_q\agln$.

In particular, if $V$ is an irreducible $U_q \agln$ from the category
$\wh{\mc F}_N$, then it is generated by a highest weight vector
satisfying \eqref{qi}, where $Q^\pm_i(u)$ are polynomials in
$u^{\pm1}$. The polynomials $Q^\pm_i(u)$ have the
property that for $i=1,\dots, N-1$, the ratio
$Q_i^\pm(u)/Q^\pm_{i+1}(uq^2)$ is also a polynomial in $u^{\pm1}$.
We denote $Q^+_i(u)$ by $Q_i(u)$ and by analogy with the $U_q \asln$
case, we call the 
polynomials $Q_i(u)$ ($i=1,\ldots,N$), the Drinfeld polynomials of $V$.

An $N$--tuple of polynomials ${\bs Q}=(Q_1(u),\dots,Q_N(u))$ with
constant terms $1$ is called {\em dominant} if for $i=1,\dots, N-1$,
the ratio $Q_i(u)/Q_{i+1}(uq^2)$ is a polynomial.

Thus the irreducible objects in the category $\wh{\mc F}_N$ are in
one-to-one correspondence with the dominant $N$--tuples of
polynomials. Denote the
irreducible polynomial $U_q \agln$--module
corresponding to a dominant $N$--tuple ${\bs Q}$ by $V({\bs Q})$.

\lemref{mult prop} implies the multiplicative property of the Drinfeld
polynomials of irreducible representations from the category $\wh{\mc
F}_N$. The pullback by the automorphism $\tau_a$ (see formula
\Ref{twist}) has the effect of changing the Drinfeld polynomials by
the rule $Q_i(u) \mapsto Q_i(ua)$.

The Hopf algebra structure on $U_q \agln$ (defined in terms of its
$R$--matrix realization) makes the category of its finite-dimensional
representations into a monoidal category. The following lemma is a
consequence of \lemref{mult prop}.

\begin{lem}    \label{full mon}
The category $\wh{\mathcal F}_N$ is a full monoidal subcategory of the
category of all finite-dimensional representations of $U_q\agln$,
i.e., if $V_1,V_2\in \wh{\mathcal F}_N$, then $V_1\otimes
V_2\in\wh{\mathcal F}_N$.
\end{lem}

The following lemma shows that the definition of $\wh{\mathcal F}_N$ is
not too restrictive.

\begin{lem}\label{gl to sl functor}
Any finite-dimensional representation of $U_q\asln$ may be obtained as
the restriction of a polynomial representation of $U_q\agln$.
\end{lem}

\begin{proof}
Let $V$ be an irreducible finite-dimensional $U_q\asln$--module. In
order to lift it to a representation of $U_q\agln$, we need to define
an action of ${\mc Q}^\pm_i(u)$ and $t_i$ on $V$ in such a way that
the relations \eqref{PandQ} and the relations $k_i = t_i t_i^{-1}$ are
satisfied. We make ${\mc Q}_N^\pm(u)$ act as the identity and define
the action of the remaining series ${\mc Q}^\pm_i(u), i<N$,
recursively using formula \eqref{PandQ}. Likewise, we make $t_N$ act
as the identity and define the remaining elements $t_i$ recursively
using the formula $t_i = t_{i+1} k_i$. The representation of $U_q
\agln$ obtained in this way satisfies the following conditions: its
restriction to $U_q \gln$ is of type 1, so that $V = \oplus V_{\bs
\la}, \bs \la \in \Z^N$; the common generalized eigenvalues of ${\mc
Q}^\pm_i(u)$ on $V_{\bs \la}$ have the form
$$
\Gamma_{i}^{\pm}(u)=\frac{\prod_{k=1}^{k_i}(1-a_{ik}^{\pm 1} u^{\pm
1})}{\prod_{r=1}^{r_i}(1-b_{ir}^{\pm 1} u^{\pm 1})},
$$
and $\la_i = k_i-r_i$.

Now observe that $U_q \agln$ has a family of one--dimensional
representation $\on{Det}_a, a \in \C^\times$ (the quantum
determinants) on which the subalgebra $U_q\asln\subset U_q\agln$ acts
trivially, ${\mc Q}^\pm_i(u)$ acts by multiplication by $(1-a^{\pm 1}
q^{\pm 2(i-1)} u^{\pm 1})$, and $t_i$ acts by multiplication by $q$.

By \lemref{mult prop}, the generalized eigenvalues of ${\mc
Q}^\pm_i(u)$ and $t_i$ satisfy the multiplicative property. This
property implies that the tensor product of $V$ with the
one-dimensional representations $\on{Det}_a$, where $a^{-1}$ runs over
the list of the poles of $\Gamma_{i}^+(uq^{2(i-1)})$ ($i=1,\dots,N$),
has the structure of a $U_q \agln$--module on $V$, which satisfies the
condition of the lemma.
\end{proof}

\subsection{The $q$--characters}\label{qchar section}

In this section we recall some facts about the $q$--characters of
$U_q\asln$--modules which can be found in \cite{FR}, \cite{FM1} and
extend the construction to the case of $U_q\agln$--modules.

Let $V$ be a finite-dimensional $U_q\asln$--module and
${v_1,\dots,v_n}$ be a basis of common generalized eigenvectors of
$\{k_i, {\mathcal P}_i^\pm(u)\}$. By Proposition 1 in
\cite{FR}, the common eigenvalues of operators ${\mathcal P}_i^+(u)$
on each vector $v_s$ are rational functions of the form 
\be
\Gamma_{i,s}^\pm(u)=\frac{\prod_{k=1}^{k_{is}} (1-a_{isk}^{\pm
1}u^{\pm 1})}{\prod_{r=1}^{r_{is}} (1-b_{isr}^{\pm 1}u^{\pm 1})},
\qquad a_{isk},b_{isr}\in\C^\times.
\ee
Moreover, 
$k_iv_s=q^{k_{is}-r_{is}}v_s$.  By definition, the
$q$--character of $V$, denoted by $\chi_q(V)$, is an element of the
commutative ring of polynomials $\Z[Y^{\pm
1}_{i,a}]_{i=1,\dots,N-1}^{a\in\C^\times}$, which is equal to \be
\chi_q(V)=\sum_{s=1}^n \prod_{i=1}^{N-1}\left(\prod_{k=1}^{k_{is}}
Y_{i,a_{isk}}\prod_{r=1}^{r_{is}}
Y_{i,b_{isr}}^{-1}\right).  \ee In the case of an irreducible
representation, the monomial in $\chi_q(V)$ corresponding to the
highest weight vector is called the highest weight monomial.

Let $\Rep \,U_q\asln$ be the Grothendieck ring of the category of
all finite--dimensional representations of $U_q\asln$. Then the map
\be \chi_q:\; \Rep\,U_q\asln\to \Z[Y^{\pm
1}_{i,a}]_{i=1,\dots,N-1}^{a\in\C^\times}, \qquad [V]\mapsto\chi_q(V),
\ee is an injective homomorphism of rings. Moreover, the image of
$\chi_q$ is given by \be {\mc K}=\bigcap_{i=1,\dots,N-1}
\left(\Z[Y_{j,a}^{\pm 1}]_{j\neq i}^{ a \in \C^\times} \otimes
\Z[Y_{i,b} + Y^{-1}_{i,bq^2} Y_{i-1,bq}Y_{i+1,bq}]_{b \in
\C^\times}\right), \ee where we set $Y_{N,a}=Y_{0,a}=1$.

Fix $a\in\C^\times$. Denote by $\Rep_aU_q\asln$ the subring
of $\Rep\,U_q\asln$ spanned by those representations whose
$q$--characters are polynomials in $Y^{\pm 1}_{i,aq^{i+2k-1}}$
($i=1,\dots,N-1$, $k\in\Z$).

By Corollary 6.15 in \cite{FM1}, the class of the irreducible module
$V$ is in $\Rep_aU_q\asln$ if and only if the roots of its
Drinfeld polynomials $P_i(u)$ are of the form $q^{-i-2k+1}$
($i=1,\dots,N-1$, $k \in\Z$) or, equivalently,
the highest weight monomial in $\chi_q(V)$ is a
monomial in $Y_{i,aq^{i+2k-1}}$ ($i=1,\dots,N-1$, $k \in \Z$). 
Moreover, we have the decomposition of rings \be \Rep\,U_q\asln \simeq
\bigotimes_{a\in\C^\times/q^{2\Z}}\Rep_aU_q\asln, \ee where
the tensor product is over the torus of bases of all lattices of the
form $\{aq^{2\Z}\}$ in $\C^\times$.

Similarly, let $V\in\wh{\mc F}_N$ be a polynomial
$U_q\agln$--module. Let ${v_1,\dots,v_n}$ be a basis of common
generalized eigenvectors of $\{t_i, {\mathcal Q}_i^\pm\}$ in $V$. Then 
the generalized eigenvalues
of ${\mathcal Q}_i^\pm(u)$ on each vector $v_s$ are polynomials
$$
\Gamma_{i,s}^{\pm}(u)=\prod_{k=1}^{\la_{is}}(1-a_{isk}^{\pm 1} u^{\pm
1}).
$$
Furthermore, we have: $t_i v_s = q^{\la_{i,s}} v_s$.

By definition, the $q$--character of $V$, denoted again by
$\chi_q(V)$, is an element of the ring of polynomials $\Z[\Lambda
_{i,a}]_{i=1,\dots,N}^{a\in\C^\times}$, given by the formula \be
\chi_q(V)=\sum_{s=1}^n\prod_{i=1}^N \prod_{k=1}^{\la_{is}}
\La_{i,a_{isk}}.  \ee In 
the case of an irreducible representation $V$, the monomial in
$\chi_q(V)$ corresponding to the highest weight vector is called the
highest weight monomial of $V$. By construction, the highest weight
monomial of the irreducible representation $V({\bs Q})$ of $U_q \agln$
records all information about the Drinfeld polynomials $Q_i(u)$.

Let $\Rep\,U_q\agln$ be the Grothendieck ring of the monoidal
category $\wh{\mc F}_N$. Using Lemma \ref{mult prop} and 
the one-to-one correspondence between irreducible
representations and highest weight monomials, we obtain the following

\begin{lem}    \label{ring str}
The map \be \chi_q:\; \Rep\,U_q\agln\to
\Z[\La_{i,a}]_{i=1,\dots,N}^{a\in\C^\times}, \qquad [V]\mapsto\chi_q(V),
\ee is an injective homomorphism of rings.
\end{lem}

Consider the homomorphism of rings
\bean\label{kappa}
\kappa_q:\Z[\La_{i,a}]_{i=1,\dots,N}^{a\in\C^\times}\to \Z[Y^{\pm
1}_{i,a}]_{i=1,\dots,N-1}^{a\in\C^\times}, \qquad \La_{i,a}\mapsto Y_{i,aq^{i-1}}Y^{-1}_{i-1,aq^{i}}.
\eean
From formula \eqref{PandQ} we obtain the following 

\begin{lem}
Let $V \in \wh{\mc F}_N$ be a $U_q\agln$--module with $q$--character
$\chi_q(V)$. Then the $q$--character of the restriction of $V$ to
$U_q\asln$ is equal to $\kappa_q(\chi_q(V))$.
\end{lem}

A monomial $m$ in $\Z[\La_{i,a}]_{i=1,\dots,N}^{a\in\C^\times}$ is
called {\em dominant} if $\kappa_q(m)$ does not contain negative powers
of $Y_{i,a}$. The highest weight monomial of any irreducible polynomial
representation is dominant.

Fix $a\in\C^\times$. Denote $\Rep_aU_q\agln$ the subring of
$\Rep\,U_q\agln$ spanned by all representations whose
$q$--characters belong to $\Z[\La_{i,aq^{2k}}]_{i=1,\dots,N}^{k \in
\Z}$. Then, we have an injective ring homomorphism
$$
\Rep_aU_q\agln \to \Z[\La_{i,aq^{2k}}]_{i=1,\dots,N}^{k \in \Z}.
$$
If $V$ is an irreducible polynomial representation, then $[V]$ belongs
to the subring $\Rep_aU_q\agln$ if and only if the highest weight
monomial of $V$ is a monomial in $\La_{i,aq^{2k}}$ ($k \in \Z$). We
have the decomposition of rings \be \Rep\,U_q\agln \simeq
\bigotimes_{a\in\C^\times/q^{2\Z}}\Rep_aU_q\agln, \ee where the tensor
product is over the set of bases of distinct lattices of the form
$\{aq^{2\Z}\}$ in $\C^\times$.

Irreducible polynomial representations of $U_q\agln$ form a basis in
$\Rep\,U_q\agln$ (as a free $\Z$--module). Those irreducible
representations, whose Drinfeld polynomials have roots in
$a^{-1}q^{2\Z}$, form a basis in $\Rep_aU_q\agln$. In Section
\ref{evaluation modules section} we describe other bases in
$\Rep\,U_q\agln$ and $\Rep_aU_q\agln$.

Note that the pullback with respect to the twist $\tau_{b/a}$ given
by formulas \Ref{twist} gives rise to an isomorphism
$\Rep_aU_q\agln$ with $\Rep_b U_q\agln$. This
isomorphism is compatible with the bases of irreducible
representations.

Finally, it is clear from the definition that the $q$--characters are
compatible with all embeddings discussed in Section \ref{inclusions
section}. For example, if $V$ is a $U_q\agln$--module with the
$q$--character $\chi_q(V)$ then the $q$--character of $V$ as the
$U_q\widehat{\mathfrak g\mathfrak l}_{N-1}$--module is obtained by replacing
all $\La_{N,a}$ in $\chi_q(V)$ by $1$; the character of $V$ as a
$U_q\gln$--module is obtained from $\chi_q(V)$ by replacing each
$\La_{i,a}$ by $\La_i$ corresponding to the weight 
$(0,\dots,0,1,0,\dots,0)$ with $1$ in the $i$th place.

\subsection{Evaluation modules}\label{evaluation modules section}
In this section we compute the $q$--characters of the evaluation modules.

Let us recall some standard definitions from the theory of partitions.
Fix $N\in\Z_{>0}$, and let
$\bs\la=(\la_1\geq\la_2\geq\dots\geq\la_N)$, $\la_i\in\Z_{\geq
0}$. Then $\bs\la$ is a {\em partition} of
$|\bs\la|=\sum_{i=1}^N \la_i$. Given a partition ${\bs \la}$, one draws a
picture called the {\em Young diagram} which consists of
$|\bs\la|$ boxes arranged in $N$ rows: the first row has length $\la_1$, the
second row has length $\la_2$, etc. We denote the box number $i$ in
the column number $j$ by $B_{ij}$.

A {\em semistandard tableau} of shape $\bs\la$ is the Young diagram of
$\bs\la$, in which each box is filled with a number from the set
$\{1,\dots,N\}$, such that the numbers in each row do not decrease
(from left to right) and the numbers in each column strictly increase
(from top to bottom). We denote the number in the box $B_{ij}$ by
$f(i,j)$.

Here is an example of a semistandard tableau with $N=4$,
$\bs\la=(4,2,2,1)$.

\begin{picture}(20,30)(0,20)
\setlength{\unitlength}{1mm}
\thicklines
\put(50,40){\framebox(5,5){1}}
\put(55.3,40){\framebox(5,5){1}}
\put(60.6,40){\framebox(5,5){2}}
\put(65.9,40){\framebox(5,5){3}}

\put(50,34.7){\framebox(5,5){2}}
\put(55.3,34.7){\framebox(5,5){2}}

\put(50,29.42){\framebox(5,5){3}}
\put(55.3,29.42){\framebox(5,5){4}}

\put(50,24.1){\framebox(5,5){4}}

\end{picture}

Recall from \cite{GZ} that the restriction of the irreducible
$\gw_N$--module with highest weight $\bs\la$ to
$\gw_{N-1}$ (with respect to the left upper corner inclusions) 
is multiplicity free. Moreover, the weights 
occurring in this decomposition,
$(\la_{N-1,1},\ldots,\la_{N-1,N-1})$,
are all the
weights that satisfy $\la_i\geq \la_{N-1,i} \geq \la_{i+1}$.
Decomposing further with respect to $\gw_{N-2}, \gw_{N-3}$, etc., we
eventually come to irreducible representations of $\gw_1$, which are
one-dimensional. Therefore we obtain a decomposition of $V_{\bs \la}$
into a direct sum of one-dimensional subspaces parameterized by the
Gelfand-Zetlin schemes:
$$
\begin{array}{ccccccccc}
\la_{1} & & \la_{2} &  & \ldots &  & \ldots &
& \la_{N} \\ & \la_{N-1,1} & & \la_{N-1,2} & & \ldots & &
\la_{N-1,N-1} & \\ & & \ldots & \ldots & \ldots & \ldots & \ldots & &
\\ & & & \la_{2,1} & & \la_{2,2} & & & \\ & & & & \la_{1,1} & & & &
\end{array}
$$
where $\la_{M,i} \geq \la_{M-1,i} \geq \la_{M,i+1}$ for all $i,M$ (we
set $\la_{N,i} = \la_i$). To such a scheme we attach a
semistandard tableau by the following rule. We embed the Young
diagram of the partition $(\la_{i,1},\ldots,\la_{i,i})$ into the Young
diagram of ${\bs \la}$ by identifying the boxes in the left top
corners. We refer to it as the $i$th subdiagram of ${\bs \la}$. Then
we fill all boxes which belong to the $i$th subdiagram but not to the
$(i-1)$th diagram with the number $i$. It is clear that this procedure
gives rise to a bijection between the Gelfand-Zetlin schemes 
and the semistandard tableaux of shape ${\bs \la}$. 
In our example, the semi-standard tableau above corresponds to the following
Gelfand-Zetlin scheme:
$$
\begin{array}{ccccccc}
4& &2& &2& &1  \\ 
 &4& &2& &1&   \\
 & &3& &2& &   \\
 & & &2& & &   
\end{array}
$$

Define the {\em content}, $c(i,j)$, of the box $B_{ij}$ in a Young
diagram by the formula $j-i$. In our example we have

\begin{picture}(20,30)(0,20)
\setlength{\unitlength}{1mm}
\thicklines
\put(50,40){\framebox(5,5){0}}
\put(55.3,40){\framebox(5,5){1}}
\put(60.6,40){\framebox(5,5){2}}
\put(65.9,40){\framebox(5,5){3}}

\put(50,34.7){\framebox(5,5){-1}}
\put(55.3,34.7){\framebox(5,5){0}}

\put(50,29.42){\framebox(5,5){-2}}
\put(55.3,29.42){\framebox(5,5){-1}}

\put(50,24.1){\framebox(5,5){-3}}

\end{picture}

Now to a semistandard tableau $\bs t$ and a number $a\in\C^\times$, we
attach a monomial ${M}_{a,\bs t} \in
\Z[\La_{i,b}]_{i=1,\dots,N}^{b\in\C^\times}$ by the formula \be
{M}_{a,\bs t}=\prod_{i,j}\La_{f(i,j),aq^{2c(i,j)}}.  \ee

For example, to the above tableau we attach the monomial 
\be
\La_{1,a}\La_{1,aq^2}\La_{2,aq^{-2}}\La_{2,a}\La_{2,aq^{4}}
\La_{3,aq^{-4}}\La_{3,aq^6}\La_{4,aq^{-6}}\La_{4,aq^{-2}}. 
\ee

Let $V$ be a representation of $U_q\gln$ from the category ${\mc
F}_N$. Then the representation of $U_q\agln$ obtained by the pullback
with respect to the evaluation homomorphism $ev_a$ defined by formula
\Ref{gl evaluation} is called {\em evaluation module} and is denoted by
$V(a)$.

\begin{lem}\label{evaluation character lemma}
Let $V_{\bs\la}$ be the irreducible representation of $U_q\gln$ with
highest weight $\bs\la=(\la_1,\dots,\la_N)$, where $\la_i\in\Z_{\geq
0}$. Then $V_{\bs\la}(a)$ is a polynomial 
$U_q\agln$--module. Moreover, the $q$--character of $V_{\bs\la}(a)$
is given by the formula \be \chi_q(V_{\bs\la}(a))=\sum_{\bs t}
M_{a,\bs t},
\ee where the sum is over all semistandard tableaux of shape
$\bs\la$.
\end{lem}
\begin{proof}
The proof proceeds along the lines of the proof of a similar statement
in the case of the Yangians given in \cite{NT}, Lemma 2.1 (see also
\cite{Ch}; a similar formula may be found in \cite{BR} in the
Bethe Ansatz context).

Since finite-dimensional representations of $U_q \gln$ are
$q$--deformations of representations of $\gln$, the patterns
for restrictions of $U_q \gln$--modules  and $\gln$--modules 
are the same. In particular, for each
Gelfand-Zetlin scheme $S$ (equivalently, a semistandard tableau $\bs t$ of
shape ${\bs \la}$), there is a vector $v_S$ of 
$V_{\bs \la}$, which is contained in the $U_q
\gw_M$--component of $V_{\bs \la}$ with highest weight
$(\la_{M,1},\ldots,\la_{M,M})$. We claim that each of these
vectors $v_S$ is a common eigenspace of $\{t_i,{\mc Q}^+_i(u)\}$ which
contributes the monomial $M_{a,\bs t}$ to $\chi_q(V_{\bs\la}(a))$.

In order to show this, we consider an action of the series
$C^\pm_{(M)}(z) = \prod_{i=1}^M\La^\pm_i(zq^{2-2i})$ on $V_{\bs
\la}(a)$. Using formulas \eqref{La and Theta} and \eqref{theta q}, we
can then recover the action of ${\mc Q}^\pm_i(u)$ and $t_i$. 

The coefficients of $C^\pm_{(M)}(z)$ are central elements of $U_q
\aglm$ and they act by scalar operators on each irreducible $U_q \aglm$
module.  In particular each $v_S$ is an eigenvector of
$C^\pm_{(M)}(z)$. 

Recall that by Theorem \ref{cube theorem}, 
the evaluation homomorphism commutes with the restrictions. Therefore, 
if $W$ is an irreducible $U_q \glm$ submodule of $V_{\bs\la}$, then $W(a)$ 
is an irreducible $U_q\aglm$ submodule in $V_{\bs\la}(a)$.

To compute the action of $C^\pm_{(M)}(z)$ on $v_S$, we consider
the lowest weight vector $w_S$ of the corresponding $U_q \glm$ module.
In other words, $w_S$ is the lowest weight vector of the
the evaluation  $U_q \aglm$ module containing $v_S$. 
Note that by
the weight decomposition of $V_{\bs \la}$, the
vector $w_S$ is an eigenvector of each $\La^\pm_i(z)$, $i=1,\dots,M$.
Let us compute these eigenvalues.

Consider any evaluation module $V_{\bs \mu}$ of $U_q \wh{\gw}_M$. 
It follows from the formulas for the triangular
decomposition of $L^\pm(z)$ that the action of $\La^\pm_i(z)$
($i=1,\dots,M$) 
on a lowest weight vector $\wt{v}_{\bs \mu}$ of $V_{\bs \mu}$ 
coincides with the action of
$L_{ii}^\mp(z)$. Note that $L_{ii}^\pm=T_i^{\pm1}$ by Lemma \ref{gl dj
to r}, therefore using \Ref{gl evaluation} we compute
$$L^\mp_{ii}(z) \wt{v}_{\bs \mu} = \frac{T_i^{\mp 1} - (z/a)^{\mp
1}T_i^{\pm 1}}{1-(z/a)^{\mp 1}} \; \wt{v}_{\bs \mu} = \frac{q^{\mp
\mu_{M+1-i}} - (z/a)^{\mp 1} q^{\pm \mu_{M+1-i}}}{1-(z/a)^{\mp 1}} \;
\wt{v}_{\bs \mu}$$  
(note that the weight of $\wt{v}_{\bs \mu}$  is
$(\mu_M,\mu_{M-1},\ldots,\mu_1)$). 
Therefore, using formula \eqref{La
and Theta} we find that the eigenvalue of $\Theta^\pm_M(z^{-1}
q^{2M-1})$ on $v_S$ is
$$
\prod_{i=1}^{M-1} \frac{q^{\mp \la_{M-1,M-i}} - (zq^{-2i}/a)^{\mp 1}
q^{\pm \la_{M-1,M-i}}}{q^{\mp \la_{M,M-i}} - (zq^{-2i}/a)^{\mp 1}
q^{\pm \la_{M,M-i}}} \frac{1-(z/a)^{\mp
1}}{q^{\mp \la_{M,M}} - (z/a)^{\mp 1} q^{\pm \la_{M,M}}} .
$$
Using formula \eqref{theta q} and the above rule assigning to $S$
a semistandard tableau, we obtain the statement of the lemma.
\end{proof}

In particular, the highest weight monomial of $V_{\bs\la}(a)$
corresponds to the tableau, in which $f(i,j)=i$ (i.e., it is filled
with $1$'s in the first row, $2$'s in the second row, etc.). This
describes the Drinfeld polynomials of the evaluation modules.

Let $\om_i=(1,\dots,1,1,0,\dots,0)$, where the last $1$ is in the $i$th
position, be the fundamental weights of $\gln$. 
The evaluation modules $V_{\om_i}(a)$ are called the
{\em fundamental} representations of $U_q \agln$. According to
\lemref{evaluation character lemma}, we have
\begin{equation}    \label{ti}
\chi_q(V_{\omega_i}(a)) = \sum_{1\leq j_1<j_2<\ldots<j_i\leq N}
\Lambda_{j_1,a} \Lambda_{j_2,aq^{-2}} \ldots
\Lambda_{j_i,aq^{-2+2i}}.
\end{equation}
In particular, the highest weight monomial of $V_{\om_i}(a)$ is equal
to $\prod_{j=1}^i\La_{j,aq^{2-2j}}$. Note that for $i=1,\dots,N-1$,
the restriction of $V_{\om_i}(a)$ to $U_q\asln$ gives us the
representation with Drinfeld polynomials ${\bs
P}=(1,\dots,1,1-aq^{1-i}u,1,\dots,1)$, where $1-aq^{1-i}u$ is in $i$th
position. The representation $V_{\om_N}(a)$ is the one-dimensional
determinant representation of $U_q \agln$ and its restriction to $U_q
\asln$ is trivial.

{}From the uniqueness of an irreducible polynomial representation of
$U_q \agln$ with a given highest weight monomial (equivalently, a set
of Drinfeld polynomials) and the multiplicative property of the
$q$--characters (see \lemref{ring str}) we obtain the following
 
\begin{lem}
Any irreducible
finite-dimensional polynomial representation $V$ of $U_q\agln$ occurs as a
quotient of the submodule of the tensor product $V_{\omega_{i_1}}(a_1)
\otimes \ldots \otimes V_{\omega_{i_n}}(a_n)$, generated by the tensor
product of the highest weight vectors. Moreover, the parameters
$(\omega_{i_k},a_k)$, where $k=1,\dots,n$, are uniquely determined by $V$ up
to a permutation.
\end{lem}

Therefore, we have a basis in $\Rep\,U_q\agln$ given by the tensor
products of the fundamental representations $V_{\om_i}(a)$
($i=1,\dots,N$, $a\in\C^\times$). We also have a basis in $\Rep_a
U_q\agln$ given by the tensor products of the fundamental
representations $V_{\om_i}(aq^{2k})$ ($i=1,\dots,N$,
$k\in\Z$). Therefore we have isomorphisms
\begin{equation}    \label{iso repa}
\Rep_a U_q\agln \simeq
\Z[V_{\omega_i}(aq^{2k})]_{i=1,\dots,N}^{k\in\Z}, \qquad \Rep U_q\agln \simeq
\Z[V_{\omega_i}(b)]_{i=1,\dots,N}^{b\in\C^\times}.
\end{equation}

\section{The inductive limit of 
$\Rep_a U_q\agln$}\label{ind limit def generic section}

In this section we define the inductive limit
$\underset{\longrightarrow}{\lim}\;\Rep \,U_q\agln$ and a Hopf
algebra structure on it. 

\subsection{Definition of the inductive limit}\label{def of ind section}

In this section we define the ring $\Rep\,U_q\agli$ together
with two linear bases, corresponding to tensor products of the
fundamental representations and to irreducible representations
respectively.

Let ${\mc T}_N=\Z[t_{i,n}]_{i=1,\dots,N}^{n\in\Z}$, ${\mathcal
T}_\infty=\Z[t_{i,n}]_{i\in\Z_{>0}}^{n\in\Z}$ be commutative rings
of polynomials. We introduce a $\Z_{\geq 0}$--gradation on them by
setting the degree of $t_{i,n}$ to be equal to $i$.

For $N\in\Z_{>0}$, we have embeddings of graded rings 
\be
{\mc T}_N\hookrightarrow {\mc T}_{N+1},\qquad t_{i,n}\mapsto
t_{i,n}.  \ee 
Then the ring ${\mc T}_\infty$ is the inductive
limit $\underset{\longrightarrow}{\lim}\;{\mc T}_N$.

Fix $a\in \C^{\times}$. According to formula \eqref{iso repa}, for
$N\in\Z_{>0}$, we have isomorphisms of rings 
\bean\label{fund basis}
\Rep_a\,U_q\agln \overset{\sim}\to {\mathcal T}_N,\qquad 
[V_{\om_i}(aq^{2n})]\mapsto t_{i,n}.  
\eean 
We call the ring ${\mc
T}_\infty \simeq \underset{\longrightarrow}{\lim}\;
\Rep_a U_q\agln$ the Grothendieck ring of $U_q\agli$
(restricted to the lattice with base $a$), and denote it by
$\Rep_a U_q\agli$. The monomial basis of ${\mathcal
T}_\infty=\Rep_a U_q\agli$ corresponds to the basis of tensor
products of the fundamental representations.

We also set \bean    \label{decomp for a}
\Rep \,U_q\agli:=\bigotimes_{a\in\C^{\times}/q^{2\Z}}
\Rep_aU_q\agli.  \eean Thus we have: \bea
\Rep_a U_q\agli&=&\underset{\longrightarrow}{\lim}\;
\Rep_a U_q\agln,\\
\Rep\,U_q\agli&=&\underset{\longrightarrow}{\lim}\;
\Rep\,U_q\agln.  \eea

\begin{remark}
Note that we may define a similar way the inductive limit
$\Rep_aU_q\widehat{\mathfrak{sl}}_\infty$ of the Grothendieck
rings $\Rep_a U_q \widehat{\mathfrak{sl}}_N$. Since
$\Rep_a U_q\asln$ is identified with the subring of
$\Rep_a U_q\agln={\mc T}_N$ generated by $t_{i,n}$,
$i=1,\dots,N-1$, we obtain that 
$$\Rep_a U_q\widehat{\mathfrak{sl}}_\infty \simeq {\mc
T}_\infty \simeq \Rep_a U_q\agli.$$
However, it is more convenient to work with the inductive limit of
$\Rep_a U_q\agln$, rather than of $\Rep_a U_q\asln$,
because for $N=N_1+N_2$ there exist an algebra homomorphism $U_q
\widehat{\mathfrak{gl}}_{N_1} \otimes U_q
\widehat{\mathfrak{gl}}_{N_2} \to U_q
\widehat{\mathfrak{gl}}_{N}$ which induces an isomorphism on the
``Cartan'' parts (generated by $t_i^{\pm1}$, $g_{i,n}$). We use this
homomorphism extensively in what follows.\qed
\end{remark}

Now we show that irreducible representations of $U_q \agln$ stabilize
in the limit $N \to \infty$. Namely, each irreducible representation
of $U_q\agln$ defines an element of
$\Rep\,U_q\agli$. \lemref{irrep to fund stable} proved
below shows that irreducible representations $V_1$ and $V_2$ of
$U_q\agln$ and $U_q\aglm$, respectively, where $N\leq M$, define the
same element of $\Rep\,U_q\agli$ if and only if the sequence
of Drinfeld polynomials of $V_2$ coincides with the sequence of the
Drinfeld polynomials of $V_1$ extended by $1$'s. This implies that
$\Rep \,U_q\agli$ has a basis corresponding to irreducible
representations.

More precisely, let ${\bs Q}=(Q_1(u),\dots,Q_N(u))$ be a dominant 
$N$--tuple of polynomials and $V({\bs Q})$ the
corresponding $U_q \agln$--module. Set \be \deg {\bs Q}=\sum_{i=1}^N \deg
Q_i.  \ee We say that the degree of $V({\bs Q})$ is $\deg{\bs Q}$.

For $M\geq N$, let ${\bs Q}^{(M)}$ be the $M$--tuple of polynomials
$(Q_1(u),\dots, Q_N(u),1,\dots,1)$.  An $M$--tuple of polynomials
${\bs Q}_1$ and an $N$--tuple of polynomials ${\bs Q}_2$ are called
equivalent if ${\bs Q}_1={\bs Q}_2^{(M)}$. Note that if ${\bs Q}_1$
and ${\bs Q}_2$ are equivalent, then $\deg {\bs Q}_1=\deg {\bs Q}_2$.

An $N$--tuple ${\bs Q}$ of polynomials with constant term $1$ is
called {\em minimal} if $Q_N(u) \neq 1$ or $N=1$. An $M$--tuple ${\bs Q}$ is
equivalent to exactly one minimal $N$--tuple with some $N\leq M$. We
call $N$ the length of the $M$--tuple ${\bs Q}$.

Let $m\in {\mc T}_N\subset{\mc T}_M$ be a monomial in $t_{i,n}$ of degree $d$. For
$M\geq N$, let $W^{(M)}(m)$ be the corresponding element of
$\Rep_a U_q\aglm$ (the class of the tensor product of the fundamental
representations) defined by the isomorphism \eqref{fund basis}. We can
represent it in $\Rep_a U_q\aglm$ uniquely as a linear
combination of irreducible representations \bean\label{fund to irrep}
W^{(M)}(m)=\sum_i c_i^{(M)}(m)\; [V({\bs Q}_i)], \eean with some non-negative
integral coefficients $c_i^{(M)}(m)$.

Introduce a $\Z$--gradation on the ring ${\mc L}_N :=
\Z[\La_{i,aq^{2k}}]_{i=1,\dots,N}^{k \in \Z}$ by setting $\deg
\Lambda_{i,aq^{2k}} = 1$. According to formula \eqref{ti}, the
$q$--character homomorphism $\chi_q: {\mc T}_N \to {\mc L}_N$
preserves the gradation. This implies that if $c_i(m)\neq 0$ in
formula \eqref{fund to irrep}, then $\deg{\bs Q}_i=d$, so in
particular the length of ${\bs Q}_i$ is at most $d$. Therefore the set
of possible ${\bs Q}_i$ appearing in formula \eqref{fund to irrep} is
the same for all $M\geq d$, and we can rewrite formula \eqref{fund to
irrep} for $M\geq d$ as
\begin{equation}    \label{fund to irrep1}
W^{(M)}(m)=\sum_i c_i^{(M)}(m)\; [V({\bs Q}_i^{(M)})],
\end{equation}
where the sum is over all $d$--tuples of monomials ${\bs Q}_i$ of
degree $d$. Note that only finitely many coefficients $c_i^{(M)}(m)$
are nonzero in this formula.

\begin{lem}\label{fund to irrep stable}
There exists $M_0\geq d$, such that $c_i^{(M)}(m)=c_i^{(M_0)}(m)$ for
all $M\geq M_0$ and all $i$.
\end{lem}

In order to prove this lemma, we need to introduce our main technical
tool: the restriction operators.

\subsection{Definition of the restriction operators}\label{def restr
section}

The algebra $U_q \widehat{\mathfrak gl}_{N}$ contains subalgebras
isomorphic to $U_q \widehat{\gw}_{N-1}$ and to $U_q
\widehat{\gw}_1$. The first one is defined by the upper left corner
inclusion in Section \ref{inclusions
section}, and the second one is the subalgebra generated (in the
``new'' realization) by $t_N$, $g_{N,n}$ ($n\in\Z\setminus
0$). According to the relations of \secref{agln}, the two subalgebras
commute inside $U_q \agln$. Hence we obtain an embedding $U_q
\widehat{\gw}_{N-1} \otimes U_q \widehat{\gw}_1 \hookrightarrow U_q
\widehat{\mathfrak gl}_{N}$.

Observe that the isomorphism \eqref{fund basis} identifies the ring
$\Rep_a U_q \widehat{\mathfrak g}{\mathfrak l}_1$ with the ring of
polynomials in variables $X_n:=t_{1,n}$  ($n\in\Z$).

The restriction of a polynomial representation of $U_q\agln$ 
to $U_q \widehat{\gw}_{N-1}\otimes U_q \widehat{\gw}_1$
belongs to the product of the categories $\wh{\mc F}_{N-1}$ and
$\wh{\mc F}_1$. Thus, we obtain a functor $\wh{\mc F}_N \to \wh{\mc
F}_{N-1} \times \wh{\mc F}_1$ and a linear map \be \Rep_a U_q\agln
\to \Rep_a U_q \wh{\gw}_{N-1} \otimes \Rep_a U_q
\widehat{\gw}_1.  \ee 
Under this map, \bean\label{restriction
decomp} [V] \mapsto \sum_{\bs m} [V_{\bs m}]\otimes X_{\bs m}, \eean
where the sum is over all ${\bs m}=(m_1\geq\dots\geq m_k)$, $m_i\in\Z$
and $X_{\bs m}:= X_{m_1}X_{m_2}\dots X_{m_k}$.

For each ordered set ${\bs m}=(m_1\geq\dots\geq m_k)$, $m_i\in\Z$
(empty if $k=0$) we define the restriction operator $\res^{(N)}_{\bs m}$
by the formula
\bean    \label{resm}
\res^{(N)}_{\bs m} :\;\Rep_a U_q\widehat{\gw}_N 
\to  \Rep_a U_q\widehat{\gw}_{N-1},\qquad {[V]}  \mapsto  [V_{\bs
m}].  \eean

Recall the notation ${\mc L}_M =
\Z[\Lambda_{i,aq^{2n}}]_{i=1,\ldots,M}^{n\in\Z}$ and consider the
natural isomorphism of rings $\al_N: {\mc L}_N \to {\mc L}_{N-1}
\otimes {\mc L}_1$. Given $P \in {\mc L}_N$, we write
\bean    \label{alN}
\al_N(P) = \sum_{\bs m} P_{\bs
m}\otimes(\Lambda_{N,aq^{2m_1}}\La_{N,aq^{2m_2}}\dots
\La_{N,aq^{2m_k}}),
\eean
where the sum is over all ${\bs m}=(m_1\geq\dots\geq m_k)$, $m_i\in\Z$.

Define a linear operator $\wt{\res}^{(N)}_{\bs m}: {\mc L}_N
\to {\mc L}_{N-1}$ sending $P$ to $P_{\bs m}$. 

The $q$--character of the $U_q \widehat{\gw}_1$--module corresponding
to $X_i$ equals $\La_{1,q^{2i}}$. Therefore
\begin{equation}    \label{restr and lambda}
\chi_q(\res^{(N)}_{\bs m} [V]) = \wt{\res}^{(N)}_{\bs m} (\chi_q(V)).
\end{equation}
In other words, the $q$--character of $\res^{(N)}_{\bs m} [V]$ is obtained
from $\chi(V)$ by collecting all monomials of $\chi_q(V)$ which are
products of $P \in {\mc L}_{N-1}$ and $\Lambda_{N,\bs m}$, and setting
$\Lambda_{N,\bs m}=1$.

In particular, we obtain the following

\begin{lem}\label{union}
\bea
\res^{(N)}_{\bs m}[V_1\otimes V_2]=\sum_{{\bs m}={\bs m_1}\sqcup {\bs
m_2}}\;\res^{(N)}_{\bs m_1}[V_1]\otimes \res^{(N)}_{\bs m_2}[V_2].
\eea
Here and below the summation is over all possible splittings of the
sequence ${\bs m}$ into a disjoint union of non-increasing
subsequences, with each splitting counted only once.
\end{lem}

Using formulas \eqref{restr and lambda} and  \eqref{ti} we find the following

\begin{lem}\label{action lemma}
For the fundamental representation we have
\be
\res^{(N)}_{\bs m} (t_{i,j})=\left\{\begin{matrix}
t_{i,j}& {\rm if}\;\;{\bs m}=\emptyset,\;\; j<N,\\
t_{i-1,j}& {\rm if}\;\;{\bs m}=(j-i+1),\\
0&{\rm otherwise}.
\end{matrix}\right.
\ee
\end{lem}

These formulas do not depend on $N$ if $N$ is large enough. Using this
fact and Lemma \ref{union} we can pass to the inductive limit.

For ${\bs m}=(m_1\geq\dots\geq m_k)$, $m_i\in\Z$, define the
operator
\be
\res_{\bs m}: \; \Rep_aU_q\agli\to \Rep_aU_q\agli
\ee
by the formulas
\bean \label{res t}
\res_{\bs m} (t_{i,j})=\left\{\begin{matrix}
t_{i,j}& {\rm if}\;\;{\bs m}=\emptyset,\\
t_{i-1,j}&{\rm if}\;\; {\bs m}=(j-i+1),\\
0&{\rm otherwise},
\end{matrix}\right.\\ \label{res tensor}
\res_{\bs m}[V_1\otimes V_2]=\sum_{{\bs m}={\bs m_1}\sqcup {\bs
m_2}}\;\res_{\bs m_1}[V_1]\otimes \res_{\bs m_2}[V_2].  \eean 

We denote the operators $\res_{(i)}$, $i\in\Z$ simply by $\res_i$.
Note that the operators $\res_i$ ($i\in\Z$) are derivations of the ring
$\Rep_aU_q\agli$.

\subsection{Stabilization of irreducible
representations}\label{stab section}

Now we are ready to prove \lemref{fund to irrep stable}.

\medskip

\noindent {\em Proof of \lemref{fund to irrep stable}.} For the sake
of brevity, we
use the simplified notation, $W^{(M)} \equiv W^{(M)}(m)$ and
$c_i^{(M)} \equiv c_i^{(M)}(m)$.

Choose a coefficient $c_i^{(M)}$ in the RHS of \Ref{fund to
irrep1}. Since $c_i^{(M)}\geq 0$, it is sufficient to prove that $\{
c_i^{(M)} \}_{M\geq N}$ form a non-increasing sequence, i.e.,
$c_i^{(M)}\leq c_i^{(N)}$, for $M\geq N$.

According to our assumption, $[W^{(M+1)}]$ is a product of
$t_{i,n}$ with $i\leq M$. By \lemref{action lemma},
$\res^{(M+1)}_{\emptyset}(t_{i,n}) = t_{i,n}$ for all $i\leq M$. Using
Lemma \ref{union} we obtain $$\res^{(M+1)}_{\emptyset} W^{(M+1)} =
W^{(M)}.$$ In other words, $\chi_q(W^{(M)})$ is exactly the sum
of all monomials in $\chi_q(W^{(M+1)})$ which do not contain
$\La_{M+1,n}$ ($n\in\Z$).

Thus, for any $d$--tuple ${\bs Q}_i$ with $d\leq M$,
we have that  $\res_{\emptyset}^{(M+1)} [V({\bs Q}_i^{(M+1)})]$ is a sum of
$[V({\bs Q}_i^{(M)})]$ with a linear combination of other irreducible
$U_q \wh{\gw}_M$--modules with non-negative coefficients. Applying
$\res_{\emptyset}^{(M+1)}$ to both sides 
of the equality
$$
W^{(M+1)}=\sum_i c_i^{(M+1)}\; [V({\bs Q}_i^{(M+1)})],
$$
we find that $W^{(M)}$ equals $\sum_i c_i^{(M+1)}\; [V({\bs
Q}_i^{(M)}]$ plus a combination with non-negative coefficients of
other irreducible $U_q \wh{\gw}_M$--modules. Therefore $c_i^{(M+1)}
\leq c_i^{(M)}$.\qed

\bigskip

Actually, one can prove that $M_0$ in Lemma \ref{fund to irrep stable}
can be chosen to be equal to $d$, but we do not need this exact bound
for the purposes of this paper.

We also have a decomposition inverse to \Ref{fund to irrep1}. Fix a
dominant 
$N$--tuple of polynomials ${\bs Q}$
with roots in $aq^{2\Z}$.
Then for all $M\geq N$ we have an equality in
$\Rep_a U_q\aglm$, \bean\label{irrep to fund} V({\bs
Q}^{(M)})=\sum_i b_i^{(M)}({\bs Q})\; W^{(M)}(m_i), \eean where
$b_i^{(M)}$ are some integers and $W^{(M)}(m_i)$ are tensor products
of fundamental representations corresponding to monomials
$m_i\in\mc T_\infty$ of degree $d_i$. For each $M$ in this sum
only finitely many terms are nonzero.

\begin{lem}\label{irrep to fund stable}
If $b_i\neq 0$, then $d_i=\deg{\bs Q}$. Moreover there exists $M_0\geq
\deg{\bs Q}$, such that $b_i^{(M)}=b_i^{(M_0)}$ for all $M\geq M_0$
and all $i$.
\end{lem}
\begin{proof}
Let $d$ be the degree of ${\bs Q}$. The transition matrix from $\{
W(m_i) \}$ to $\{ V({\bs Q}_j) \}$ with entries $b_i^{(M)}({\bs Q}_j)$
is inverse to the transition matrix from $\{ V({\bs Q}_i) \}$ to $\{
W(m_j) \}$ with entries $\{c_i^{(M)}(m)\}$. Since $V({\bs Q})$ is a
homogeneous element of the ring $\Rep_a U_q\aglm$ of degree $\deg
{\bs Q}$, these matrices are unions of blocks corresponding to fixed
degrees. These blocks are still of infinite size. However, we can
further decompose them into blocks of finite size.

Fix $n_1,n_2\in\Z$, $n_1\leq n_2$ and consider the subring ${\mathcal
X}^{n_1,n_2}_N$ of ${\mathcal T}_N$ generated by $t_{i,n}$, where
$i\in \Z_{>0}$, $n_1\leq n\leq n_2$. We claim that the image of the
subring ${\mc X}^{n_1,n_2}_N$ in $\Rep_a U_q\agln$ has a basis of
irreducible representations whose $q$--character contains $\La_{1,n}$
only with $n\in[n_1,n_2]$.
 
Indeed, let $m$ be a monomial in ${\mc X}^{n_1,n_2}_N$. By formula
\eqref{ti}, the $q$--character of the corresponding tensor product of
fundamental representations contains $\La_{1,n}$ only with
$n\in[n_1,n_2]$. In particular, the $q$--characters of all irreducible
representations in the RHS of \Ref{fund to irrep1} contain $\La_{1,n}$
only with $n\in[n_1,n_2]$.

On the other hand, let $V$ be an irreducible representation whose
$q$--character contains $\La_{1,n}$ only with
$n\in[n_1,n_2]$. Represent $V$ as a linear combination of tensor
products of the fundamental representations as in formula \Ref{irrep
to fund}. Then all $W(m_i)$'s in the RHS belong to ${\mc X}^{n_1,n_2}_N$.
Otherwise, the highest weight monomial of the $q$--character of at
least one such $W(m_i)$'s would appear with a nonzero coefficient in
the $q$--character of the RHS of \eqref{irrep to fund}. This monomial
would then contain $\La_{1,n}$ with some $n\not\in[n_1,n_2]$. This
contradicts our assumption and hence proves the claim.

Finally, observe that the subspace of ${\mc X}^{n_1,n_2}_\infty$ of
degree $d$ is finite-dimensional, so our transition matrices are
unions of blocks of finite size. The statement of the lemma now
follows from \lemref{fund to irrep stable}.
\end{proof}

From Lemma \ref{irrep to fund stable} we obtain the following
\begin{cor}\label{irrep form basis}
The ring $\Rep \,U_q\agli$ has a basis corresponding to irreducible
polynomial representations of $U_q \agln$ with minimal dominant
highest weights ($N\in\Z_{\geq 0}$).
\end{cor}

In particular, the ring
$\Rep_a U_q\agli$ has two linear bases, 
corresponding to all tensor products of the
fundamental representations and to the irreducible
representations respectively. 
The former basis is just the monomial basis of the
ring ${\mc T}_\infty$. We call it the {\em PBW basis}. The second basis is
labeled by minimal dominant $N$--tuples of polynomials, $N\in\Z_{>0}$,
with roots in $\{aq^{2\Z}\}$. We call it the {\em canonical
basis}. The 
transition matrix from the PBW to the canonical basis is a
block-triangular matrix with non-negative integer entries and $1$'s on
the diagonal.

\subsection{A Hopf algebra structure on $\Rep\, U_q
\agli$}\label{hopf str section}

In this section we introduce operations of multiplication and
comultiplication on $\Rep\, U_q \agli$ and show that they combine
into the structure of a Hopf algebra.

First we define the operation of multiplication in $\Rep\, \agli$.
For each finite $N$ we have an associative multiplication on
$\Rep\,\agln$ induced by the operation of tensor product of
representations. It is clear from the definitions that the operation
of multiplication of tensor products of the fundamental
representations is stable in the limit $N\to \infty$. Therefore we
obtain a well-defined operation of multiplication on $\Rep\, \agli$.

In particular, this implies that the operation of multiplication of
the irreducible representations is also stable. Indeed, fix an
$N$--tuple ${\bs P}$ and an $L$--tuple ${\bs Q}$ of polynomials with
constant term $1$. Then for all $M\geq \max\{N,L\}$ we have an
equality in $\Rep \,U_q\aglm$, \bean\label{mult of irrep} [V({\bs
P}^{(M)})\otimes V({\bs Q}^{(M)})] =\sum a_i^{(M)}({\bs P},{\bs Q})
[V({\bs R}_i^{(M)})], \eean where $a_i^{(M)}({\bs P},{\bs Q})$ are
some non-negative integers. Here the sum is over minimal
$S$--tuples ${\bs R}_i$ with $S\leq M$. For each $M$, only finitely
many terms are non-zero in this sum.

\begin{cor}\label{tensor prod stable}
If $a_i^{(M)}({\bs P},{\bs Q})\neq 0$, then $\deg{\bs P}+\deg{\bs
Q}=\deg{\bs R}_i$. Moreover there exists $M_0\geq \deg{\bs P}+\deg{\bs
Q}$, such that $a_i^{(M)}({\bs P},{\bs Q})=a_i^{(M_0)}({\bs P},{\bs
Q})$ for all $M\geq M_0$ and all $i$.
\end{cor}

\begin{proof}
We express $V^{(M)}({\bs P})$ and $V^{(M)}({\bs Q})$ as linear
combinations of tensor products of the fundamental representations,
then multiply them out and express the result as a linear combination
of irreducible representations.  All of these operation are stable (see
Lemma \ref{fund to irrep stable} and Lemma \ref{irrep to fund
stable}).
\end{proof}

Now we define the comultiplication on $\Rep\, U_q \agli$.

Let $N,N_1,N_2$ be natural numbers and $N=N_1+N_2$. Using the ``new``
realization of $U_q\agln$, we define its subalgebras
$U_1=U_q{\widehat{\gw}_{N_1}}$ 
generated by $x_{i,n}^\pm$ ($i=1,\dots,N_1-1$, $n\in\Z$), $t_i^{\pm
1}$, $g_{i,n}$, ($i=1,\dots,N_1$, $n\in\Z\setminus 0$) and
$U_2=U_q{\widehat{\gw}_{N_2}}$ generated by
$x_{i,n}^\pm$ ($i=N_1+1,\dots,N-1$, $n\in\Z$), $t_i^{\pm 1}$,
$g_{i,n}$ ($i=N_1+1,\dots,N$, $n\in\Z\setminus 0$).

It follows from the defining relations that the algebras $U_1$ and
$U_2$ commute inside $U_q\agln$. Therefore, at the level of the
Grothendieck rings, the restriction of a representation $V$ of
$U_q\agln$ to $U_1\otimes U_2$ is uniquely written in the form
\bean\label{hopf decomp} [V] = \bigoplus\limits_i\; [V_i^{(1)}]\otimes
[V_i^{(2)}], \eean where $V_i^{(j)}$ are irreducible
$U_j$--modules ($j=1,2$). Define a linear map 
\bea
\Delta_{N_1, N_2}^N:\; \Rep_a U_q\agln &\to&  \Rep_a
U_q{\widehat{\gw}_{N_1}}  \otimes \Rep_a
U_q{\widehat{\gw}_{N_2}},\\
{[V]} &\mapsto& \sum_i [V_i^{(1)}]\otimes [V_i^{(2)}].
\eea

In the language of $q$--characters, $\Delta_{N_1, N_2}^N$ may be
described as follows: we decompose the $q$--character of $V$, which is
a polynomial in $\La_{i,a}$ ($i=1,\dots,N,$ $a\in\C^\times$) into
the sum of the products of polynomials in $\La_{i,a}$
($i=1,\dots,N_1$, $a\in\C^\times$) and $\La_{j,b}$
($j=N_1+1,\dots,N$, $b\in\C^\times$), which are the $q$--characters of
representations $V_i^{(1)}$ and $V_i^{(2)}$ of $U_1$ and $U_2$,
respectively.

In particular, we obtain that for any $a\in \C^\times$,
$$
\Delta_{N_1, N_2}^N(\Rep_a \agln) \subset \Rep_a U_1 \otimes
\Rep_a U_2.
$$

We list simple properties of the maps $\Delta^{N}_{N_1,N_2}$, which
follow directly from the definitions.

\begin{lem}\label{compatibility} For $U_q\agln$ modules $V$ and $W$ we
have \bean \label{delta of product}
\Delta^N_{N_1,N_2}([V\otimes W])=\sum_{i,j} \left(
[V_i^{(1)}\otimes W_j^{(1)}] \right) \otimes \left( [V_i^{(2)}\otimes
W_j^{(2)}] \right).  \eean
\end{lem}

\begin{lem}\label{coass}
Let $N=N_1+N_2+N_3$. Then we have
\be 
\left(\Delta^{N_1+N_2}_{N_1,N_2}\otimes \on{id}\right) \circ
\Delta^N_{N_1+N_2,N_3} =
\left(\on{id}\otimes \Delta^{N_2+N_3}_{N_2,N_3} \right)\circ
\Delta^N_{N_1,N_2+N_3}.
\ee
\end{lem}

Now we compute $\Delta^N_{N_1,N_2}$ explicitly on the generators
$t_{i,n}$ of $\Rep_a U_q \agln$ (i.e., the fundamental representations).
Using formula \eqref{ti}, we obtain

\begin{lem} Let $N_1,N_2\geq i$. Then
\bean\label{delta fund}
\Delta^N_{N_1,N_2}(t_{i,n})=\sum_{j=0}^i t_{j,n}\otimes t_{i-j,n-j}.
\eean
\end{lem}
Therefore, the result does not depend on $N,N_1,N_2$ when they are
large enough. \lemref{compatibility} implies that $\Delta^N_{N_1,N_2}$
of a tensor product of fundamental representations $t_{i,n}, i\leq
N_1,N_2$ can be written as a linear combination of terms $W^{(1)}
\otimes W^{(2)}$, where each $W^{(j)}$ is a tensor product of the
fundamental representations. Moreover, the result does not depend on
$N,N_1,N_2$ when they are large. Therefore we obtain a well-defined
map
$$
\Delta: \Rep_a U_q \agli \to \Rep_a U_q \agli \otimes \Rep_a U_q
\agli.
$$
In more concrete terms the map $\Delta$ may be described as
follows. We have
\bean
\label{deltati}
\Delta(t_{i,n})=\sum_{j=0}^i t_{j,n}\otimes t_{i-j,n-j}. \eean Since
$\Rep_a U_q \agli$ is the free polynomial algebra in $t_{i,j}$, the
following formula then allows us to extend $\Delta$ inductively to all
other elements of $\Rep_a U_q \agli$ (cf. \eqref{delta of product}):
$$
\Delta([V\otimes W])=\sum_{i,j} \left(
[V_i^{(1)}\otimes W_j^{(1)}] \right) \otimes \left( [V_i^{(2)}\otimes
W_j^{(2)}] \right),
$$
where
$$
\Delta([V]) = \sum_i [V_i^{(1)}]\otimes [V_i^{(2)}], \qquad
\Delta([W]) = \sum_i [W_i^{(1)}]\otimes [W_i^{(2)}].
$$

Using the decomposition \eqref{decomp for a}, we also obtain a map:
$$
\Delta: \Rep\, U_q \agli \to \Rep\, U_q \agli \otimes \Rep\, U_q
\agli.
$$

\begin{thm}\label{hopf str} 
There is a unique antipode $S$, such that $(\Rep\,
U_q\agli,\otimes,\Delta,S)$ is a Hopf algebra, and $(\Rep_a
U_q\agli,\otimes,\Delta,S)$ is its Hopf subalgebra for any $a \in
\C^\times$. Moreover, $\Rep_aU_q\agln$ ($N \in \Z_{>0}$) is a Hopf
subalgebra of $\Rep_aU_q\agli$, and the Hopf algebra $\Rep_aU_q\agli$
is an inductive limit of these subalgebras.

Formula \eqref{decomp for a} gives us a decomposition of
the Hopf algebra $\Rep \, U_q\agli$ into the tensor product of its
Hopf subalgebras $\Rep_a U_q\agli, a \in \C^\times/q^{2\Z}$.
\end{thm}

\begin{proof}
We have already proved the associativity of $\otimes$.
Coassociativity of $\Delta$ follows from Lemma \ref{coass} and
stabilization of $\Delta_{N_1,N_2}^N$. The compatibility of $\otimes$
and $\Delta$ follows from Lemma \ref{compatibility} and stabilization
of $\Delta_{N_1,N_2}^N$.

Thus, we obtain the structure of a bialgebra on
$\Rep_aU_q\agli$. However, we will show in Lemma \ref{repa to
uminus} below that this bialgebra is isomorphic to the bialgebra of a
Hopf algebra. Therefore it admits a unique antipode (recall that if an
antipode exists, then it is unique). This proves the first assertion
of the theorem.

It is clear from the definition of multiplication on
$\Rep_aU_q\agli$ that $\Rep_aU_q\agln$ is its
subalgebra. Furthermore, formula \eqref{deltati} implies that
$\Rep_aU_q\agln$ is stable under comultiplication. It is also
stable under the antipode, as we will see below. Therefore we obtain
the second assertion.
\end{proof}

Note that Hopf algebra $\Rep_aU_q\agli$ is commutative but
not cocommutative. 

\subsection{The Hopf algebra $\Rep\,U_q\gli$} The definition of the Hopf
algebra $\Rep\,U_q\gli$ is parallel to that of
$\Rep\,U_q\agli$. Explicitly, we have 
\be
\Rep\,U_q\gli=\Z[r_i]_{i\in\Z_{>0}},\qquad \Delta(r_i)=\sum_{j=0}^i
r_j\otimes r_{i-j}.
\ee
Here $r_i$ ($i\in\Z_{>0}$)
are the classes of the $i$th fundamental modules and we set $r_0=1$.
In particular, the Hopf algebra $\Rep\,U_q\gli$ is commutative and
cocommutative.

\section{A Hopf algebra structure on $\Rep_aU_q\agli$}
\label{ind lim structure generic}

In this section we identify the Hopf algebra $\Rep_aU_q\agli$ with the
Hopf algebra of functions on the group $\wSLi$ of lower unipotent
infinite matrices. Then we describe the corresponding action of the
Lie algebra $\sw_\infty$ on $\Rep_aU_q\agli$ in representation
theoretic terms. Finally, we identify the subspace of $\Rep_aU_q\agli$
spanned by the evaluation representations taken at a fixed set of
points of $aq^{2\Z}$ with an integrable representation of
$\sw_\infty$.

\subsection{Hopf algebras $\Z[\wSLi]$ and $\Z[\wSLin]$}\label{Minfty}

For $N \in \Z_{>0}$, let $\wSLin$ be the proalgebraic group (over
$\Z$) of matrices $M = (M_{i,j})_{i,j \in \Z}$ satisfying the
conditions: $M_{i,i} = 1$ for all $i$, and $M_{i,j} = 0$ whenever
$i<j$ or $i>j+N$. We define a product of two such matrices as their
usual product, in which all $(i,j)$ entries with $i>j+N$ are set to be
equal to $0$. It is clear that this operation makes $\wSLin$ into a
prounipotent proalgebraic group. Furthermore, we have natural
surjective group homomorphisms $\wt{SL}^{(N+1),-}_\infty
\twoheadrightarrow \wSLin$. Let
$$
\wSLi = \underset{\longleftarrow}{\lim} \; \wSLin.
$$
This is also a prounipotent proalgebraic group, which consists of
all matrices $M = (M_{i,j})_{i,j \in \Z}$, such that $M_{i,i} = 1$ for
all $i$, and $M_{i,j} = 0$ for all $i<j$.

The algebra $\Z[\wSLi]$ of regular functions on $\wSLi$ is by
definition the algebra of polynomials in commuting variables
$M_{i,j}$, where $i,j\in\Z$, $i>j$.  The algebra $\Z[\wSLin]$ of
functions on $\wSLin$ is the algebra of polynomials in commuting
variables $M_{i,j}$, where $i,j\in\Z$ and $j< i\leq j+N$. Both are
Hopf algebras, with the comultiplication $\Delta$ given on the
generators by the formula \bean
\label{delta of M} \Delta
(M_{i,j})=\sum_{k=j}^i M_{i,k}\otimes M_{k,j}. \eean The antipode $S$
is an anti-automorphism given on the generators by the formula
\bean\label{S of M} S(M)=M^{-1}, \eean where $M$ is the matrix
$M=(M_{i,j})_{i,j\in\Z}$. Furthermore, the Hopf algebra $\Z[\wSLi]$ is
the inductive limit of its Hopf subalgebras $\Z[\wSLin]$, $N \in
\Z_{>0}$.

Recall the isomorphism $\Rep_a U_q \agli \simeq
\Z[t_{i,n}]_{i\in\Z_{>0}}^{n\in\Z}$, $\Rep_a U_q \agln \simeq
\Z[t_{i,n}]_{i=1,\ldots,N}^{n\in\Z}$.  Define homomorphisms of
polynomial algebras 
\bean\label{repa = dual sli} \Z[\wSLi] &\to&
\Rep_a U_q \agli,\qquad M_{i+j,j}\mapsto t_{i,j}\qquad
(j\in\Z,\;i\in \Z_{>0}),\\  
\Z[\wSLin] &\to&
\Rep_a U_q \agln,\qquad M_{i+j,j}\mapsto t_{i,j}\qquad
(j\in\Z,\;i=1,\dots,N).\label{repa = dual sln}
\eean 
Thus, in the case $\Rep_a U_q\agli$ we have 
\be M= \left(
\begin{matrix}
&&\dots&&\dots&\\
&1&0&0&0&\\
&t_{1,i-1}& 1 & 0 & 0 &\\
\dots&t_{2,i} & t_{1,i}&1&0&\dots\\
&t_{3,i+1}&t_{2,i+1}&t_{1,i+1}&1& \\
&&\dots&&\dots&
\end{matrix}
\right), \ee where $t_{i,j}\in \Rep_a U_q\agli$ corresponds to the
$i$-th fundamental representation evaluated at $aq^{2j}$. In the case
$\Rep_a U_q\agln$ we have the same matrix $M$ with all entries below
the $N$-th diagonal replaced by zero.

\begin{thm}\label{repa to uminus}
The maps \Ref{repa = dual sli} and \Ref{repa = dual sln} are
isomorphisms of Hopf algebras. Moreover, we have a commutative diagram
$$
\begin{CD}
\Z[\wSLi] @>{\sim}>> \Rep_a U_q\agli \\
@AAA @AAA \\
\Z[\wSLin] @>{\sim}>> \Rep_a U_q\agln \\
\end{CD}
$$
\end{thm}

\begin{proof}
We only need to check the compatibility of the maps \Ref{repa = dual
sli}, \Ref{repa = dual sln} with the comultiplication. 
This follows immediately by comparison of
formulas \eqref{delta fund} and \eqref{delta of M}.
\end{proof}

\subsection{The Hopf algebra $U_\Z{\mathfrak s}{\mathfrak l}_\infty^-$}
\label{dlinnaya}

In this section we describe the Hopf algebra $U_\Z{\mathfrak
s}{\mathfrak l}_\infty^-$ and identify it as restricted dual to
$\Z[\wSLi]$ (and hence to $\Rep_a U_q\agli$).

Let $\sw_\infty^-$ be the Lie algebra of strictly lower triangular
matrices of infinite size (with rows and columns both labeled by $\Z$)
with finitely many nonzero complex entries.  Its completion
$\wt{\sw}_\infty^-$ consisting of arbitrary strictly lower triangular
matrices is the Lie algebra of the complex version of the group
$\wSLi$.

The Lie algebra $\sw_\infty^-$ is generated by $f_i$
($i\in\Z$), subject to the Serre relations: \bean\label{serre}
[f_i,[f_i,f_{i\pm1}]]=0 \quad (i\in\Z),\qquad \;[f_i,f_j]=0 \quad
(|i-j|>1).
\eean 
Define the root vectors $f_{i,j}$ ($i,j\in\Z$, $i>j$) by the
formula \bean\label{root vectors}
 f_{i,j}= [\dots
[[f_{i-1},f_{i-2}],f_{i-3}],\;\dots\;,f_{j}].  \eean In particular,
$f_{i+1,i}=f_i$. The element $f_{i,j}$ corresponds to the matrix with
the $(i,j)$th entry $1$ and all other entries $0$. The root vectors
form a basis of $\sw_\infty^-$, and we have the commutation relations
\bean\label{commutation in agli} [f_{i,j},f_{k,l}]=\delta_{j,k}
f_{i,l}-\delta_{i,l}f_{k,j}.  \eean

The algebra $U{\mathfrak s}{\mathfrak l}_\infty^-$ is by definition
the universal enveloping algebra of Lie algebra $\sw_\infty^-$. Thus,
it is an associative algebra over $\C$ generated by $f_i$ ($i\in\Z$),
satisfying the relations \eqref{serre}. The algebra $U{\mathfrak
s}{\mathfrak l}_\infty^-$ is a cocommutative Hopf algebra with the
comultiplication $\Delta$ and the antipode $S$ given on the generators
by \be \Delta(f_i)=f_i\otimes 1+ 1\otimes f_i,\qquad S(f_i)=-f_i.  \ee

Note that the root vectors are primitive, i.e.,
\bean\label{roots are primitive}
\Delta(f_{i,j})=f_{i,j}\otimes 1 + 1\otimes f_{i,j},\qquad S(f_{i,j})
= - f_{i,j}.
\eean

It is well known that the algebra $U{\mathfrak s}{\mathfrak
l}_\infty^-$ has an integral form generated by the divided powers
$f_i^{(n)}=f_i^n/n!$.  It is a Hopf algebra over $\Z$ which we denote
by $U_\Z{\mathfrak s}{\mathfrak l}_\infty^-$.  Note that the divided
powers of the root vectors $f_{i,j}^{(n)}$ belong to $U_\Z{\mathfrak
s}{\mathfrak l}_\infty^-$.

We fix the following order on the set of root vectors: $f_{i,j} \succ
f_{k,r}$, if $j>r$ or $j=r$ and $i>k$. We have a PBW $\Z$--basis in
$U_\Z{\mathfrak s}{\mathfrak l}_\infty^-$, which consists of the
ordered products of divided powers of the root vectors of the form
\bean\label{pbv vector} f_{i_1,j_1}^{(n_1)} f_{i_2,j_2}^{(n_2)}\dots
f_{i_m,j_m}^{(n_m)}, \eean where $f_{i_s,j_s} \succ
f_{i_{s-1},j_{s-1}}$ $(s=2,\dots,m)$.

Now recall the Hopf algebra $\Z[\wSLi]$ defined in the previous
section.  Fix the order on the generators $M_{i,j}$ similar to the one
in the case of $U_\Z{\mathfrak s}{\mathfrak l}_\infty^-$. Namely,
$M_{i,j} \succ M_{k,r}$, if $j>r$ or $j=r$ and $i>k$. We have a
$Z$--basis in $\Z[\wSLi]$, which consists of monomials of the form
\bean\label{pbv dual vector} M_{k_1,r_1}^{p_1} M_{k_2,r_2}^{p_2}\dots
M_{k_n,r_n}^{p_n}, \eean where $M_{k_s,r_s} \succ M_{k_{s-1},r_{s-1}}$
$(s=2,\dots,n)$.

There is the following nondegenerate pairing
$$\langle,\rangle:\;U{\mathfrak s}{\mathfrak l}_\infty^-\otimes
\C[\wSLi] \to \C.$$ Each element of $U{\mathfrak s}{\mathfrak
l}_\infty^-$ defines a differential operator on the group $\wSLi$. We
then set
$$
\langle P,f \rangle = (P \cdot f)(1),
$$ where $1$ denotes the unit element of $\wSLi$. 

It is straightforward to check that this pairing is compatible with
the integral forms and that the $\Z$--bases described above are dual
to each other with respect to this pairing: \bean\label{pairing def}
\left\langle f_{i_1,j_1}^{(s_1)}\dots
f_{i_m,j_m}^{(s_m)},M_{k_1,r_1}^{p_1} \dots M_{k_n,r_n}^{p_n}
\right\rangle =1, \eean if $m=n$ and $i_t=k_t$, $j_t=r_t$, $s_t=p_t$
($t=1,\dots,m$), and the pairing of basis elements is zero otherwise.

\begin{lem}    \label{two dual}
The pairing $\langle,\rangle$ is a nondegenerate Hopf algebra
pairing, i.e., for all $f,g\in U_\Z{\mathfrak s}{\mathfrak
l}_\infty^-$, $T,K \in \Z[\wSLi]$, we have \bea \langle
f,TK\rangle&=&\langle \Delta (f), T\otimes K\rangle,\quad \langle
gf,T\rangle=\langle g\otimes f,\Delta (T) \rangle, \quad \langle S(f),
T \rangle= \langle f, S(T)\rangle. \eea
\end{lem}
\begin{proof}
The verification of the first equality on the basis elements is
straightforward using \eqref{roots are primitive}.

If the second equality is true for $f=h_1$ and $f=h_2$ (and for all
$g\in U_\Z{\mathfrak s}{\mathfrak l}_\infty^-$) then it is also true for
$f=h_1 h_2$. Therefore it suffices to check the second equality in
the case when $f$ is one of the generators $f_i$. This is done by a direct
computation using \eqref{commutation in agli}.

The third equality now follows automatically, because if it exists, the
antipode is unique.
\end{proof}

We call two vector spaces $V,W$ {\em restricted dual} of each other,
if we are given a pairing $V\otimes W\to\C$ and linear bases
$\{v_i\}_{i\in I}$ of $V$ and $\{w_i\}_{i\in I}$ of $W$, labeled by
the same set $I$, on which the pairing has the form $v_i\otimes
w_j\mapsto \delta_{i,j}$. In other words, the pairing is
non-degenerate and the bases $\{v_i\}$, $\{w_i\}$ are dual to each
other.

In this paper we always work with vector spaces which have
graded components of finite rank, and all pairings respect the
gradations. Then by the restricted dual space of $V=\oplus_i V_i$,
where $\dim V_i<\infty$, we understand $V^*:=\oplus_i V_i^*$. For
instance, we have a $\Z$--gradation on $U{\mathfrak s}{\mathfrak
l}_\infty^-$ defined by the formula $\deg f_i = 1$ for all $i
\in \Z$, and a $\Z$--gradation on $\C[\wSLi]$ defined
by the formula $\deg M_{i,j} = i-j$ for all $i,j \in \Z$, $i>j$. The
above pairing respects these gradations. Therefore \lemref{two dual}
implies that the Hopf algebras $U{\mathfrak s}{\mathfrak
l}_\infty^-$ and $\C[\wSLi]$ are restricted dual of each other.

\medskip

\begin{remark}
The topologically dual Hopf algebra of $U{\mathfrak s}{\mathfrak
l}_\infty^-$ (in the sense of \cite{CP}, Section 4.1D), is larger than
$\C[\wSLi]$.  The reason is that the algebra $U{\mathfrak
s}{\mathfrak l}_\infty^-$ contains many primitive elements $f_{i,j}$, so
that the elements $\exp(M_{i,j})=\exp(f^*_{i,j})$ are group-like and
therefore belong to the topological dual.\qed
\end{remark}

\medskip

Now suppose that $V$ and $W$ have integral forms $V_\Z$ and $W_\Z$,
which are free $\Z$--modules with $\Z$--bases dual to each
other under the above pairing. Then we will say that $V_\Z$ and $W_\Z$
are (restricted) dual to each other.

In our case $\C[\wSLi]$ has the integral form $\Z[\wSLi]$. The
corresponding dual integral form of $U{\mathfrak s}{\mathfrak
l}_\infty^-$ is its $\Z$--subalgebra $U_{\Z}{\mathfrak s}{\mathfrak
l}_\infty^-$. Formula \eqref{pairing def} shows that the monomial
$\Z$--basis of $\Z[\wSLi]$ is dual to the PBW $\Z$--basis of divided
powers in $U_{\Z}{\mathfrak s}{\mathfrak l}_\infty^-$. Therefore the
Hopf algebras $\Z[\wSLi]$ and $U_{\Z}{\mathfrak s}{\mathfrak
l}_\infty^-$ are dual to each other.

\subsection{Restriction operators}\label{restrictions section} 

In this section we compare the restriction operators $\res_{\bs m}$ on
$\Rep_a U_q\agli$ with the operators of right multiplication by the
products of the divided powers of the generators $f_i$ on $\Z[\wSLi]$.

\begin{lem}\label{f is dual to res} 
For all $g\in U\sli^-$, $T\in\Rep_aU_q\agli$, we have
\bean\label{dual pairing}
\langle gf_i,T\rangle = \langle g, \res_i\;T\rangle.
\eean
Moreover, if ${\bs m}$ has components $n_1>\dots>n_s$, and $n_i$
occurs $k_i$ times,
then
\bean\label{dual pairing general}
\langle gf_{\bs m},T\rangle = \langle g,
\res_{\bs m}T\rangle,
\eean
where
$
f_{\bs m}=f_{n_1}^{(k_1)}\dots f_{n_s}^{(k_s)}$. 
\end{lem}
\begin{proof}
Note that \Ref{dual pairing} is a special case of \Ref{dual pairing
general} (actually, the two formulas are equivalent, but we prefer to
prove \Ref{dual pairing general}, since we will need it in Section
\ref{restrictions section root}
in a different context).

First we prove \Ref{dual pairing general} when $T=t_{j,n}$. If $s=1$,
we obtain
\bea
\langle gf_i,t_{j,n}\rangle=\langle g\otimes
f_i,\Delta(t_{j,n})\rangle=
\langle g\otimes f_i, \sum_{k=0}^j t_{j-k,n}\otimes t_{k,n-j+k} \rangle
\\ 
=\sum_{k=0}^j\langle g,t_{j-k,n}\rangle\langle f_i,t_{k,n-j+k}\rangle=
\delta_{i, n-j+1} \langle g, t_{j-1,n}\rangle=\langle g, \res_i\;
t_{j,n}\rangle.
\eea
For $s\geq 2$ we have
\be
\langle g,\res_{\bs m} t_{j,n}\rangle = 0 =\langle g f_{\bs
m},t_{j,n}\rangle.
\ee

Now suppose that \Ref{dual pairing general} is proved for $T=T_1$ and
$T=T_2$. Then for $T=T_1T_2$, we obtain
\begin{align*}
&\langle g,\res_{\bs m}(T_1T_2)\rangle \\ &=\langle g, \sum_{{\bs
m}={\bs m_1}\sqcup{\bs m_2}}(\res_{\bs m_1} T_1)(\res_{\bs
m_2}T_2)\rangle =\langle \Delta g, \sum_{{\bs m}={\bs m_1}\sqcup{\bs
m_2}}(\res_{\bs m_1} T_1)\otimes (\res_{\bs m_2}T_2)\rangle \\ &=\langle
\Delta(g) \sum_{{\bs m}={\bs m_1}\sqcup{\bs m_2}}f_{\bs m_1}\otimes
f_{\bs m_2},T_1\otimes T_2\rangle =\langle \Delta (gf_{\bs
m}),T_1\otimes T_2\rangle=\langle gf_{\bs m},T_1T_2\rangle.
\end{align*}
This proves \Ref{dual pairing general} in general.
\end{proof}

The algebra $U_\Z {\sw}_\infty^-$ acts naturally on $\Z[\wSLi]$. The
above lemma shows that under the isomorphism $\Z[\wSLi] \simeq \Rep_a
U_q\agli$ the generators $f^{(n)}_i$ of $U_\Z{\sw}_\infty^-$ act as the
operators $\res_{(i,\ldots,i)}$, where $i$ is repeated $n$
times. In particular, we obtain the following 

\begin{thm}    \label{action of sl}
The restriction operators $\res_i$ satisfy the Serre relations, \bea
\res_i^2\,\res_{i+1}-2\res_i\,
\res_{i+1}\,\res_i+\res_{i+1}\,\res_i^2=0&\qquad& (i\in\Z),\\
\res_i\,\res_j=\res_j\,\res_i&\qquad& (|i-j|>1).  \eea The space
$\Rep_a U_q\agli \otimes \C$, equipped with the action of the
operators $\res_i$ is a left $U\mathfrak s\mathfrak
l_\infty^-$--module, which is (restricted) dual to the free right
$U\mathfrak s\mathfrak l_\infty^-$--module with one generator.
\end{thm}

The categories of the right and left $U\mathfrak s\mathfrak
l_\infty^-$ modules are equivalent due to the existence of the
antipode. Hence we will identify $\Rep_aU_q\agli \otimes \C$ with the
restricted dual of a free left $U\mathfrak s\mathfrak
l_\infty^-$--module with one generator.

\subsection{Integrable representations of $\sw_\infty$}
\label{int section}

We recall some facts about the Lie algebra $\sw_\infty$ (see \cite{K}
for details). It has generators $e_i$, $h_i$, $f_i$ ($i \in \Z$), satisfying
the standard relations. In particular, the Lie subalgebra of
$\sw_\infty$ generated by $f_i$ ($i \in \Z$) is the Lie algebra
$\sw_\infty$ introduced in \secref{dlinnaya}.

The weight lattice of $\sw_\infty$ is spanned by the fundamental
weights $\varpi_i$ ($i \in \Z$). It contains the root lattice spanned by
$\al_i = 2 \varpi_i - \varpi_{i-1} - \varpi_{i+1}$ ($i \in \Z$). For
weights $\la,\mu$ we write $\la \geq \mu$ if $\la = \mu + \sum_i n_i
\al_i$, where all $n_i \in \Z_{\geq 0}$.

Introduce the category ${\mc O}$ of $\sw_\infty$--modules. Its objects
are $\sw_\infty$--modules $L$ on which the operators $h_i$ ($i \in \Z$)
act diagonally, and the weights occurring in $L$ are less than or equal
to some fixed weight $\la$. The morphisms in ${\mc O}$ are
$\sw_\infty$--homomorphisms. A representation $L$ of $\sw_\infty$ is
called integrable if it belongs to ${\mc O}$ and all generators $e_i$
and $f_i$ act on $L$ locally nilpotently. It is called a highest
weight representation if it is generated by a vector $v$, such that
$e_i \cdot v = 0$ for all $i \in \Z$. 

The category of integrable
representations of $\sw_\infty$ is semisimple. Irreducible integrable
representations of $\sw_\infty$ are the irreducible highest weight
representations with the highest weights $\bs\nu = \sum_i \nu_i
\varpi_i$, where $\nu_i \in \Z_{\geq 0}$ for all $i \in \Z$, and
$\nu_i=0$ for all but finitely many $i$'s. Such weights are called
dominant integral weights.

The irreducible representation corresponding to a dominant integral
weight $\bs\nu = \sum_i \nu_i \varpi_i$ is generated by a vector
$v_{\bs\nu}$ such that
$$
e_i v_{\bs\nu} = 0, \qquad h_i v_{\bs\nu} = \nu_i v_{\bs\nu}, \qquad
f_i^{\nu_i+1} v_{\bs\nu} = 0 \qquad (i \in \Z).
$$
It is denoted by $L(\bs\nu)$. The number $k =
\sum_i \nu_i$ is called the level of $L(\bs\nu)$. Since $\sw_\infty$ is
a Lie subalgebra of $\sw_\infty$, its universal enveloping algebra
$U \sw_\infty^-$ acts on any $\sw_\infty$--module. In particular, we
have a canonical surjective homomorphism of $U\sw_\infty^-$--modules
\bean
\label{from U to L} p_{\bs\nu}: U \sw_\infty^- \to
L(\bs\nu) \qquad g \mapsto g \cdot v_{\bs\nu}. 
\eean

Level one representations $L(\varpi_n)$ ($n \in \Z$) have the following
explicit realization, due to Hayashi \cite{Ha}. Consider the {\em Fock
space} $\bigwedge$ with the basis of vectors $|{\bs \la} \rangle$ labeled
by all partitions ${\bs \la}$. In order to describe the action of the
generators $e_i$ and $f_i$ in this basis, we need to introduce some
terminology. We say that a box $B_{j,r}$ can be removed from (resp.,
added to) the Young diagram ${\bs \la}$ if removing this box from
this diagram (resp., adding this box to this diagram) one obtains
another Young diagram. We denote the partition corresponding to the
new Young diagram by $\bs\la\setminus(i,j)$ (resp., $\bs\la \cup (i,j)$).

Fix $n\in\Z$. Define the operators $\wt{e}_m$, $\wt{f}_m$ and
$\wt{h}_m$ on $\bigwedge$ by the following rules. Set $\wt{e}_m |{\bs \la}
\rangle$ equal $-|\bs\la\setminus (i,j)\rangle$ if the Young diagram
of $\la$ contains a removable box $B_{i,j}$ and $m=n+j-i$, and equal
$0$ otherwise. Set $\wt{f}_m |{\bs \la} \rangle$ equal $-|\bs \la\cup
(i,j)\rangle$ if the box $B_{i,j}$ can be added to the Young diagram
of ${\bs \la}$ and $m=n+j-i$, and equal $0$ otherwise. Finally, set $\wt{h}_m
|{\bs \la} \rangle$ equal $|{\bs \la} \rangle$, if a box $B_{i,j}$
with $j-i=m-n$ can be added to the Young diagram of ${\bs \la}$,
equal $-|{\bs \la} \rangle$, if a box $B_{i,j}$ with $j-i=m-n$ can be
removed from the Young diagram of ${\bs \la}$, and equal $0$
otherwise.

\begin{lem}[\cite{Ha}]    \label{Hayashi}
The operators $\wt{e}_m$, $\wt{f}_m$ and $\wt{h}_m$ define a
representation of $\sw_\infty$ on $\bigwedge$ denoted by
${\bigwedge}(\varpi_n)$. 
Moreover, $\bigwedge(\varpi_n) \simeq L(\varpi_n)$.
\end{lem}

\subsection{Evaluation modules and integrable modules}\label{int mod section}

Recall that we have the evaluation homomorphism, $ev_{a}: U_q \agln
\to U_q \gln$. Fix $n\in\Z$. By Lemma \ref{evaluation character
lemma}, $ev_{aq^{2n}}$ induces an embedding
\bean\label{gl in agl}
\iota^{(N)}(\varpi_n):\;\Rep\, U_q\gln \to
\Rep_a U_q\agln,\qquad [V_{\bs\la}]\mapsto
[V_{\bs\la}(aq^{2n})],
\eean 
where $V_{\bs\la}$ denotes the
irreducible representation of $U_q\gln$ of highest weight $\bs\la$. By
Lemma \ref{irrep to fund stable}, the embeddings \eqref{gl in agl}
stabilize and hence induce an embedding
\bean\label{gli in agli} \iota(\varpi_n):\;\Rep\, U_q\gli \to
\Rep_a U_q\agli. \eean

Recall from \thmref{repa to uminus} and \lemref{two dual} that we have
a pairing of Hopf algebras
$$
\langle,\rangle: U_\Z\mathfrak s\mathfrak l_\infty^- \otimes \Rep_a
U_q\agli \to \Z
$$
and that $\Rep_a U_q\agli$ is naturally a $U_\Z\mathfrak
s\mathfrak l_\infty^-$--module, which is restricted dual to
$U_\Z\mathfrak s\mathfrak l_\infty^-$.

The following lemma is obtained from Lemma \ref{evaluation character
lemma}.

\begin{lem}    \label{iotan}
The $\Z$--submodule $\iota(\varpi_n)(\Rep\, U_q\gli) \subset \Rep_a
U_q\agli$ is a right comodule of $\Rep_a U_q\agli$ and a
$U_\Z\mathfrak s\mathfrak l_\infty^-$--submodule of $\Rep_a U_q\agli$.
\end{lem}

We now identify the $U\mathfrak s\mathfrak l_\infty^-$--module
$\iota(\varpi_n)(\Rep\, U_q\gli) \otimes \C$ with the restricted dual
of an integrable level 1 representation.

Define the non-degenerate pairing 
\bean\label{pairing of level 1}
\langle,\rangle_n:\;  \bigwedge(\varpi_n) \otimes \left(\iota(\varpi_n)(\Rep\,
U_q\gli) \otimes \C\right)\to \C,
\eean 
by declaring the bases $\{ | {\bs \la}
\rangle\}$ and $\{ [V_{\bs \la}(aq^{2n})] \}$ to be dual to each
other.

\begin{prop}    \label{level one}
The $U\mathfrak s\mathfrak l_\infty^-$--module $\iota(\varpi_n)(\Rep\,
U_q\gli) \otimes \C$ is restricted dual to $L(\varpi_n)$ with respect
to the pairing \eqref{pairing of level 1}.  Furthermore, the embedding
$$\iota(\varpi_n):\;\Rep\, U_q\gli \otimes \C \hookrightarrow\Rep_a
U_q\agli \otimes \C$$ is dual to the surjective $U\mathfrak s\mathfrak
l_\infty^-$--homomorphism $$p_{\varpi_n}: U\mathfrak s\mathfrak
l_\infty^- \twoheadrightarrow L(\varpi_n).$$
\end{prop}

\begin{proof}
We consider the case $n=0$; all other cases are obtained by a shift of
all indices.

First, we check the formula \bean \langle f_i x,y \rangle_0 + \langle
x,\res_i y \rangle_0 = 0.\label{module-module} \eean For that, we
compute explicitly the action of operators $\res_m$ on
$[V_{\bs\la}(a)]$. Recall that the content of the box $B_{i,j}$ is
defined to be $j-i$. From the formula for the $q$--character of
$V_{\bs \la}(a)$ given in Lemma \ref{evaluation character lemma} we
obtain that if the partition $\bs\la$ has a removable box of content
$m$, then \be \res_m V_{\bs\la}(a) = V_{\bs\la\setminus (i,j)}(a). \ee 
 If the partition $\bs\la$ has no
removable boxes of content $m$, then $\res_m V_{\bs\la}(a) = 0$. 
Comparing this with the action of $f_m$ in $\bigwedge(\varpi_0)$ (given by
$\wt{f}_m$), we obtain formula \eqref{module-module}. This proves the
first statement of the proposition.

Now we have the diagram of $U\mathfrak s\mathfrak l_\infty^-$--modules
and homomorphisms $$\begin{CD} \Rep\,U_q\gli \otimes \C
@>{\iota(\varpi_0)}>> \Rep_a U_q\agli \otimes \C \\
@V{\langle,\rangle_0}VV @VV{\langle,\rangle}V\\
\bigwedge(\varpi_0)^* @>p_{\varpi_0}^*>> (U\mathfrak s\mathfrak
l_\infty^-)^*
\end{CD}$$
To complete the proof, we need to show that this diagram is
commutative. Note that the vertical arrows are isomorphisms. Consider
the homomorphism $\bigwedge(\varpi_0)^* \to (U\mathfrak s\mathfrak
l_\infty^-)^*$, obtained by composing these isomorphisms with
$\iota(\varpi_0)$. Its dual sends $1 \in U\mathfrak
s\mathfrak l_\infty^-$ to the vector $|0 \rangle \in
\bigwedge(\varpi_0)$. Such a homomorphism is unique, and hence
coincides with $p_{\varpi_0}^*$. Therefore the above diagram is
commutative.
\end{proof}

\propref{level one} implies that as a $\Rep_a U_q\agli$--comodule,
$\iota(\varpi_n)(\Rep\, U_q\gli)$ is an integral form of the
restricted dual to the $U\mathfrak s\mathfrak l_\infty^-$--module
$L(\varpi_n)$, which is preserved by $U_\Z\sw_\infty^-$.

Now we describe the integrable modules of higher levels $k$.
Given a dominant integral weight $\bs\nu = \sum_i \nu_i \varpi_i$, set
$\bigwedge(\bs\nu) = \otimes_i (\bigwedge({\varpi_i})^{\otimes \nu_i})$. This
is an integrable module of $\sw_\infty$ of level $k= \sum_i
\nu_i$. Denote the tensor product of the highest weight vectors of
$\bigwedge({\varpi_i})$ by $|0\rangle_{\bs\nu}$. This vector has weight
${\bs\nu}$, which is the highest among all weights occurring in
$\bigwedge(\bs\nu)$, and satisfies:
$$
e_i |0\rangle_{\bs\nu} = 0, \qquad f_i^{\nu_i+1} |0\rangle_{\bs\nu} =
0 \qquad (i\in\Z).
$$
Hence we have a homomorphism of $\sw_\infty$--modules
$\gamma_{\bs\nu}: L(\bs\nu) \to \bigwedge(\bs\nu)$ sending
$v_{\bs\nu}$ to $|0\rangle_{\bs\nu}$. Since $L(\bs\nu)$ is
irreducible, this homomorphism is injective. The semi-simplicity of
the category of integrable $\sw_\infty$--modules implies the
following

\begin{lem}\label{fock as cyclic submodule}
The irreducible module $L(\bs\nu)$ is a direct summand in
$\bigwedge(\bs\nu)$, which is equal to the cyclic $U\mathfrak s\mathfrak
l_\infty^-$--submodule of $\bigwedge(\bs\nu)$ generated by
$|0\rangle_{\bs\nu}$.
\end{lem}

Denote by $\wt{p}_\nu$ the homomorphism $U\mathfrak s\mathfrak
l_\infty^- \to \bigwedge(\bs\nu)$ obtained by composing $p_{\bs\nu}$ given by
formula \eqref{from U to L} and $\gamma_{\bs\nu}$. \lemref{fock as cyclic
submodule} means that the image of $\wt{p}_{\bs\nu}$ is a direct summand of
$\bigwedge(\bs\nu)$ isomorphic to $L(\bs\nu)$.

Write $\bs\nu$ in the form $\sum_{j=1}^k \varpi_{s_j}$. Using
\propref{level one}, we identify the restricted dual
$\bigwedge(\bs\nu)^*$ with $(\Rep\, U_q\gli)^{\otimes k} \otimes \C$
as $U \sw_\infty^-$--modules. Denote by $\langle,\rangle_{\bs\nu}$ the
pairing $\bigwedge(\bs\nu) \otimes (\Rep\, U_q\gli)^{\otimes k} \to
\C$.

According to formula \eqref{res tensor}, the map \bean\label{many gl
to agl} \iota({\bs\nu}):\; (\Rep\, U_q\gli)^{\otimes k} \otimes \C
&\to& \Rep_a U_q\agli \otimes \C,\\
{[V_{\bs\la_{1}}]\otimes\dots\otimes [V_{\bs\la_{k}}]}
&\mapsto& [V_{\bs\la_{1}}(aq^{2s_1})\otimes\dots\otimes
V_{\bs\la_{k}}(aq^{2s_k})]\notag \eean is a homomorphism of $U
\sw_\infty^-$--modules. The image of $\iota({\bs\nu})$ is spanned by all
possible tensor products of $k$ evaluation modules, evaluated at the
points $aq^{2s_1},\dots,aq^{2s_k}$. We denote the image of
$\iota({\bs\nu})$ by $\Rep_a^{\bs\nu} U_q \gli$.

We have a diagram of $U \sw_\infty^-$--modules and $U
\sw_\infty^-$--homomorphisms
$$\begin{CD} (\Rep\, U_q\gli)^{\otimes k} \otimes \C
@>{\iota(\bs\nu)}>> \Rep_a U_q\agli \otimes \C \\
@V{\langle,\rangle_{\bs\nu}}VV @VV{\langle,\rangle}V\\
\bigwedge(\bs\nu)^* @>\wt{p}_{\bs\nu}^*>> (U\mathfrak
s\mathfrak l_\infty^-)^*
\end{CD}$$

\begin{prop}\label{many gl to agl lemma}
The above diagram is commutative. The image of the map $\iota(\bs\nu)$ is
a $U\mathfrak s\mathfrak l_\infty^-$--submodule of $\Rep_a U_q\agli
\otimes \C$ isomorphic to the restricted dual of $L(\bs\nu)$.
\end{prop}

\begin{proof}
The first statement follows from the uniqueness of the homomorphism
$U\mathfrak s\mathfrak l_\infty^- \to \bigwedge(\bs\nu)$ sending $1$ to
$|0\rangle_{\bs\nu}$, as in the proof of \propref{level one}. The second
statement follows from \lemref{fock as cyclic submodule}.
\end{proof}

In particular, \propref{many gl to agl lemma} implies that as a $\Rep_a
U_q\agli$--comodule, the image of $(\Rep\, U_q\gli)^{\otimes n}$ under
$\iota(\bs\nu)$ is an integral form of the restricted dual to the
$U\mathfrak s\mathfrak l_\infty^-$--module $L(\bs\nu)$, which is
preserved by $U_\Z\sw_\infty^-$.

\subsection{Combinatorial identity}\label{comb section}
A dominant $N$--tuple of polynomials ${\bs Q}$ is called {\em
$\bs\nu$--admissible} if ${\bs Q}$ is the highest weight of a
representation of the form $$V_{\bs\la_{
1}}(aq^{2s_1})\otimes\dots\otimes V_{\bs\la_{k}}(aq^{2s_k}).$$

\begin{conj}\label{can conj}
The classes of irreducible representations $V({\bs Q})$ whose 
Drinfeld polynomials ${\bs Q}$ are
$\bs\nu$--admissible, form a basis in
$\Rep_a^{\bs\nu} U_q \gli$.
\end{conj}

Clearly, Conjecture \ref{can conj} holds for $k=1$. In general,
Conjecture \ref{can conj} can be derived from results in
representation theory of affine Hecke algebras (see, e.g.,
\cite{BZ,A,Gr}) via the quantum affine Schur duality functor
\cite{J,CP:schur}. It would be interesting to give a direct proof.
 
We end this section by describing a combinatorial corollary of
Proposition \ref{many gl to agl lemma} and Conjecture \ref{can conj}.

Fix ${\bs s}=(s_1\geq\dots\geq s_k=0$), such that $s_i\in\Z$.
By definition, an ${\bs s}$--diagram is a union of boxes
(with multiplicities) which can be represented as a union of $k$ Young
diagrams where the $i$-th Young diagram is shifted in the horizontal
direction by $s_i$.

For example, let $k=2$, $s_1=1$, and $s_2=0$. Then the union of the
Young diagrams of shapes $(2,1)$ (shifted by $1$) and $(3,1)$ looks as
follows:

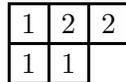
\begin{picture}(0,20)(0,20)
\setlength{\unitlength}{1mm}
\thicklines
\put(50,30){\framebox(5,5){1}}
\put(55.3,30){\framebox(5,5){2}}
\put(60.6,30){\framebox(5,5){2}}

\put(50,24.7){\framebox(5,5){1}}
\put(55.3,24.7){\framebox(5,5){1}}

\end{picture}

Here the number inside each box shows its multiplicity.

Note that the same diagram can be obtained as a union of two other
Young diagrams corresponding to shapes $(2,0)$ (shifted by $1$) and
$(3,2)$.

In the case $k=1$, we obtain in this way the usual Young diagrams. If
all $s_i=0$, then we obtain the so-called plain partitions of
magnitude $k$.

Let $N_n({\bs s})$ be the number of distinct ${\bs s}$--diagrams with
$n$ boxes and let the series 
$\chi_{\bs s}(\xi)=\sum_{i=0}^\infty N_i({\bs s}) \xi^i$
be the corresponding generating function.

Conjecture \ref{can conj} and Proposition \ref{many gl
to agl lemma} imply that the function $\chi_{\bs s}(\xi)$ coincides with the
character of the integrable $\sw_\infty$--module of weight
$\sum_{i=1}^k\varpi_{s_i}$ in the principal gradation. This character
may be found, e.g., in \cite{FKRW}. Therefore we obtain

\begin{lem} Conjecture \ref{can conj} implies that
\be \chi_{\bs
s}(\xi)=\prod_{i=1}^\infty(1-\xi^i)^{-k}{\prod_{1\leq
i<j\le k}(1-\xi^{s_i-s_j+j-i})}.  \ee
\end{lem}

It would be interesting to find a purely combinatorial proof of this
identity.  In the case of plane partitions (i.e., when all $s_i=0$),
such a proof is given in \cite{M} \footnote{We thank I. Gessel for
pointing out this reference to us.} (our formula is obtained from
formula in Section 429 of \cite{M} by taking the limit $l,m\to
\infty$).

\section{Quantum algebras at roots of unity}\label{algebras at roots}

In this section we collect information about the restricted
specializations of the algebras $U_q\sln$, $U_q\gln$, $U_q\asln$,
$U_q\agln$ to roots of unity. In the case of $U_q\asln$ these results
were obtained in \cite{CP:root,FM2}.

From now on we consider the algebras $U_q\sln$, $U_q\gln$,
$U_q\asln$, $U_q\agln$ as defined over the field $\C(q)$ 
of rational functions in $q$.

\subsection{Restricted integral form of $U_q \asln$} 

If $A$ is an invertible element of a $\C(q)$--algebra $\mc A$, then we
will denote by $A^{(r)}$, $\bin[A;r]$ ($r\in\Z_{\geq 0}$) the
following elements of $\mc A$: \be A^{(r)}=\frac{A^r}{[r]_q!},\qquad
\bin[A;r]=\prod_{i=1}^r\frac{Aq^{1-i}-A^{-1}q^{i-1}} {q^i-q^{-i}}.
\ee

Let $\uqr\asln$ be the $\C[q,q^{-1}]$--subalgebra of $U_q^{\rm
DJ}\asln$ generated by $K_i^{\pm1}$ and $X^{\pm (r)}_i$ ($i=0,1\dots,N-1$,
$r\in\Z_{>0}$). We call $\uqr\asln$ the {\em restricted integral form
of $U_q\asln$}.

We obtain from formula \Ref{sln comult} the following formula for the
comultiplication of divided powers in $U_q^{\rm DJ}\asln$:
\begin{align*}
\Delta(X_i^{+(r)})&=\sum_{j=0}^r q^{j(r-j)}\;X_i^{+(r-j)}K_i^j\otimes
X_i^{+(j)}, \\
\Delta(X_i^{-(r)})&=\sum_{j=0}^r q^{j(r-j)}\;X_i^{-(j)}\otimes
X_i^{-(r-j)}K^{-j}, 
\end{align*}
($i=0,1,\dots,N-1$, $r\in\Z_{>0}$). Therefore we have:
\begin{lem}
The algebra $\uqr\asln$ is a Hopf subalgebra of $U_q^{\rm DJ}\asln$.
\end{lem}

We have the following description of the algebra $\uqr\asln$ in the
``new'' realization.

Denote by $U_q^{{\rm res},\pm}\asln$ the $\C[q,q^{-1}]$--subalgebra of
$U_q^{\rm 
new}\asln$ generated by $x_{i,n}^{\pm(r)}$ ($n\in\Z$, $r\in\Z_{>0}$,
$i=1,\dots,N-1$). Denote $U_q^{{\rm res},0}\asln$ the $\C[q,q^{-1}]$ subalgebra
of $U_q^{\rm new}\asln$ generated by $k_i^{\pm 1}$, $\bin[k_i;r]$,
${\mc P}_{i,n}$ ($i=1,\dots,N-1$, $n\in\Z$, $r\in\Z_{>0}$).
 
\begin{prop}[\cite{CP:root}, Proposition 6.1] \label{sln tri} We have
the following embeddings:
\be
U_q^{{\rm res},\pm}\asln\subset\uqr\asln,\qquad U_q^{{\rm res},0}\asln
\subset\uqr\asln.
\ee
Moreover, the multiplication gives rise to
the following isomorphism of vector spaces: 
\bean\label{triangular}
\uqr\asln=U_q^{{\rm res},-}\asln\otimes U_q^{{\rm res},0}\asln\otimes U_q^{{\rm res},+}\asln.
\eean
\end{prop}

Note that $\bin[k_i^{-1};r]\in\uqr\asln$ and
moreover, $\bin[k_i^{-1};r]$ belongs to the algebra generated by
$k_i^{\pm1}$ and $\bin[k_i;s]$
($i=1,\dots,N-1$, $s,r\in\Z_{>0}$, $s<r$).

By \cite{L}, \S 37.1.3, the braid group automorphisms
${\mc B}_i^{(\pm, j)}$ ($i=0,1,\dots,N-1$, $j=1,2$) 
given by formulas \Ref{braid1}, \Ref{braid2}
preserve the integral form $\uqr\asln\subset U_q\asln$.

The shift automorphisms $\tau$ (see formula \Ref{taua}) also preserve
the integral form $\uqr\asln$.

Let $\uqr\sln$ be the $\C[q,q^{-1}]$--subalgebra of $\uqr\asln$
generated by $K_i^{\pm1}$, $X^{\pm (r)}_i$ ($i=1,\dots,N-1,$ $r\in\Z_{>0}$).
It is also a Hopf subalgebra preserved by ${\mc B}_i^{(\pm, j)}$
($i=1,\dots,N-1$, $j=1,2$).

\subsection{Restricted integral forms of $U_q\agln$ and $U_q\gln$}

We define the {\em restricted integral form $\uqr\agln$} as the
$\C[q,q^{-1}]$--subalgebra of $U_q^{\rm new}\agln$ generated by
$t_i^{\pm1}$, $\bin[t_i;r]$, ${\mc Q}_{i,n}$,
$x_{i,n}^{\pm(r)}$ ($i=1,\dots,N$, $n\in\Z$, $r\in\Z_{>0}$). Note
that this set is not minimal, for example $\mc Q_{i,n}$ ($i=2,\dots,N$)
can be expressed via other generators. 

Introduce the notation
\be
g_{i,(n)}=g_{i,n}/[n]_q.
\ee
Then the generators ${\mc Q}_{i,n}$ may be replaced by $g_{i,(n)}$
($i=1,\dots,N$, $n\in\Z\setminus 0$).

We expect that the algebra $\uqr\agln$ is a Hopf subalgebra of
$U_q\agln$, but we do not have a proof of this statement. Note,
however, that $\uqr\agln$ contains a ``large'' subalgebra isomorphic to
$\uqr\asln\otimes {\mc C}^{\rm res}$ which is a Hopf algebra, see
formula \Ref{det decom root}. 

Next, we define $\uqr\gln$ as the $\C[q,q^{-1}]$--subalgebra of
$U_q^{\rm DJ}\gln$ 
generated by $T_i^{\pm1}$, $\bin[T_i;r]$, $X_j^{\pm(r)}$
($i=1,\dots,N$, $j=1,\dots,N-1$, $r\in\Z_{>0}$).

\begin{lem}
The algebra  $\uqr\gln$ is a Hopf subalgebra of $U_q\gln$.
\end{lem}

\begin{proof}
Since we already know that $U_q^{\on{res}} \sw_N$ is a Hopf algebra,
it remains to show that the $\C[q,q^{-1}]$--subalgebra
of $U_q\gln$ generated by $T_i^{\pm
1},\bin[T_i;r]$ ($i=1,\ldots,N$, $r \in \Z_{>0}$), is a Hopf
subalgebra.

Let $A$ be an invertible element of some
$\C(q)$--Hopf algebra $\mc A$
such that $\Delta(A)=A\otimes A$. Then we claim that the
$\C[q,q^{-1}]$--subalgebra of ${\mc A}$ generated by
$\bin[A;r],A^{\pm1}$ 
($r\in\Z_{>0}$) is a Hopf subalgebra. It is sufficient
to check this for some algebra $\mc A$ with an element $A$, such that
$A^{\pm 1}$ generate $\C(q)[A,A^{-1}]\in{\mc A}$. An example of such algebra is
${\mc A}=U_q{\mathfrak {sl}}_2$ with $A=K$. Then $\Delta(\bin[K;r])$ 
is in $(\uqr{\mathfrak {sl}}_2)^{\otimes 2}$ (because
$\uqr{\mathfrak {sl}}_2$ is a Hopf subalgebra) and in $(\C(q)[K^{\pm
1}])^{\otimes 2}$. Our claim follows then because the intersection of
$\uqr{\mathfrak {sl}}_2$ and $\C(q)[K^{\pm 1}]$ in $U_q{\mathfrak
{sl}}_2$ is precisely the $\C[q,q^{-1}]$--subalgebra generated by
$K^{\pm 1}$ and $\bin[K;r]$ ($r\in\Z_{>0})$.
\end{proof}

Denote by $U_q^{{\rm res},\pm}\agln$ the $\C[q,q^{-1}]$--subalgebra of $U_q^{\rm
new}\agln$ generated by $x_{i,n}^{\pm(r)}$ ($n\in\Z$, $r\in\Z_{>0}$,
$i=1,\dots,N-1$). Denote by $U_q^{{\rm res},0}\agln$ the
$\C[q,q^{-1}]$--subalgebra of $U_q^{\rm new}\agln$ generated by
$t_i^{\pm 1}$, $\bin[t_i;r]$, $\bin[t_jt_{j+1}^{-1};r]$, ${\mc Q}_{i,n}$
($i=1,\dots,N$, $j=1,\dots,N-1$, $n\in\Z$, $r\in\Z_{>0}$).

Using the defining relations of \secref{agln}, we obtain \bean
\label{x and g} \left[{g_{i,(n)}},x^{\pm(r)}_{j,m}\right] =
(\pm\delta_{i,j}\mp \delta_{i,j+1})\frac{ q^{nj}}{n[r]_q!}
\sum_{s=0}^r(x^\pm_{j,m})^sx_{j,n+m}^\pm(x^\pm_{1,m})^{r-1-s}, \eean where
the right hand side belongs to $\uqr\asln$ by Proposition 4.7 of
\cite{CP:root}. These relations, together with \propref{sln tri},
imply the following result.

\begin{lem}    \label{tr dec gln}
The multiplication map gives rise to an isomorphism of vector spaces:
\bea
\uqr\agln=U_q^{{\rm res},-}\agln\otimes U_q^{{\rm res},0}\agln\otimes
U_q^{{\rm res},+}\agln.
\eea
\end{lem}

Now we can specialize $q$ to any non-zero complex number.  Fix
$\ep\in\C^\times$. Define the structure of a module over the ring
$\C[q,q^{-1}]$ on $\C$
by the formula
\be p(q,q^{-1})\cdot
c=p(\epsilon,\ep^{-1})c,\qquad p(q,q^{-1})\in \C[q,q^{-1}],\;\;c\in\C.
\ee 
Denote this module by $\C_\ep$. 
Set 
\bea \ur\sln=\uqr\sln\otimes_{\C[q,q^{-1}]}\C_\ep,\qquad
\ur\asln=\uqr\asln\otimes_{\C[q,q^{-1}]}\C_\ep,\\
\ur\gln=\uqr\gln\otimes_{\C[q,q^{-1}]}\C_\ep,\qquad 
\ur\agln=\uqr\agln\otimes_{\C[q,q^{-1}]}\C_\ep.  \eea 
We will use the
same notation for the elements of the specialized algebras, e.g., we
will write simply $x_i^{\pm(r)}$ for $x_i^{\pm(r)}\otimes 1\in\ur\agln$.

\subsection{Small quantum affine algebras}

Denote by $U^{\on{fin}}_\ep \asln$ the subalgebra of $\ur\asln$
generated by $k_i^{\pm1}$, $h_{i,n}$ and $x_{i,r}^\pm$ ($i\in
1,\dots,N-1$ $r\in \Z$, $n\in\Z$).

From the formulas for the isomorphism $\on B$ given in 
Proposition \ref{B lemma} we obtain the following

\begin{lem}
In the Drinfeld-Jimbo realization, the algebra $U^{\on{fin}}_\ep
\asln$ is isomorphic to the subalgebra of $\ur\asln$ generated by
$K_i^{\pm 1}$ and $X_i^\pm$ ($i=0,\dots,N-1$).
\end{lem}

Denote by $U^{\on{fin}}_\ep \agln$ the subalgebra of $\ur\agln$
generated by $t_i^{\pm1}$, $g_{i,n}$ and $x_{j,r}^\pm$ ($i\in
1,\dots,N$, $j=1,\dots,N-1$, $n\in\Z\setminus 0$). We call the algebras
$U^{\on{fin}}_\ep \asln$ and $U^{\on{fin}}_\ep \agln$, {\em small}
quantum affine algebras.

If $\ep$ is not a root of unity, then the algebras $\fu\asln$ and
$\fu\agln$ coincide with $\ur\asln$ and $\ur\agln$. 

Let $\ep$ be a primitive root of unity (of order greater than 2) 
and let $l$ be the smallest natural number with
the property $\ep^{2l}=0$ or, equivalently, $[l]_\ep=0$.
Then $h_{i,jl}=0$ (resp. $g_{i,jl}=0$) and $(x^\pm_{i,k})^l=0$
($k\in\Z$, $j\in\Z\setminus 0$) in $\fu\asln$ (resp. $\fu\agln$). Moreover,
$\fu\asln$ (resp. $\fu\agln$) is a specialization of the subalgebra of
$\uqr\asln$ (resp. $\uqr\agln$) generated by $k_i^{\pm1}$, $h_{i,n}$ and
$x_{i,r}^\pm$ (resp. $t_i^{\pm1}$, $g_{i,n}$ and $x_{i,r}^\pm$),
modulo these relations. In particular, $\fu\asln$ and
$\fu\agln$ inherit a Hopf algebra structure from $\uqr\asln$
and $\uqr\agln$. However, we will not consider this
Hopf algebra structure in our paper.

Let $\fu\gln$ be the subalgebra of $\ur\gln$ generated by
$T_i^{\pm1}$, $X_j^\pm$ ($i=1,\dots,N-1$, $j=1,\dots,N$).

\subsection{Some maps}\label{inclusions section root}
The maps of quantum algebras from Section \ref{inclusions section}
induce maps at the level of integral forms and their specializations. 

We have the following obvious embeddings of algebras:
\bea
\ur\widehat{\sw}_{N-1}&\to &\ur\asln,\qquad
\fu\widehat{\sw}_{N-1}\to \fu\asln,\\
\ur\widehat{\gw}_{N-1}&\to &\ur\agln,\qquad
\fu\widehat{\gw}_{N-1}\to \fu\agln,\\
\ur{\gw}_{N-1}&\to& \ur\gln, \qquad \fu{\gw}_{N-1}\to \fu\gln,
\eea
given by the ``upper left corner'' rule (the first map is defined in
the new realization of $\ur\asln$).

We also have an embedding of algebras
\bea
\ur\gln &\to & \ur\agln,\qquad \fu\gln\to\fu\agln
\eea
mapping $X_i^{\pm(r)}$ to $x_{i,0}^{\pm (r)}$ and
$T_i^{\pm1}$ to $t_i^{\pm1}$ ($i=1,\dots, N$; $r\in\Z_{>0}$ for the
first map and $r=0$ for the second one).

{}From formula \Ref{det decom} we obtain an embedding of algebras
\bean\label{det decom root} \uqr\asln\otimes {\mathcal C}^{\rm
res}\to\uqr\agln, \eean 
where ${\mc C}^{\on{res}}$ denotes the
(central) subalgebra of $\uqr\agln$ generated by the elements $C_n, n
\in \Z \backslash 0$ defined by formula \eqref{Cn}, and by
$\prod_{i=1}^Nt_i^{\pm1}$.
However, if $\ep$ is a root of unity then 
there is no analogue \Ref{det decom2}
at the level of restricted integral forms (due to the different
choice of generators in the Cartan parts of the algebras).

In addition, we obtain from \Ref{det decom} and \Ref{det decom2} the
embeddings \bean\label{det decom root fin} \fu\asln\otimes {\mathcal
C}^{\on{fin}}\to\fu\agln,\qquad \fu\agln\to\fu\asln\otimes {\mathcal
C}^{\on{fin}},\eean where ${\mc C}^{\on{fin}}$ is the
commutative algebra generated by $\sum_{i=1}^Ng_{i,n}$ ($n\in\Z$,
$\ep^{2n}\neq 1$) and $\prod_{i=1}^N t_i$ (i.e., ${\mc C}^{\on{fin}}$
is the image of ${\mc C}$ in $\ur\agln$). In the second embedding we
extend $\fu\asln$ by $k_i^{1/N}$ and ${\mc C}^{\on{fin}}$ by
$\prod_{i=1}^N t_i^{1/N}$.

\subsection{Evaluation homomorphism}\label{evaluation section root}
Recall there is the 
evaluation algebra homomorphism $ev_a:\;U_q\agln \to
U_q\gln$ defined by the formula \eqref{gl evaluation}.  In this
section we show that $ev_a$ descends to a homomorphism of restricted
specializations and of the small quantum algebras.

\begin{lem} 
The image of $\uqr\agln \subset U_q\agln$ under $ev_a$ is equal to
$\uqr\gln$, for any $a\in\C^\times$.
\end{lem}

\begin{proof}
The algebra $\uqr\agln$ is generated by $\uqr\asln$, ${\mc
Q}^\pm_1(u)$ and $t_i^{\pm 1}$, $\bin[t_i;r]$ ($i=1,\dots, N$,
$r\in\Z_{>0}$).  Moreover, by Theorem \ref{cube theorem}, the
evaluation map commutes with the embedding $U_q\asln\to
U_q\agln$. Therefore we can use formulas \Ref{ev sln} to compute the
evaluation map on $\uqr\asln$.

For $i=1,\dots,N$, $j=1,\dots,N-1$, $r\in\Z_{>0}$,
we have \be
ev_a(x_{j,0}^{\pm(r)})=ev_a(X_j^{\pm(r)})=X_j^{\pm(r)},\qquad
ev_a(t_i^{\pm1})=T_i^{\pm1}, \qquad
ev_a(\bin[t_i;r])=\bin[T_i;r].  \ee In particular, we
obtain that the image of $\uqr\gln$ under $ev_a$ contains $\uqr\gln$.

In order to find $ev_a(X_{0}^\pm)$, we use the presentation \Ref{ev
braid}. We then obtain \bean\label{ev via braid}
ev_a(X_0^{+(n)})&=&((-q)^Naq^{-1}T_1T_N)^n {\mc B}_1^{(+,2)} {\mc
B}_2^{(+,2)}\dots {\mc B}_{N-2}^{(+,2)}(X_{N-1}^{-(n)}),\notag\\
ev_a(X_0^{-(n)})&=&(a^{-1}q^{-1}T_1^{-1}T_N^{-1})^n{\mc
B}_1^{(-,1)} {\mc B}_2^{(-,1)}\dots {\mc
B}_{N-2}^{(-,1)}(X_{N-1}^{+(n)}).  \eean Therefore, the result is in
$\uqr\gln$.

Finally we prove that $ev_a({\mc Q}^\pm_1(u))\in\uqr\gln$. Note that
the evaluation map commutes with the left upper corner restrictions to
$U_q\widehat{\mathfrak gl}_M$ ($M<N$). Therefore, we are reduced to the
case of $U_q\widehat{\mathfrak g}{\mathfrak l}_1$.  By formula \Ref{gl
evaluation}, we obtain \be ev_a(\Theta_1^\pm(u)) = ev_a(L^\pm(u)) =
\frac{L^\pm-(u/a)^{\pm1}L^\mp}{1-(u/a)^{\pm1}} =
\frac{T_1^{\pm1}-(u/a)^{\pm1}T_1^{\mp1}}{1-(u/a)^{\pm1}}.  \ee
Therefore using formula \Ref{theta q} we find \bean\label{ev Q}
ev_a({\mc Q}_1^\pm(u)) = \prod_{i=1}^\infty
\frac{1-(u/a)^{\pm1}q^{-2i-1}}{1-(u/a)^{\pm1}T_1^{\mp2}q^{-2i-1}}=
\sum_{i=0}^\infty \left(\prod_{j=0}^{i-1}\frac{T_1^{\mp
(2j+2)}-q^{-2j}}{1-q^{-(2j+2)}} \right)(u/aq)^{\pm i}\notag\\
=\sum_{i=0}^\infty \bin[T_1^{\pm1};i] T_1^{\mp i}q^{i}(u/a)^{\pm
i}\in\uqr\gln.  \eean This completes the proof.
\end{proof}
Therefore we have a surjective homomorphism of algebras which we also
denote by $ev_a$:
\be
ev_a:\;\ur\agln\to \ur\gln.
\ee

Similarly, we obtain a surjective homomorphism of algebras
denoted in the same way:
\be
ev_a:\;\fu\agln\to \fu\gln.
\ee

For $a\in\C^\times$, the 
shift automorphism $\tau_a$ defined by \Ref{twist} also clearly
induces an automorphism of $\ur\agln$ which we also denote by
$\tau_a$. We have
\bean\label{twist root}
\tau_a(g_{i,n})=a^ng_{i,n},\qquad \tau_a(t_i^{\pm1})=t_i^{\pm1},\qquad
\tau_a(x_{i,n}^{\pm(r)})=a^{nr}x_{i,n}^{\pm(r)}.
\eean
The same formulas with $r=0$
define the shift automorphism $\tau_a$ of $\fu\agln$.

In both cases, we have $ev_a\circ\tau_b=ev_{ab}$.

\section{Finite-dimensional representations and $q$-characters at
roots of unity}\label{rep at roots}

From now on, we fix $\ep$ to be a primitive root of unity of order
$s$. We set $l$ to be the order of $\ep^2$ (equivalently, the order of
$2$ in $\Z/s\Z$). Thus, \bea l=\left\{
\begin{matrix} 
s & s \;\;  {\rm is \;\; odd},\\
s/2 & s \;\; {\rm is \;\; even}.
\end{matrix}
\;\right.
\eea
Note the identities
\be
q^s=1,\qquad [l]_\ep=0,\qquad (-\ep)^l=-\ep^{l^2}.
\ee

\subsection{Finite-dimensional representations of $\ur\asln$}
According to the results of \cite{CP:root}, the
classification of finite-dimensional representations of $\ur\asln$ is
similar to that of $U_q\asln$ (presented in Section \ref{fd}). We
briefly review these results.

A vector $w$ in a $\ur\sln$--module $W$ is called a vector of weight
$\bs \la$, if \be K_i \cdot w = \ep^{(\bs \la,\al_i)} w, \qquad
\bin[K_i;l]\cdot w=\left[\begin{matrix}(\bs
\la,\al_i)\\l\end{matrix}\right]_\ep w, \qquad (i=1,\dots,N-1).\ee  A
representation $W$ of $U_q\sln$ is said to be of type 1 if $W$ is the
direct sum of its weight spaces $W = \oplus_{\bs \la} W_{\bs\la}$,
where $W_{\bs \la}$ is the space of all vectors in $W$ of weight
$\bs\la$.

A representation $V$ of $\ur\asln$ is said to be of type 1 if $V$ is
of type 1 as a module over $\ur\sln\subset\ur\asln$. In what follows,
we consider only finite-dimensional type 1 representations.

Let $V$ be a representation of $\ur\asln$. A vector $v\in V$ is
called a highest weight vector if \be x_{i,n}^{+(r)} \cdot v=0,\qquad 
{\mc P}_{i,n}\cdot v = P_{i,n}v, \qquad k_i \cdot v =
\ep^{\mu_i} v, \qquad \bin[k_i;r]v=b_{i,r}v \ee ($i=1,\dots,N-1$, $n
\in \Z$, $r\in\Z_{>0}$), where $P_{i,n},b_{i,r} \in \C$ and $\mu_i\in
\Z_{\geq 0}$. A representation $V$ is called a highest weight
representation if $V=\ur\asln\cdot v$, for some highest weight vector
$v$.

Every finite-dimensional irreducible representation $V$ of $\ur\asln$
of type 1 is a highest weight representation and \be P_i^\pm(u):=
\sum_{n\geq 0} P_{i,\pm n} u^{\pm n}=\prod_{j=1}^{\mu_i}(1-a_j^{\pm
1}u^{\pm 1}) \ee is a polynomial. In particular, the polynomials
$P_i^+(u)$ determine the polynomials $P_i^-(u)$ the numbers
$\mu_i=\deg P_i^+(u)$ and $b_{i,r} ={\left[\begin{matrix}{\deg
P_i^+(u)}\\{r}\end{matrix}\right]_\ep}$. For any $(N-1)$--tuple ${\bs
P} = (P_1(u),\ldots,P_{N-1}(u))$ of polynomials with constant term 1,
there exists a unique representation of $\ur \asln$ with such Drinfeld
polynomials. It is denoted by $\wt{V}_\ep({\bs P})$.

For $a\in\C^\times$, the Drinfeld polynomials of the pullback
$\tau_a^*(\wt{V}_\ep({\bs P}))$ of $\wt{V}_\ep({\bs P})$ by the
shift automorphism $\tau_a$ are given by $P_i(au)$ ($i=1,\dots,N-1$).

We denote by $\Rep\,\ur\asln$ the Grothendieck ring of all
finite-dimensional $\ur\asln$ modules of type 1. 

\subsection{Finite-dimensional representations of
$\ur\agln$}\label{fin dim root}
We introduce the notion of polynomial finite-dimensional
representations of $\ur\agln$ analogous to the case of generic $q$
presented in Section \ref{agln rep}.

A vector $w$ in a $\ur\gln$--module $W$ is called a vector of weight
$\bs\la=(\la_1,\dots,\la_N) \in \Z^n$, if 
\be
T_i \cdot w = \ep^{\la_i} w, \qquad\bin[T_i;r]\cdot
w=\bin[\la_i;r]_\ep w,\qquad (i=1,\dots,N,\;\; r\in\Z_{>0}).
\ee
A representation $W$ of $\ur\gln$ is said to be
of type 1 if it is the direct sum of its weight spaces $W =
\oplus_{\bs\la} W_{\bs\la}$, where $W_{\bs\la}$ is the space of all
vectors in $W$ of weight $\bs \la$.

Let ${\mc F}_N^{\ep}$ 
be the category of finite-dimensional representations
$W$ of $\ur\gln$ of type 1, such that all weights ${\bs \la} =
(\la_1,\dots,\la_N)$ occurring in $W$ have non-negative integral
components, i.e., $\la_i\in\Z_{\geq 0}$ for all $i=1,\dots,N$.
The simple objects of category ${\mc F}_N^{\ep}$ 
are the irreducible representations with highest
weights $\bs\la=(\la_1,\dots,\la_N)$, where
$\la_1\geq\la_2\geq\dots\geq \la_N\geq 0$. We denote the corresponding
Grothendieck ring by $\Rep\,\ur\gln$.

Let $\wh{\mathcal F}_N^\ep$ be the category of all finite-dimensional
representations $V$ of $\ur\agln$, such that

\begin{itemize}
\item the restriction of $V$ to $\ur\gln$ is in ${\mathcal F}_N^\ep$;

\item each generalized common eigenvalue of ${\mathcal Q}^\pm_i(u)$ 
in $V_{\bs\la}$, where ${\bs \la} = (\la_1,\ldots,\la_N)$,
is a polynomial $\Gamma^\pm_i(u)$ in $u^{\pm 1}$ of degree
$\la_i$, and the zeroes of the functions $\Gamma^+_i(u)$ and
$\Gamma^-_i(u)$ coincide.
\end{itemize}

We will call such representations {\em polynomial}.

Let $V$ be a representation of $\ur\agln$.  A vector $v\in V$ is
called a highest weight vector if \be x_{j,n}^{+(r)} \cdot v=0,\quad
\quad {\mc Q}_{i,n}\cdot v = Q_{i,n}v, \quad \quad t_i \cdot v =
\ep^{\mu_i} v,\quad \bin[t_i;r]v=b_{i,r}v \ee ($i=1,\dots,N$,
$j=1,\dots,N-1$, $n \in \Z$, $r\in\Z_{>0}$), where $Q_{i,n},b_{i,r}
\in \C$ and $\mu_i\in \Z_{\geq 0}$. A representation $V$ is called a
highest weight representation if $V=\ur\agln\cdot v$, for some highest
weight vector $v$. If $V$ is irreducible, then the polynomials
$Q_i(u): = \sum_{n\geq 0} Q_{i,n} u^{n}$ are called the Drinfeld
polynomials of $V$, and the $N$--tuple ${\bs Q} =
(Q_1(u),\ldots,Q_N(u))$ is called the highest weight of $V$.

An $N$--tuple ${\bs Q}$ of polynomials with constant terms $1$, ${\bs
Q}=(Q_1(u),\dots,Q_N(u))$, is called {\em dominant}
if the ratio $Q_i(u)/Q_{i+1}(u\ep^2)$ ($i=1,\dots, N-1$)
is also a polynomial.

\begin{lem}\label{class gln root}
Every irreducible polynomial representation $V$ of $\ur\agln$ is a
highest weight representation and \be Q_i(u)= \sum_{n\geq 0} Q_{i,n}
u^{n}=\prod_{j=1}^{\mu_i}(1-a_ju)\qquad (i=1,\dots,N) \ee is a polynomial
such that $(Q_1(u),\dots,Q_N(u))$ is a dominant $N$--tuple. Then
\be
\sum_{n\leq 0} Q_{i,n}u^{n}=\prod_{j=1}^{\mu_i}(1-a_j^{-1}u^{-1}),\qquad
\mu_i=\deg Q_i(u),\qquad b_{i,r}
={\left[\begin{matrix}{\deg Q_i(u)}\\{r}\end{matrix}\right]_\ep}.
\ee
Moreover, for any dominant $N$--tuple ${\bs Q} =
(Q_1(u),\ldots,Q_{N}(u))$ of polynomials with constant term 1, there
exists a unique irreducible representation of $\ur \agln$ with highest
weight ${\bs Q}$.
\end{lem}

\begin{proof}
This lemma is proved in exactly the same way as the corresponding
statement in the $\ur\asln$ case using the triangular decomposition of
\lemref{tr dec gln}.
\end{proof}

Let $\Rep\,\ur\agln$ be the Grothendieck group of the 
category $\wh{\mc F}^\ep_N$. Since we have not shown that the
comultiplication on $U_q \agln$ maps $\uqr\agln$ to
$(\uqr\agln)^{\otimes 2}$, we do not immediately have a ring structure
on $\Rep\,\ur\agln$. But we will define it by other means in Section
\ref{spec section}.

One checks readily that the restriction of the irreducible polynomial
representation of $\ur\agln$ corresponding to a dominant $N$--tuple
${\bs Q}$ to $\ur\asln$ is the irreducible representation of
$\ur\asln$ with $P_i(u)=Q_i(uq^{i-1})/Q_{i+1}(uq^{i+1})$. Therefore we
have an analogue of Lemma \ref{gl to sl functor}:

\begin{lem}
The restriction map $\Rep\,\ur\agln\to\Rep\,\ur\asln$ is
surjective.
\end{lem}

\subsection{The $\ep$--characters}

The theory of $\ep$--characters for $\ur \widehat{\g}$, where
$\widehat\g$ is an untwisted affine Lie algebra (in particular,
$\wh{\sw}_N$), has been developed in \cite{FM2}. Here we review this
theory and extend it to the case of $\ur\agln$.
 
The definition of the $\ep$--characters is parallel to that in the
case of generic $q$. Let $V$ be a finite-dimensional
$\ur\asln$--module and let ${v_1,\dots,v_m}$ be a basis of common
generalized eigenvectors of $k_i,\bin[k_i;r], {\mathcal
P}_i^\pm(u)$. It follows from Lemma 3.1 of \cite{FM2} that the common
eigenvalues of operators ${\mathcal P}_i^+(u)$ on each vector $v_s$
are rational functions of the form \be
\Gamma_{i,s}^\pm(u)=\frac{\prod_{k=1}^{k_{is}} (1-a_{isk}^{\pm
1}u^{\pm 1})}{\prod_{r=1}^{r_{is}} (1-b_{isr}^{\pm 1}u^{\pm 1})}, \ee
where $a_{isk},b_{isr}$ are some nonzero complex numbers. Moreover,
$k_iv_s=\ep^{k_{is}-r_{is}}v_s$ and
$\bin[k_i;t]v_s=\bin[k_{is}-r_{is};t]_\ep v_s$.  By definition, the
$\ep$--character of $V$, denoted by $\chi_\ep(V)$, is an element of
the commutative ring of polynomials $\Z[Y^{\pm
1}_{i,a}]_{i=1,\dots,N-1}^{a\in\C^\times}$, which is equal to \be
\chi_q(V)=\sum_{s=1}^m \prod_{i=1}^{N-1}\left(\prod_{k=1}^{k_{is}}
Y_{i,a_{isk}}\prod_{r=1}^{r_{is}}Y_{i,b_{isr}}^{-1}\right).  \ee

The map \be \chi_\ep:\;\Rep\,U_q\asln\to \Z[Y^{\pm
1}_{i,a}]_{i=1,\dots,N-1}^{a\in\C^\times}, \qquad [V]\mapsto\chi_q(V),
\ee is an injective homomorphism of rings. Moreover, the image of
$\chi_\ep$ is given by \be {\mc K}=\bigcap_{i=1,\dots,N-1}
\left(\Z[Y_{j,a}^{\pm 1}]_{j\neq i}^{ a \in \C^\times} \otimes
\Z[Y_{i,b} + Y^{-1}_{i,b\ep^2} Y_{i-1,b\ep}Y_{i+1,b\ep}]_{b \in
\C^\times}\right), \ee where we set $Y_{N,a}=Y_{0,a}=1$ see
\cite{FM2}, Proposition 3.6.

Fix $a\in\C^\times$. Denote by $\Rep_a\ur\asln$ the subring
of $\Rep\,\ur\asln$ spanned by those representations whose
$q$--characters are polynomials in $Y^{\pm 1}_{i,a\ep^{i+2k-1}}$
($i=1,\dots,N-1$, $k=1\dots,l$).

The class of an irreducible module $V$ is in $\Rep_aU_q\asln$
if and only if the roots of the Drinfeld polynomial $Q_i(u)$ are of
the form $\ep^{-i-2k+1}$ ($k\in\{1,\dots,l\}$) or, equivalently, the
highest weight monomial in $\chi_\ep(V)$ is a monomial in
$Y_{i,a\ep^{i+2k-1}}$ ($k \in \{1,\dots,l\}$). Moreover, we have the
decomposition of rings \be \Rep\,\ur\asln \simeq \bigotimes_{a\in\C^\times/\ep^{2\Z}}\Rep_aU_q\asln, \ee where the tensor product is over the angle of
bases of all lattices of the form $\{a\ep^{2\Z}\}$ in $\C^\times$.

Similarly, let $V\in\wh{\mc F}_N$ be a polynomial
$\ur\agln$--module. Let ${v_1,\dots,v_m}$ be a basis of common
generalized eigenvectors of $t_i$, $\bin[t_i;r]$, ${\mathcal Q}_i^\pm$
in $V$. Then 
the common generalized eigenvalues
of ${\mathcal Q}_i^\pm(u)$ on each vector $v_s$ are polynomials
$$
\Gamma_{i,s}^{\pm}(u)=\prod_{k=1}^{\la_{is}}(1-a_{isk}^{\pm 1} u^{\pm
1}).
$$
Furthermore, we have: $t_i v_s = \ep^{\la_{i,s}} v_s$,
$\bin[t_i;r]v_s=\bin[\la_i;r]_\ep v_s$.

By definition, the $\ep$--character of $V$, denoted again by
$\chi_\ep(V)$, is an element of the ring of polynomials $\Z[\Lambda
_{i,a}]_{i=1,\dots,N}^{a\in\C^\times}$, given by the formula \be
\chi_\ep(V)=\sum_{s=1}^m\prod_{i=1}^N \prod_{k=1}^{\la_{is}}
\La_{i,a_{isk}}.  \ee

Using Lemma \ref{class gln root} 
and  the one-to-one correspondence between irreducible
representations and highest weight monomials, we obtain the following

\begin{lem}    \label{ring str root}
The map \be \chi_\ep: \;\Rep\,U_\ep\agln\to
\Z[\La_{i,a}]_{i=1,\dots,N}^{a\in\C^\times}, \qquad [V]\mapsto\chi_q(V),
\ee 
is injective.
\end{lem}
Note that since we do not have a comultiplication structure in
$\ur\agln$, we cannot say that the map $\chi_\ep$ is a ring
homomorphism. Instead, we will prove that the image of the map
$\chi_\ep$ is closed under multiplication, and define the product
structure on $\Rep\,U_\ep\agln$ by taking the pullback of the product in
$\Z[\La_{i,a}]_{i=1,\dots,N}^{a\in\C^\times}.$

Recall the homomorphism $\kappa_\ep:\Z[\La_{i,a}]_{i=1,\dots,N}^{a\in\C^\times}\to \Z[Y^{\pm
1}_{i,a}]_{i=1,\dots,N-1}^{a\in\C^\times}$ defined by \Ref{kappa}.
From formula \eqref{PandQ} we obtain the following
\begin{lem}
Let $V \in \wh{\mc F}_N^\ep$ be a $\ur\agln$--module with $\ep$--character
$\chi_\ep(V)$. Then the $\ep$--character of the restriction of $V$ to
$\ur\asln$ is equal to $\kappa_\ep(\chi_\ep(V))$.
\end{lem}

A monomial $m$ in $\Z[\La_{i,a}]_{i=1,\dots,N}^{a\in\C^\times}$ is
called {\em dominant} if $\kappa_\ep(m)$ does not contain negative powers
of $Y_{i,a}$.

Fix $a\in\C^\times$. Denote by $\Rep_a\ur\agln$ the subspace of
$\Rep\,\ur\agln$ spanned by all representations whose
$\ep$--characters belong to
$\Z[\La_{i,a\ep^{2k}}]_{i=1,\dots,N}^{k=1,\dots,l}$.  Thus, we have an
injective map \be \Rep_a\ur\agln \to \Z[\La_{i,a\ep^{2k}}]_{i=1,\dots,N}^{k=1,\dots,l}.  \ee
If $V$ is an irreducible representation in the category $\wh{\mc
F}_N^\ep$, then $[V]$ belongs to the subring $\Rep_a\ur\agln$
if and only if the highest weight monomial of $V$ is a monomial in
$\La_{i,a\ep^{2k}}$ ($k\in \{1,\dots,l\}$).  We have the decomposition of
rings
\bean \label{dec res}
\Rep\,\ur\agln \simeq
\bigotimes_{a\in\C^\times/\ep^{2\Z}}\Rep_a\ur\agln. \eean

The pullback with respect to the homomorphism $\tau_{b/a}$
(see formula \Ref{twist root}) gives rise to an isomorphism
$\Rep_a\ur\agln$ with $\Rep_b\ur\agln$.
Finally, it is clear from the definition that the homomorphism
$\chi_\ep$ is
compatible with all embeddings discussed in Section \ref{inclusions
section root}.

\subsection{Specialization of representations}\label{spec section}

Let $V$ be an irreducible finite-dimensional polynomial representation
of $U_q\agln$ with highest weight vector $v$. Set 
\be V^{\rm res}=\uqr\agln\cdot v\subset V.\ee

The following lemma is proved in the same way as Lemma 4.6 (i) in
\cite{CP:weyl}.

\begin{lem}
$V^{\rm res}$ is a free $\C[q,q^{-1}]$--module of rank $\dim V$.
\end{lem}

Denote by $V^\ep$ the specialization of $V$ at $q=\ep$: \be
V^\ep=V^{\rm res} \otimes_{\C[q,q^{-1}]}\C_\ep.  \ee The next lemma
follows from the definitions.

\begin{lem}    \label{ep-char of spec}
$V^{\ep}$ is a polynomial representation of $\ur\agln$. Moreover, the
$\ep$--character of $V^\ep$ is obtained from the $q$--character of $V$
via replacing $q$ by $\ep$.
\end{lem}

Set \bea \wh{\mc S}^l_N:\;\Z[\La_{i,aq^{2k}}]_{i=1,\dots,N}^{k\in\Z} \to\Z[\La_{i,a\ep^{2k}}]_{i=1,\dots,N}^{k=1,\dots,l},\qquad
\La_{i,aq^{2k}}\mapsto \La_{i,a\ep^{2k}}.  \eea Then $\wh{\mc
S}^l_N$ is a surjective ring homomorphism. We call $\wh{\mc S}^l_N$
the {\em specialization} map.

We also denote by $\wh{\mc S}^l_N$ the specialization map at the level
of the Grothendieck groups, 
\bea
\wh{\mc S}^l_N:\;\Rep\,\uqr\agln&\to& \Rep\,\ur\agln, \eea
sending $[V]$ to $[V^\ep]$ for all irreducible $V$.

The next lemma follows from Proposition 2.5 and Theorem 3.2 of
\cite{FM2}.

\begin{lem}    \label{surje}
The specialization map $\wh{\mc S}^l_N:\;\Rep\,\uqr\agln\to \Rep\,\ur\agln$ 
is surjective. Moreover,  $\wh{\mc
S}^l_N(\Rep_a\uqr\agln)=\Rep_a\ur\agln$ for any $a\in\C^\times$, 
and the following diagram
is commutative: 
\be
\begin{CD} \Rep_aU_q\agln
@>{\chi_q}>> \Z[\La_{i,aq^{2k}}]_{i=1,\dots,N}^{k\in\Z} \\
@V{\wh{\mc S}^l_N}VV @V{\wh{\mc S}^l_N}VV\\
 \Rep_a\ur\agln
@>{\chi_\ep}>>\Z[\La_{i,a\ep^{2k}}]_{i=1,\dots,N}^{k=1,\dots,l}.
\end{CD}
\ee
\end{lem}

In particular, the image of $\chi_\ep$ is closed under
multiplication. We define the ring structure on $\Rep\,\ur\agln$
by the formula
\bean \label{product ep}
[V]\cdot[W]=\chi_\ep^{-1}(\chi_\ep(V)\chi_\ep(W)).
\eean
In other words, we define the product structure on $\Rep\, \ur\agln$
in such a way that
$\chi_\ep$ becomes a ring homomorphism. Then for any $a\in\C^\times$,
the specialization map $\wh{\mc S}^l_N:\;
\Rep_aU_q\agln\to\Rep_a\ur\agln$ is a surjective
homomorphism of rings. We extend the ring structure on
$\Rep_a\ur\agln$ to $\Rep\, \ur\agln$ using the decomposition
\eqref{dec res}.

The weight lattice $P = \Z^N$ of $\gln$ carries the standard partial
order. Namely, for $i=1,\dots,N-1$, set $\al_i=(0,\dots,0,1,-1,0,\dots,0)$
, where $1$ is in the $i$th place. Then for
$\bs\la,\bs\mu\in P$, we have $\bs\la\geq\bs\mu$ if and only if
$\bs\la-\bs\mu=\sum_{i=1}^{N-1}n_i\al_i$, where $n_i\geq 0$ ($i=1,\dots,N-1$).
This induces a partial order
on the set of $N$--tuples of polynomials with constant coefficient
$1$: we say that ${\bs Q} \geq {\bs Q}'$ if $(\deg Q_i)_{i=1,\dots,N} \geq
(\deg Q'_i)_{i=1,\dots,N}$. Proposition 2.5 of \cite{FM2} implies that
\begin{equation}    \label{tri matr}
[V({\bs Q})^\ep] = [V({\bs Q}_\ep)] + \sum_{{\bs Q}'<{\bs
Q}_\ep} m_{{\bs Q}',{\bs Q}_{\ep}} [V({\bs Q}')], \qquad m_{{\bs Q}',{\bs
Q}_{\ep}} \in \Z_{\geq 0},
\end{equation}
where ${\bs Q}_\ep$ is obtained from ${\bs Q}$ by substituting
$q=\ep$. Thus, in general the specialization of an irreducible module
is not irreducible, but the corresponding decomposition matrix is
triangular with respect to the above partial order. 

The following lemma shows that the specializations of the fundamental
representations $V_{\om_i}(a)$ remain irreducible.

\begin{lem}\label{spec fund lem}
$\wh{\mc S}_N^l(V_{\om_i}(a))$ is an irreducible $\ur\agln$--module.
\end{lem}

\begin{proof}
It follows from the explicit formula for $\chi_q(V_{\om_i}(a))$ given
in formula \eqref{ti} that the only dominant monomial in the
specialized character $\chi_\ep(S_N(V_{\om_i}(a)))=\wh{\mc
S}_N^l(\chi_q(V_{\om_i}(a)))$ is the highest weight monomial.
\end{proof}

{}From the isomorphisms \eqref{iso repa}, Lemma \ref{surje} and
formula \Ref{tri matr} we obtain isomorphisms
\begin{equation}    \label{iso repa root}
\Rep_a \ur\agln \simeq
\Z[V_{\omega_i}(a\ep^{2k})]_{i=1,\dots,N}^{k=1,\dots,l}, \qquad \Rep \ur\agln \simeq
\Z[V_{\omega_i}(b)]_{i=1,\dots,N}^{b\in\C^\times}.
\end{equation}

Similarly, we have a map ${\mc
S}_l:\;\Rep\,\uqr\gln\to\Rep\,\ur\gln$.

\subsection{Finite-dimensional representations of the small quantum
affine algebras}\label{fin modules section}

The finite-dimensional representations of $\fu\asln$ are essentially
classified in \cite{BK} (see also \cite{FM2}). We now recall this
classification.


A vector $v$ in a $\fu\asln$--module $V$ is called a highest weight
vector if \be x_{i,n}^+v=0,\qquad {\Phi}^{\pm}_i(u) v
=\Psi_i^\pm(u)v\qquad (i=1,\dots,N-1,\;\; n \in \Z), \ee where
$\Psi_i^\pm(u)\in \C[[u^{\pm 1}]]$ ($i=1,\dots,N-1$).  If in addition
$V=\fu\asln\cdot v$, then $V$ is called a highest weight
representation with highest weight
$\bs\Psi^\pm(u)=(\Psi_i^\pm(u))_{i=1}^{N-1}$.

We will say that a polynomial $P(u) \in \C[u]$ is $l$--{\em acyclic}
if it is not divisible by $(1-au^l)$, for all $a\in\C^\times$.
In other words, the set of roots of 
an acyclic polynomial does not contain subsets of the form $\{
a,a\ep^2,\ldots,a\ep^{2l-2} \}$, $a \in \C^\times$.

According to Theorem 2.6 of \cite{FM2}, every irreducible
finite-dimensional representation of $\fu\asln$ is a highest weight
representation.  The irreducible finite-dimensional representation $V$
with highest weight $\bs\Psi^\pm(u)$ exists if and only
$\Psi_i^\pm(u)$ has the form
\be \pm \ep^{\deg P_i}\frac{P_i(u\ep^{-1})}{P_i(u\ep)},
\ee 
where $P_i(u)$ is an
$l$--acyclic polynomial with constant term $1$ ($i=1,\dots,N$). In
this case all 
common generalized eigenvalues of $\{\Phi_i(u^{\pm1})\}$ on $V$ are of
this form.

\medskip

The description of finite-dimensional representations of $\fu\agln$
is obtained from the description of finite-dimensional representations of
$\fu\asln$ using the decomposition \Ref{det decom
root fin}.

Let ${\mc F}^{\on{fin}}_N$ be the category of all all finite-dimensional
representations $V$ of $\fu\gln$, such that $T_i$ ($i=1,\ldots,N$) act
by semi-simple operators with eigenvalues in $\ep^\Z$. 
We denote the corresponding Grothendieck
group by $\Rep\,\fu\gln$. It is a finite abelian group 
spanned by $[V_{\bs\la}^\ep]$ with $\bs\la=(\la_1,\dots,\la_N)$ such that
$0\leq\la_i-\la_{i+1}<l$ ($i=1,\dots,N$) where $\la_{N+1}=0$. 
Since $\fu\gln$ is a
Hopf algebra, $\Rep\,\fu\gln$ is a ring. However, we will
not use this structure.
Instead we define a different multiplication on
$\Rep\,\fu\agln$ in Section \ref{decomp theorem section}.

Let $\wh{\mc F}^{\on{fin}}_N$ be the category of all finite-dimensional
representation $V$ of $\fu\agln$, such that $V$ is in ${\mc
F}^{\on{fin}}_N$ as a $\fu\gln$--module and
 all generalized eigenvalues of the series
$\Theta_i^\pm(u)$ are of the form \bean\label{form of weight}
\ep^{\deg Q_i}\frac{Q_i(u\ep^{-1})}{Q_i(u\ep)}, \eean where $Q_i(u)$
are polynomials with constant term $1$, such that $Q_i(u)/Q_{i+1}(\ep^2u)$
are $l$--acyclic polynomials. We call such representations {\em polynomial}.

Let $V$ be a representation of $\fu\agln$. A vector $v\in V$ is called
a highest weight vector if 
\be
x_{i,n}^{+} \cdot v=0,\qquad {\Theta^\pm_i(u)}\cdot v =\Omega^\pm_i(u)v
\qquad (i=1,\dots,N,\;\; n\in \Z),
\ee
where $\Omega^\pm_i(u)\in \C[[u^{\pm 1}]]v$. If an addition
$V=\fu\agln\cdot v$, then $V$ is called the highest weight
representation with highest weight
$\bs\Omega^\pm(u)=(\Omega^\pm_i(u))_{i=1,\dots,N}$.
 
The next lemma is obtained from Theorem 2.6 of \cite{FM2} and
decomposition  \Ref{det decom root fin}.

\begin{lem}\label{class gln root factor}
Every irreducible polynomial representation $V$ of $\fu\agln$
is a highest weight representation. The irreducible finite-dimensional
representation with a given highest weight $\bs\Omega^\pm(u)$
exists (and in this case is
unique) if and only if $\Omega_i(u)$ is of the form
\Ref{form of weight}.
\end{lem}

Let $\Rep\,\fu\agln$ be the abelian group spanned by all polynomial
representations of $\fu\agln$.  We will define a ring structure on
$\Rep\,\fu\agln$ in Section \ref{decomp theorem section}.

The restriction of an irreducible $\fu\agln$--representation is an
irreducible $\fu\asln$--module. Similar to Lemma \ref{gl to sl
functor} we obtain

\begin{lem}
Any finite-dimensional representation of $\fu\asln$ may be obtained as
the restriction of a polynomial representation of $\fu\agln$.
\end{lem}

\section{Frobenius homomorphism}\label{frob section}

Let $\ep^*=\ep^{l^2}$. Then $\ep^*\in\{\pm1\}$ and $-\ep^*=(-\ep)^l$.
In this section we construct a ring homomorphism $\on{Fr}^*:\;
\Rep\,\urs\agln\to \Rep\,\ur\agln$. For this we need to introduce an
auxiliary algebra $\dur\agln$, following G. Lusztig \cite{L} (see also
\cite{FM2}).

\subsection{Definition of the Frobenius homomorphism}

Let $P=\Z^N$ be the weight
lattice of $\gln$. For ${\bs\la}\in P$, let $U_{{\bs\la}}$ be the quotient of
$\ur\agln$ by the left ideal generated by 
\be
t_i-\ep^{\la_i},\qquad \bin[t_i;r] -
\bin[{\la_i};r]_\ep,\qquad
\bin[t_jt_{j+1}^{-1};r]-\bin[{\la_j-\la_{j+1}};r]_\ep
\ee 
($i=1,\dots,N$, $j=1,\dots,N-1$, $r\in\Z_{>0}$).
Let $1_{\bs\la}\in U_{\bs\la}$ be the image of $1\in\ur\agln$
in $U_{\bs\la}$. The space $U_{\bs\la}$ is a left $\ur\agln$--module. For
$g\in\ur\agln$, we denote by $g1_{\bs\la}$ the image of $g$ in $U_{\bs\la}$.

Set 
\be 
\dur\agln=\oplus_{{\bs\la}\in P}U_{\bs\la}. 
\ee 
This is a $\C$--algebra with multiplication given by
\bea 
1_{\bs\la}
1_{\bs\mu}=\delta_{{\bs\la},{\bs\mu}}1_{\bs\la},\qquad 
x^{\pm (r)}_{i,n} 1_{\bs\la}=1_{{\bs\la}\pm r\al_i}x^{\pm (r)}_{i,n}, 
\qquad {\mc Q}^\pm_i(u)1_{\bs\la}=1_{\bs\la}{\mc Q}^\pm_i(u)
\eea
($i=1,\dots,N$, $r\in\Z_{>0}$, $n\in\Z$), where 
$\al_i=(0,\dots,0,1,-1,\dots,0)\in P$ with $1$ at the $i$-th place.

A representation of $\dur\agln$ is called unital if $\sum_{{\bs\la}\in
P} 1_{\bs\la}$ acts on it as the identity (note that the infinite sum
$\sum_{\bs\la}1_{\bs\la}$ is a welldefined operator on any
$\dur\agln$--module).

Let $V$ be a finite-dimensional unital representation of
$\dur\agln$. Then using the projectors $1_{\bs\la}$, we obtain a
decomposition $V=\oplus_{{\bs\la}\in P} V_{\bs\la}$. We define the
$\ur\agln$--module structure on $V$ by the rule $g\cdot
v=g1_{\bs\la}\cdot v$ for $g\in\ur\agln$, $v\in
V_{\bs\la}$. Conversely, if $V$ is a type 1 finite-dimensional
$\ur\agln$--module, then we have the weight decomposition
$V=\oplus_{{\bs\la} \in P} V_{\bs\la}$. Define the action of
$\dur\agln$ on $V$ by letting $1_{\bs\la}$ act as the identity on
$V_{\bs\la}$ and by zero on $V_{\bs\mu}$, ${\bs\mu}\neq
{\bs\la}$. Hence we obtain the following lemma (cf. \cite{L}, \S
23.1.4):

\begin{lem}\label{dot cat} The category of finite-dimensional
$\ur\agln$--modules of type 1 is equivalent to the category of
finite-dimensional unital $\dur\agln$--modules.
\end{lem}

A weight ${\bs\la}\in P=\Z^N$ is called $l$--admissible if $\la_i$
($i=1,\dots,N$) are divisible by $l$. We denote the set of all
$l$--admissible weights by $P_l$. For ${\bs\la}\in P_l$ denote by ${\bs\la}/l$
the weight $(\la_1/l,\dots,\la_N/l)$.

Following \cite{L}, Theorem 35.1.9, we define a homomorphism \be \on
{Fr}:\;\dur\agln \to \durs\agln \ee as follows. Note that $\dur\agln$
is generated by $(x^{\pm}_{i,r})^{(n)}1_{\bs\la}$ and $\mc Q^\pm_1(u)1_{\bs\la}$
or, equivalently, by $(x^{\pm}_{i,r})^{(n)}1_{\bs\la}$ and
$g_{1,(r)}1_{\bs\la}$ ($i=1,\dots, N$, $n\in\Z_{>0}$, $r\in\Z$, ${\bs\la}\in
P$).

We set $\on{Fr}((x^{\pm}_{i,r})^{(m)}1_{\bs\la})=0$ and
$\on{Fr}({g_{1,(m)}}1_{\bs\la})=0$ 
if $m$ is not
divisible by $l$ or if ${\bs\la}\not\in P_l$. 
For ${\bs\la}\in P_l$, we set
\bean \label{Fr1}
\on{Fr}((x^{\pm}_{i,r})^{(nl)}1_{\bs\la})&=&\left\{
\begin{matrix} 
(\bar x^{\pm}_{i,r})^{(n)}\bar 1_{\bs\la/l} & (l {\rm \;\; is \;\; odd}),\\
(-1)^{nri}(\bar x^{\pm}_{i,r})^{(n)}\bar 1_{\bs\la/l} & (l {\rm \;\; is
\;\; even}),
\end{matrix}\right.\\ \label{Fr2}
\on{Fr}(g_{1,(nl)}1_{\bs\la})&= &
\left\{
\begin{matrix}
{\bar g_{1,(n)}}\bar 1_{{\bs\la}/l}& (l {\rm \;\; is \;\; odd}),\\
{(-1)^n\bar g_{1,(n)}}\bar 1_{{\bs\la}/l}& (l {\rm \;\; is \;\; even}).
\end{matrix}\right.
\eean
Here and below we put a bar over the elements of the target algebra to
avoid confusion.

Set $u^*=-(-u)^l$. Then for ${\bs\la}\in P_l$, we have
\be
\on{Fr}(\mc Q^\pm_1(u)1_{\bs\la})= \bar{\mc Q}^\pm_1(u^*)\bar 1_{{\bs\la}/l}.
\ee

\begin{lem} The map $\on{Fr}$ extends to an algebra homomorphism.
\end{lem}

\begin{proof}
The restriction of $\on{Fr}$ to $\ur\asln\subset\ur\agln$ coincides
with the Frobenius homomorphism $\ur\asln \to \urs\asln$, as defined
in Lemmas 4.6, 4.7 of \cite{FM2}. The latter is an algebra
homomorphism. It is also clear that the restriction of $\on{Fr}$ to
the subalgebra generated by the elements $g_{1,(r)}1_{\bs\la}$ is a
homomorphism. Thus, we need to check that the map $\on{Fr}$ respects
the relations between $g_{1,(n)}$ and $\ur\asln$.

We have $[g_{1,(n)},(x^{\pm}_{i,s})^{(m)}]=0$ 
($i=2,\dots,N-1$, $s\in\Z$, $n\in\Z\setminus 0$) 
and these relations are preserved by $\on{Fr}$: \be
\on{Fr}([g_{1,(kl)},x^{\pm (nl)}_{i,s}]1_{\bs\la}) = \pm [\bar
g_{1,(k)},\bar x^{\pm (n)}_{i,s}]\bar1_{\bs\la}=0 \qquad (i=2,\dots,N-1), \ee
where ${\bs\la}\in P_l$. Thus, it remains to check that $\on{Fr}$ preserves
the relations between $g_{1,(r)}$ and $(x^{\pm}_{1,s})^{(m)}$. In
particular, we are reduced to the case $N=2$.

So let $N=2$. Due to the triangular decomposition of \lemref{tr dec
gln}, the relations between $g_{1,(r)}$ and $(x^{\pm}_{1,s})^{(m)}$
are generated by the relations \eqref{x and g}.

Let $S$ be the automorphism of $\ur\wh{\mathfrak{gl}}_2$ defined by
\be S(x^{\pm (r)}_{1,n})=x^{\pm (r)}_{1,n\pm1},\qquad
S(g_{i,(n)})=g_{i,(n)},\qquad S(t_i)=t_i, \qquad
S\left(\bin[t_i;r]\right)=\bin[t_i;r] \ee ($n\in\Z$,
$r\in\Z_{>0}$, $i=1,2$). 

Suppose $\on{Fr}$ is well defined. Then
for odd $l$, the automorphism $S$ commutes with
$\on{Fr}$. For even $l$, we have
$\on{Fr}\circ S(y)=(-1)^nS\circ\on{Fr}(y)$ for $y\in\ur\wh{\mathfrak{gl}}_2$ of
degree $nl$ ($n\in\Z$), where the degree is defined by setting
\be \deg(x^{\pm (r)}_{1,n})=nr,\qquad
\deg(g_{i,(n)})=n,\qquad \deg(t_i)=0, \qquad
\deg\left(\bin[t_i;r]\right)=0. \ee

Therefore, applying $S$, it is sufficient to consider only the
specializations of the relations \eqref{x and g} between $g_{1,(r)}$
and $(x^{+}_{1,0})^{(m)}$ (for $(x^{-}_{1,0})^{(m)}$ the proof is the
same):
\be [g_{1,(n)},(x^{+}_{1,0})^{(r)}] =
\frac{q^{-n}}{n[r]_q!}
\sum_{s=0}^r(x^+_{1,0})^sx_{1,n}^+(x^+_{1,0})^{r-1-s} =:
\frac{q^{-n}}{n}B_{r-1,n}.  \ee For $r\in\Z_{\geq 0}$, $n\in\Z$,
define $A_{r,n}$ inductively by \be A_{0,n}=x_n^+,\qquad
A_{r,n}=\frac{A_{r-1,n}x_0^+-q^{2r}x_0^+A_{r-1,n}}{[r+1]_q}.  \ee

It is shown in the proof of Proposition 4.7 in \cite{CP:root} that
$B_{r,n}, A_{r,n} \in \uqr\asln$. Denote by $B^\ep_{r,n}$ and
$A^\ep_{r,n}$ the images of $B_{r,n}$ and $A_{r,n}$ in
$\ur\asln$. Then  we
obtain from the identities (21), (22) of \cite{CP:root}: \bea
A^\ep_{r-1,n} &=& \sum_{s=0}^{r-1}
(-1)^s\ep^{(r-1)(r-s)}(x_{1,1}^+)^{(s)}S(B^\ep_{r-s-1,n-r}),\\
B^\ep_{r-1,n} &=&
\sum_{s=1}^{r}\ep^{(r-s)s}(x_{1,0}^+)^{(r-s)}A^\ep_{s-1,n}. \eea
By induction on $n$, we obtain that
$\on{Fr}(B^\ep_{r-1,n}1_{\bs\la})=\on{Fr}(A^\ep_{r-1,n}1_{\bs\la})=0$ if ${\bs\la}\not\in
P_l$ or $r$ is not divisible by $l$ or $n$ is not divisible by
$l$. For ${\bs\la}\in P_l$ we obtain \bea
\on{Fr}(B^\ep_{rl-1,nl}1_{\bs\la})&=&
l B^{\ep^*}_{r-1,n} 1_{{\bs\la}/l},\\ 
\on{Fr}(A^\ep_{rl-1,nl}1_{\bs\la})&=&
l A^{\ep^*}_{r-1,n} 1_{{\bs\la}/l}
\eea
(at the initial step of our induction, when $n=0$, we have
$B^\ep_{r,0} = (r+1) (x^+_{1,0})^{(r+1)}$, so the above formulas
obviously hold). This concludes the proof of the lemma.
\end{proof}

The homomorphism $\on{Fr}$ is called the {\em quantum Frobenius
homomorphism}.  It induces a map \be \on{Fr}^*:\; \Rep\,\urs\agln\to
\Rep\,\ur\agln, \ee which we call the Frobenius map.

\subsection{Properties of the Frobenius map}

The Frobenius map $\on{Fr}^*$ as defined above clearly commutes with
the left upper corner embeddings from Section \ref{inclusions section
root}. Therefore we have

\begin{lem} The following diagram is commutative:
\be
\begin{CD} \Rep\,\urs\agln
@>{\on{Fr}^*}>>\Rep\,\ur\agln  \\
@VVV @VVV\\
 \Rep\,\urs\wh{\mathfrak{gl}}_{N-1}
@>{\on{Fr}^*}>>\Rep\,\ur\wh{\mathfrak{gl}}_{N-1}
\end{CD}
\ee
\end{lem}

We also have
\begin{lem}    \label{Fr is homomorphism}
The map $\on{Fr}^*$ is a ring homomorphism.
\end{lem}

\begin{proof}
Let $V, W$ be two $\urs\agln$--modules. By Lemmas \ref{mult prop} and
\ref{ep-char of spec}, the eigenvalues of $\bar{\mc Q}^\pm_i(u)$ on $V
\otimes W$ are products of their eigenvalues on $V$ and $W$. Since the
multiplication on $\Rep\,\ur\agln$ is defined via the
$\ep$--character homomorphism (see formula \eqref{product ep}), which
records the generalized eigenvalues of ${\mc Q}^\pm_i(u),
i=1,\ldots,N$, the assertion of the lemma may be reformulated as
follows: the eigenvalues of ${\mc Q}^\pm_i(u)$ on $\on{Fr}^*(V \otimes
W)$ are products of their eigenvalues on $V$ and $W$.

Now observe that the eigenvalues of ${\mc Q}^\pm_i(u)$  ($i=1,\ldots,N$)
are completely determined by the eigenvalues of ${\mc P}^\pm_i(u)$
($i=1,\ldots,N-1$) and ${\mc Q}^\pm_1(u)$.

According to \cite{L}, \S 35.1.10, the restriction of the Frobenius
homomorphism to $\ur\asln$ is compatible with the operation of tensor
products. Applying Lemmas \ref{mult prop} and \ref{ep-char of spec}
again, we obtain that the eigenvalues of ${\mc P}^\pm_i(u)$
($i=1,\ldots,N-1$) satisfy the desired multiplicative property.

Next, we obtain from formula \eqref{Fr2} that the pullback of a
generalized eigenspace of $\bar{\mc Q}^\pm_1(u)$ with eigenvalue
$Q^\pm_1(u)$ is a generalized eigenspace of ${\mc Q}^\pm_1(u)$ with
eigenvalue $Q^\pm_1(u^*)$. Therefore the eigenvalues of ${\mc
Q}^\pm_1(u)$, also satisfy this property. This completes the proof.
\end{proof}

Finally, we show that Frobenius and evaluation maps commute with each
other. For $a\in\C^\times$, set $a^*=-(-a)^l$. In other words
$a^*=a^l$ if $l$ is odd and $a^*=-a^l$ if $l$ is even.

\begin{lem}
The following diagram is commutative:
\be
\begin{CD} \Rep\,\urs\gln
@>{\on{Fr}^*}>>\Rep\,\ur\gln  \\
@V ev^*_{a^*}VV @V ev_a^* VV\\
 \Rep\,\urs\wh{\mathfrak{gl}}_{N}
@>{\on{Fr}^*}>>\Rep\,\ur\wh{\mathfrak{gl}}_{N}
\end{CD}
\ee
\end{lem}

\begin{proof}
The proof is similar to that of \lemref{Fr is homomorphism}. Namely,
we first check the commutativity of this diagram with $\gln$ replaced
by $\sln$. It follows from formula \Ref{ev via braid}, using the fact
that Frobenius map commutes with the braid group action (see \cite{L},
\S 41.1.9). We also use the identities $\ep^l=\ep^*$ for odd $l$ and
$(-\ep)^l=-\ep^*$.

Second, we check commutativity of the analogous diagram for the
generalized eigenvalues of ${\mc Q}^\pm_1(u)$ using formula \Ref{ev
Q}. Combining these results, we obtain the assertion of the lemma.
\end{proof}

\subsection{Frobenius pullback modules}\label{decomp theorem section}

Let $V$ be a ${U}^{\on{res}}_{\ep^*}\agln$--module. The
$\ur\agln$--module $\on{Fr}^*(V)$ is called the {\em Frobenius
pullback} of $V$.

Introduce the following notation 
\be {\bs \La}_{i,a}:=\prod_{j=0}^{l-1}\La_{i,a\ep^{2j}}.
\ee
The monomial ${\bs \La}_{i,a}$ corresponds to the polynomial
$(1-a^lu^l)$.

\begin{lem}\label{char of frob}
Let $V$ be an irreducible representation of
${U}^{\on{res}}_{\ep^*}\agln$ with Drinfeld polynomials $Q_i(u)$,
$i=1,\dots,N$. Then the $i$-th Drinfeld polynomial of $\on{Fr}^*(V)$
is equal to $Q_i(u^*)$.
Moreover, $\chi_\ep(\on{Fr}^*(V))$ is obtained from
$\chi_{\ep^*}(V)$ by replacing $\bar \La_{i,a^l}^{\pm 1}$ with ${\bs
\La}_{i,a}^{\pm 1}$.
\end{lem}

\begin{proof}
In the case of $\ur\asln$ the statement of this lemma is proved in
\cite{FM2}, Theorem 5.7. This gives us the eigenvalues of ${\mc
P}_i^+(u)$ ($i=1,\ldots,N-1$) on $\on{Fr}^*(V)$. Next we compute the
eigenvalues of ${\mc Q}^\pm_1(u)$ using the definition of the
Frobenius map. Combining these two results, we obtain the lemma.
\end{proof}

\begin{cor}
Let $V$ be a polynomial representation of $\urs\agln$. Then
$\on{Fr}^*(V)$ is a polynomial representation of $\ur\agln$.
\end{cor}

For a dominant $N$-tuple of polynomials $\bs Q$, let ${\bs Q}^0$ and ${\bs
Q}^1$ be the dominant $N$-tuples of polynomials such that
\begin{itemize}
\item  $Q_i(u)=Q_i^0(u)Q_i^1(u)$; 
\item  $Q_i^0(u)/Q_{i+1}^0(u)$ is $l$--acyclic;
\item  $Q_i^1(u)/Q_{i+1}^1(u)$ is a polynomial of $u^l$
\end{itemize}
($i=1,\dots,N$), where we set $Q_{i+1}^0(u)=Q_{i+1}^1(u)=1$. 
Clearly, for any dominant $N$-tuple $\bs Q$, the $N$-tuples ${\bs Q}^0$,
${\bs Q}^1$ are uniquely determined.

The next theorem follows from the corresponding result for $\ur\asln$
(see Theorem 9.3 of \cite{CP:root} and Theorem 5.4 of \cite{FM2}).

\begin{prop}\label{decomp thm}
Let $V({\bs Q})$ be an irreducible $\ur\agln$ module with Drinfeld
polynomials ${\bs Q}$. Then $[V({\bs Q})]=[V({\bs Q}^0)][
V({\bs Q}^1)]$ in $\Rep\,\ur\agln$. Moreover, 
$V({\bs Q}^1)$ is the Frobenius pullback of
an irreducible $U^{\on{res}}_{\ep^*}$--module, and $V({\bs Q}^0)$ is
irreducible over $U^{\on{fin}}_\ep \agln$.
\end{prop}

In the case of $\ur\asln$, \propref{decomp thm} is true at the
level of the category of representations. In the case of $\ur\agln$,
it is probably also true at the categorical level. We state it at the
level of Grothendieck rings due to the lack of comultiplication
structure on $\ur\agln$.

Using \propref{decomp thm}, we obtain a natural isomorphism of modules
over $\Z$, \be \Rep\, \ur\agln \simeq \Rep\, U_{\ep^*}^{\on{res}}
\agln\otimes \Rep\, \fu\agln, \ee and therefore the identification \be
\Rep \,\fu\agln\simeq \Rep\, \ur\agln /(\Rep\, U_{\ep^*}^{\on{res}}
\agln)_+, \ee where $(\Rep\, U_{\ep^*}^{\on{res}} \agln)_+$ denotes
the augmentation ideal of $\Rep\, U_{\ep^*}^{\on{res}} \agln$. The
Frobenius map $\on{Fr}^*:\;\Rep\, U_{\ep^*}^{\on{res}}\agln\to \Rep\,
\ur\agln$ is an injective ring homomorphism. Therefore we have a
natural ring structure on the quotient $\Rep\, \ur\agln /(\Rep\,
U_{\ep^*}^{\on{res}} \agln)_+$. We call the induced multiplication on
$\Rep\, \fu\agln$ the {\em factorized tensor product}.

Similarly, we use
\be \Rep\, \ur\gln \simeq \Rep\, U_{\ep^*}^{\on{res}}
\gln\otimes \Rep\, \fu\gln \ee
to define the factorized tensor product structure on $\Rep\,\fu\agln$. 

\section{Hall algebras}\label{hall section}

In this section we describe the Hall algebras of the cyclic quivers
and the infinite linear quiver, and their properties. Most of this can
be found in \cite{Ri,G,VV,Sh}. However, we only need a special case of
the Hall algebra, namely, its specialization at $q=1$ (in the notation
of \cite{Ri}). In this case all the facts we need become quite
elementary. Because of that we give a short self-contained
introduction to Hall algebras in this special case, with the emphasis
on the integral forms of these Hall algebras. Thus, this section may
be read independently of the rest of this paper.

\subsection{Definition of the Hall bialgebra}\label{Hall def section}
Let $\mc A$ be a quiver, i.e., a directed graph. We denote by ${\mc V}$
and ${\mc E}$ the sets of vertices and edges of ${\mc A}$. We will
concentrate on two cases: ${\mc A}=A_l^{(1)}$ and ${\mc
A}=A_\infty$. For $A_l^{(1)}$, we have ${\mc
V}=\{v_0,v_1,\dots,v_{l-1}\}$, ${\mc
E}=\{(v_0,v_1),(v_1,v_2),\dots,(v_{l-1},v_0)\}$. For $A_\infty$, we
have ${\mc V}=\{v_i, i\in\Z\}$, ${\mc E}=\{(v_i,v_{i+1}),
i\in\Z\}$. Therefore $A_l^{(1)}$ is the cyclic quiver with $l$
vertices and $A_\infty$ is the infinite linear quiver.

An ${\mc A}$--{\em set} is by definition a finite set $S$ together
with a decomposition $S=\sqcup_{v\in\mc V}S_v$ and a collection of
maps $\{ e_S \}_{e \in {\mc E}}$.  We call the set $S_v$ the {\em
fiber} of $S$ over $v$. For $e=(v,w)\in\mc E$ the corresponding map
$e_S$ maps $S_{v}$ to $S_w\cup\{\emptyset\}$ in such a way that the
preimage of $t\in S_w$ has at most one element: $\sharp(e^{-1}_St)\leq
1$. When no confusion may arise, we denote $e_S$ simply by $e$.

Two ${\mc A}$--sets $S,L$ are called equivalent if there exist
bijections $b_v:S_v\to L_v$, $v\in\mc V$, intertwining the maps
$e=(v,w)\in\mc E$: $e\circ b_v=b_w\circ e$. We then write $S\simeq L$.

An equivalence class of ${\mc A}$--sets can be represented by the
following picture. We draw the graph $\mc A$, then we position the
elements of the fibers $S_{v}$ over the corresponding vertices and
finally connect the elements $s_v$ and $s_w$ of fibers $S_v$ and $S_w$
by an edge whenever $e(s_v)=s_w$.

Here is an example in the case when ${\mc A}=A_\infty$.
\be
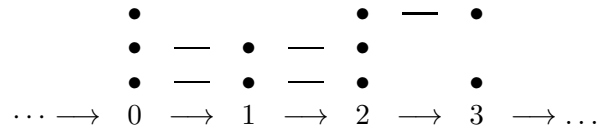
\begin{matrix}
&\bullet&&&&\bullet&\frac{\quad\;}{\quad\;}&\bullet& \\
&\bullet&\frac{\quad\;}{\quad\;}&\bullet&\frac{\quad\;}{\quad\;}&\bullet&&&
\\
&\bullet&\frac{\quad\;}{\quad\;}&\bullet&\frac{\quad\;}{\quad\;}&\bullet&&\bullet&
\\ \dots\longrightarrow& 0& \longrightarrow & 1&\longrightarrow& 2&
\longrightarrow& 3&\longrightarrow \dots
\end{matrix}
\ee
We call the connected components of this picture {\em snakes}. 

An ${\mc A}$--set is called {\em nilpotent} if none of the snakes in
it is a loop. Then each snake has a tail (the vertex in $\mc A$ which is
the projection of the first element in the snake), a head (the
projection of the last element) and the length (number of edges). In the
example above we have two snakes of length 2 with tails at 0, one
snake of length 1 with tail at 2 and two snakes of length 0 with tails
at 0 and 3.

We say that $L$ is an ${\mc A}$--subset of an ${\mc A}$--set $S$, if
all maps $e_S$, $e\in\mc E$, preserve $L$, i.e.,
$e_S(L)\subset(L\cup\emptyset)$. In this case $L$ is clearly an ${\mc
A}$--set itself and we simply write $L\subset S$. The quotient ${\mc
A}$--set $S/L$ is the set $S\setminus L$ and for $s\in S\setminus L$,
such that $e_S(s)=t$ we let $e_{S/L}(s)=t$ if $t\not\in L$, and we let
$e_{S/L}(s)=\emptyset$ otherwise.

In other words, snakes of an ${\mc A}$--subset $L$ of $S$ are subsets
of some snakes of $S$, such that all heads of snakes of $L$ are also
heads of snakes of $S$. Then the snakes of $S/L$ are complimentary
parts of the snakes of $L$; in particular all tails of snakes of $S/L$
are tails of snakes in $S$. Snakes in $S$ go to snakes either in $L$
or in $S/L$ and some snakes in $S$ get ``cut'' into two snakes, one in
$L$ and one in $S/L$.

We say that $S=K\oplus L$ if both $L$ and $K$ are ${\mc A}$--subsets
of $S$ and $L \simeq S/K$, $K \simeq S/L$. In this case the snakes of
$S$ are divided into two groups between $L$ and $K$, and no cutting
occurs.

Let $H=H(\mc A)$ be the set of equivalence classes of nilpotent ${\mc
A}$--sets. Let $\mc H=\mc H(\mc A)$ be the free $\Z$--module with a
set of generators labeled by elements of $H$. Abusing terminology, we
call the generators ${\mc A}$--sets and denote them in the same way as
the corresponding ${\mc A}$--sets. In particular, the empty set
corresponds to the unit element of ${\mc H}$.

We extend the operation of the direct sum on $H$ to ${\mc H}$ by
linearity.

For ${\mc A}$--sets $K,L\in H$, define their product by the formula
\be m(K,L)=K\cdot L=\sum_{S\in H}g_{K,L}^S\;S, \qquad g_{K,L}^S=\sharp\{M
\subset S,\;\;M\simeq K,\;\; S/M\simeq L\}.  \ee

For an ${\mc A}$--set $S\in H$, define its comultiplication by  
\be
\Delta(S)=\sum_{K,L\in H,\;\; S\simeq K\oplus L}K\otimes L.
\ee
We extend both $m$ and $\Delta$ to $\mc H$ by linearity.

\begin{lem} $({\mc H},m,\Delta)$ is a cocommutative bialgebra.
\end{lem}
\begin{proof}
Associativity, coassociativity and cocommutativity are clear. We
only have to check that comultiplication is the algebra homomorphism,
$\Delta(K\cdot L)=\Delta(K)\Delta(L)$. 

Note that if $S,S_1,S_2$ are ${\mc A}$--sets, then the element
$S_1\otimes S_2$ is present in $\Delta(S)$ if and only if $S=S_1\oplus
S_2$.  Clearly, $S_1\otimes S_2$ is present in $\Delta(K)\Delta(L)$ if
and only if $S_1\oplus S_2$ is present in $K \cdot L$. Moreover, the
equality of corresponding coefficients follows from the binomial
identity $$\left( \begin{matrix} a\\ c \end{matrix} \right) =
\sum_{c_1+c_2=c} \left( \begin{matrix} a\\ c_1 \end{matrix} \right)
\left( \begin{matrix} b\\ c_2  \end{matrix} \right).$$
\end{proof}

We show later that in the cases of $A_\infty$ and $A_l^{(1)}$, the
bialgebra $\mc H$ admits an antipode and is therefore a Hopf
algebra. So in what follows we use the term Hopf algebra rather than
bialgebra.

The Hall algebra $\mc H$ has a natural $\Z_{\geq 0}^{\mc
V}$--gradation. Namely, for $S \in H$, we set $\deg S = (\# S_v)_{v
\in {\mc V}}$. Clearly, this grading is compatible with multiplication
and comultiplication. 

For ${\bs d}\in Z_{\geq 0}^{\mc V}$ we denote by $f_{\bs d}$ the only
${\mc A}$--set of degree ${\bs d}$ all snakes of which have length 0.

\subsection{The Hall algebra associated to $A_\infty$}

Let ${\mc H}_\infty$ be the Hall algebra corresponding to $A_\infty$.
Recall that the vertices of $A_\infty$ are labeled by $\Z$.  We write
simply $f_i^{(j)}$ for the element $f_{\bs d}$ with ${\bs
d}=(\dots,0,0,j,0,0,\dots)$, where $j$ is put at the $i$-th place.
We obviously have $j!f_i^{(j)}=f_i^j$.

\begin{lem}
There exists an isomorphism of Hopf algebras 
\bea
U_\Z\sli^-\overset{\sim}\to\mc H_\infty\qquad f_i^{(j)}\mapsto f_i^{(j)}.
\eea
\end{lem}
\begin{proof}
Consider the homomorphism $U \sli^- \to\mc H_\infty \otimes \C$
sending the generators $f_i, i \in \Z$ to $f_i$. The verification of
the Serre relations is straightforward, so this homomorphism is
welldefined. Moreover, it is clearly compatible with comultiplication
on the generators of $U\sli^-$. Therefore it is a homomorphism of Hopf
algebras.

Note that the image of root vectors $f_{i,k}$, $i>k$ (see \Ref{root
vectors}) is a single snake of length $i-k-1$ with head at $i-1$ and
tail at $k$. Also note that if $K,L$ are two ${\mc A}$--sets such that
all heads of $L$ are located to the right of all heads of $L$, then $K
\cdot L=K\oplus L$. Therefore, the image of an element of the PBW
$\Z$--basis \eqref{pbv vector} of products of divided powers is an
${\mc A}$--set consisting of snakes corresponding to elements
$f_{i_s,j_s}$. It follows that the map in the lemma is surjective and
injective.
\end{proof}

Note that in fact we proved that the PBW $\Z$--basis \eqref{pbv
vector} corresponds to the $\Z$--basis of ${\mc A}$--sets in $\mc
H_\infty$.

\subsection{The Hall algebra associated to $A^{(1)}_l$}
\label{Hall l}

Let $\mc H_l$ be the Hall algebra corresponding to $A_l^{(1)}$. Recall
that we label the vertices of $A_l^{(1)}$ by $0,1,\dots,l-1$. We write
simply $f_i^{(j)}$ for $f_{{\bs d}}$ where ${\bs
d}=(0,\dots,0,j,0,\dots,0)$ and $j$ is in place $i$. We also write
$x_i$ for $f_{\bs d}$, where $\bs d=(i,i,\dots,i), i\geq 0$.

Let $U$ be the $\Z$--subalgebra of $\mc H_l$ generated by $f_i^{(j)}$
($i=0,\dots,l-1, \; j\in\Z_{>0})$ and let $W$ be the $\Z$--subalgebra
of $\mc H_l$ generated by $x_i$ ($i\in\Z_{>0})$.

The subalgebra $U$ is a Hopf subalgebra of $\mc H_l$, but the algebra
$W$ is not. The algebra $W$ is clearly commutative. We will prove that
$W = \Z[x_i]_{i \in \Z_{>0}}$.

Let $\asll^-$ be the Lie algebra generated by $f_i$
($i=0,\dots,l-1$) subject to the Serre relations 
\bea
&[f_i,[f_i,f_{j}]]=0 \qquad  & (i,j {\rm \;\; are \;\; adjacent\;\;
in}\;\; {A_l^{(1)}}),\\
&\;[f_i,f_j]=0 \qquad  & (i,j{\rm \;\; are \;\;not\;\; adjacent\;\;
in}\;\; {A_l^{(1)}}).
\eea 
The universal enveloping algebra $U\asll^-$ is a
cocommutative Hopf algebra generated by $f_i$ ($i=0,\dots,l-1$)
with comultiplication given by
\be
\Delta(f_i)=f_i\otimes 1+ 1\otimes f_i.
\ee
and the antipode
\be
S(f_i)=-f_i.
\ee
The algebra $U\asll^-$ has an integral form $U_\Z\asll^-$ generated by
the divided powers $f_i^{(j)}$ ($i=0,\dots,l-1; j \in \Z_{>0}$), which
is a Hopf $\Z$--subalgebra of $U\asll^-$.
We have a surjective homomorphism of Hopf algebras $U_\Z\asll^-\to
U\subset\mc H_l$ mapping $f_i^{(j)}$ to $f_i^{(j)}$. We will prove
that this homomorphism is an isomorphism.

Let $\mc L$ be the $\Z$--subalgebra of $\mc H_l$ generated by
$f_i^{(j)}$ ($i=0,\dots, l-1, \; j\in\Z_{>0})$ and $x_j$
($j\in\Z_{>0}$).  We will prove that in fact $\mc L=\mc H_l$.

\begin{lem}\label{all gen} The Hall algebra $\mc H_l$ is generated by
the elements $f_{\bs d}$ (${\bs d}\in\Z_{>0}^{l}$).
\end{lem}
\begin{proof}
Define a partial order on the set of all ${\mc A}$--sets. Let $K,L$ be
${\mc A}$--sets with snakes of lengths $a_1\geq a_2\geq\dots\geq a_k$ and
$b_1\geq b_2\geq \dots b_l$. We say $K>L$ if $a_1>b_1$ or $a_1=b_1$,
$a_2>b_2$ or $a_1=b_1$, $a_2=b_2$, $a_3>b_3$, etc. 

Let the head of the longest snake of $K$ be located at 1 and $a_1>1$. 
Let we have $m$ snakes of length $a_1$ with heads at $1$.
Let $K'$ be
the same as $K$ except that all $m$ longest snakes with heads at $1$ are 
cut into two snakes of length $a_1-1$ and length 0 (the zero length
snakes are at 1). Then $K>K'$.

We have \be f_i^{(m)}K'=K+..., \ee where the dots stand for the terms
smaller than $K$. Therefore, every ${\mc A}$--set can be reduced to
product of ${\mc A}$--sets which have only snakes of length 0.
\end{proof}

\begin{lem}\label{gen lem}
The Hall algebra $\mc H_l$ is generated by $f_i^{(j)}$
($i=0,1\dots,l-1, \; j\in\Z_{>0}$) and $x_j$ ($j\in\Z_{>0}$), i.e.,
$\mc H_l=\mc L$.
\end{lem}

\begin{proof}
By Lemma \ref{all gen}, it is sufficient to express $f_{\bs d}$ in
terms of $f_i$ and $x_j$. The set $\Z_{>0}^{\mc V}=\Z_{>0}^{l}$ has a
natural partial order: ${\bs d_1}\geq {\bs d_2}$ if all components of
${\bs d_1}$ are greater than or equal to the corresponding components
of ${\bs d_2}$.

Note that the homogeneous component of ${\mc H}_l$ of degree ${\bs d}$
is generated by $f_{\bs d'}$ with ${\bs d'}\leq{\bs d}$.

Fix ${\bs d}$ and suppose $f_{\bs d}\neq x_j$. Then without loss of
generality we may assume that $d_1>d_0$ and that $f_{{\bs d}'}\in\mc
L$ for all ${\bs d}'<{\bs d}$. Then it is sufficient to show that
$f_{\bs d}\in{\mc L}$.

In particular our assumption implies that
all $\mc A$-sets of degree ${\bs d}'<{\bs d}$ are in $\mc L$.

For $s=1,\dots,d_0$, denote $y_s$
the ${\mc A}$--set of degree ${\bs d}$ 
which has $s$ snakes of length 1 with tails at 1 and all other
snakes are 0. Let $y_s'$ be the ${\mc A}$--set of degree
${d_0,s,d_2,d_3,\dots}$ which has $s$ snakes of length $1$ with head
at 1 and all other snakes of length 0. We note that by our assumption 
$y_s'\in\mc L$ and compute
\be
f_1^{(d_1-s)}y_s'=y_s+\sum_{i=s+1}^{d_0}a_{is}y_i,
\ee
for some nonnegative integers $a_{is}$. Therefore $y_s\in\mc L$. In
particular $y_0=f_{\bs d}$ is in $\mc L$.
\end{proof}

\begin{lem} 
The multiplication map $U\otimes W\to {\mc H}_l$ is surjective.
\end{lem}

\begin{proof}
We need to show that elements of the form $u w$, where $u \in U$ and
$w \in W$, span ${\mc H}_l$ over $\Z$. We denote the span of elements
of this form by $U \cdot W$. Because of \lemref{gen lem}, it suffices
to show that all products $x_n f_j^{(i)}$ belong to $U \cdot W$. The
Hall algebra ${\mc H}_l$ has an automorphism induced by the rotation
of the quiver $A_l^{(1)}$. This automorphism sends $f_j^{(i)}$ to
$f_{j+1}^{(i)}$ and $x_n$ to itself. Therefore we only need to show
that $x_nf_1^{(i)}\in U\cdot W$ for all $i\in\Z_{>0}, n\in\Z_{\geq
0}$. We will show this by induction on $n$. For $n=0$ the statement is
trivial. We will prove the induction step by another induction on $i$.

Let $\bs d=(n,n+i,n,n,\dots$ be the degree of $x_nf_1^{(i)}$. By our
assumptions all $\mc A$-sets of degree smaller then $\bs d$ are in
$U\cdot W$. Then following the proof of Lemma \ref{gen lem} we obtain
that $f_{\bs d}\in U\cdot W$.

Let $y_s$, $s=0,\dots,n$ be the 
$\mc A$-sets of degree $\bs d$ which have $s$ snakes of length
1 with head at 2 and all other snakes of lengths 0. Then
$x_nf_1^{(i)}$ is a linear combination of $y_s$ and it is enough to
prove that $y_s\in U\cdot W$. 

Let $y_s'$ be be the 
$\mc A$-sets of degree $(n,n+i,n-s,n,n,\dots)$ with all snakes of
length $0$. We have $y's\in U\subset U\cdot W$.

We compute
\be
f_2^{(s)}y_s'=y_s+\sum_{j=0}^{s-1}a_{js}y_j
\ee
for some nonnegative integers $a_{js}$. 
We already showed that $y_0=f_{\bs d}$
is in $U\cdot W$, therefore by induction on $s$ all $y_s$ are in
$U\cdot W$.
\end{proof}

\begin{lem}
The map $U_\Z\asll^-\otimes\Z[x_i]_{i\in\Z_>0}\to \mc H_l$
is an isomorphism of $\Z$--modules.
\end{lem}

\begin{proof}
We have already proved that this map is surjective. Hence it is
sufficient to compare (graded) ranks of both $\Z$--modules.

For an ${\mc A}$--set $S$ of degree ${\bs d}$ define the total degree
by the formula $d_0+d_1+\dots+d_{l-1}$. Set the degree of each $f_i\in
U\asll^-$ to be 1 and the degree of $x_i$ to be $li$. Then the map in
lemma is degree preserving. Moreover, the homogeneous $\Z$--modules of
a given degree $d$, $({\mc H_n})_d$, $(U\asll^-)_d$ and
$(\Z[x_i]_{i\in\Z_>0})_d$, all have finite rank, so we compute the
corresponding formal characters: \bea {\rm ch}({\mc
H_l})(\xi)&:=&\sum_{d=0}^\infty ({\mc
H_l})_d\xi^d=\prod_{n=1}^\infty\frac{1}{(1-\xi^n)^{l}},\\ {\rm
ch}(\Z[x_i]_{i\in\Z_{>0} })(\xi)&:=&\sum_{d=0}^\infty
(\Z[x_i]_{i\in\Z_{>0}
})_d\xi^d=\prod_{n=1}^\infty\frac{1}{(1-\xi^{nl})},\\ {\rm
ch}(U_\Z\asll^-)(\xi)&:=&\sum_{d=0}^\infty
(U\asll^-)_d\xi^d=\prod_{n=1}^\infty\frac{(1-\xi^{nl})}{(1-\xi^n)^{l}}.
\eea The first two formulas are obvious and (the complexification of)
the last one is well-known (e.g., it follows from \cite{K}, Lemma
14.2). Therefore the lemma follows.
\end{proof}

\subsection{The center of the Hall algebra associated to $A^{(1)}_l$}

For $i\in\Z_{>0}$, let $z_i\in\mc H_l$ be the sum of $l$ distinct $\mc
A$--sets which consist of a single snake of length $li-1$. The degree
of $z_i$ is $(i,i,\dots,i)$ and \be \Delta(z_i)=z_i\otimes 1 +1\otimes
z_i.  \ee


Let $p_i\in\mc H_l$ be the sum of $l^i$ distinct $\mc A$-sets which
consist of $i$ snakes of length $l-1$. We set $p_0=1$. Then the 
degree of $p_i$ is $(i,i\dots,i)$ and we have
\be
\Delta(p_i)=\sum_{j=0}^i p_j\otimes p_{i-j}.
\ee

\begin{lem}
We have the identity of formal power series in $u$:
\begin{equation}    \label{pn}
\sum_{n=0}^\infty p_n (-u)^n = 
\exp\left(-\sum_{j=1}^\infty z_ju^j/j\right).
\end{equation}
\end{lem}

\begin{proof}
The following formula is checked by direct computation:
\bean\label{center iden}
(s+1)p_{s+1}=\sum_{i=0}^s (-1)^i p_{s-i}z_{i+1}.
\eean

It follows that 
\be
\sum_{s=0}^\infty (s+1)p_{s+1}(-u)^s=\left(\sum_{i=0}^\infty
p_i(-u)^i\right) \left(\sum_{j=1}^\infty -z_ju^{j-1}\right).
\ee
The lemma is obtained by integrating this identity.
\end{proof}

Denote by $Z_l$ the $\Z$--algebra generated by $p_i$ ($i\in\Z_{>0}$).
The proof of the following lemma is straightforward using the identity
\Ref{center iden}.

\begin{lem}
The algebra $Z_l$ is a commutative and cocommutative polynomial Hopf
algebra $\Z[p_i]_{i\in\Z_{>0}}$ which contains $z_j$ $(j\in\Z_{>0})$.
Moreover, $Z_l$ is in the center of $\mc H_l$.
\end{lem}

\begin{prop}    \label{tensor product}
We have the following tensor product
decomposition of Hopf algebras over $\Zl$ \bean \label{decom of Hall}
{\mc H_l} \underset{\Z}\otimes \Zl \simeq (U_\Z\asll^-
\underset{\Z}\otimes Z_l) \underset{\Z}\otimes \Zl. \eean
\end{prop}

\begin{proof}
Denote by $B$ the Hopf algebra appearing in the right hand side of
this formula. It is sufficient to prove that each $x_i$ belongs to
$B$. We will do it by induction on $i$.  Suppose that $x_j\in B$ for
all $j<i$. Then in particular all $\mc A$--sets of degree
$(d_0,d_1,\dots)$ such that $d_j\leq i$ and this inequality is not an
equality for at least one vertex belong to $B$.

Let $y_{j,s}$ ($j=0,\dots,l-1,\; s=0,\dots,l-1$) be the $\mc A$-set of
the same degree as $x_i$ which has $i$ snakes of length $s$ with tail
at $j$ with all other snakes of length $0$. In particular, we have
$y_{j,0}=x_i$.

First we claim that for all $j,s>0$ we have $y_{j,s}=b_s
y_{j,s+1}+\dots$ for some integers $b_s$ (here and below the dots
denote terms in in $B$). Indeed, let $y_{j,s,t}$ where $t=0,\dots,i$,
be the $\mc A$-set of the same degree as $x_i$ which has the following
snakes of positive length: $t$ snakes of length $s+1$ and $i-t$ snakes
of of length $s$ all with tails at $j$. Let $y_{j,s,t}'$ be the $\mc
A$-set obtained from $y_{j,s,t}$ obtaining by removing all snakes of
length 0 at vertex $j+s+1$. We have $y_{j,s,t}'\in B$ if $t<i$. Now we
compute \be y_{j,s,t}'f_{j+s+1}^{(i-t)}=y_{j,s,t}+\sum_{p>t}a_{jstp}
y_{j,s,p}, \ee for some nonnegative integers $a_{jstp}$. Therefore we
can compute $y_{j,s}=y_{j,s,0}$ via $y_{j,s+1}=y_{j,s,i}$ modulo terms
in $B$ and the claim follows.

Thus, we have $lx_i=a\sum_{j=0}^{l-1}y_{j,l-1}+\dots$ for some integer
$a$. Consider $z=p_i\in \mc Z_l\subset B$.
In particular, the degree of $z$ is the same
as the degree of $x_i$. To complete the proof it is sufficient to show
that $z-\sum_{j=0}^{l-1}y_{j,l-1}$ is in $B$.

The difference $z-\sum_{j=0}^{l-1}y_{j,l-1}$ is a sum of $\mc A$-sets
with positive integer coefficients. Each such $\mc A$-set $y$ has $i$
snakes of length $l-1$. Moreover, not all of the snakes of $y$ have the
same tail. Suppose that $y$ has exactly $s$ snakes with tail at
$0$. Then we can write $y=y_0\oplus y'$ where $y_0$ has exactly $s$
snakes of length $l-1$ with tail at 0. In particular $y_0,y'\in B$.
Note that no snake of $y'$ has a head at $l-1$.  Therefore we have
$y=y_0y'\in B$. This completes the proof.
\end{proof}

Note that since the center of the algebra $U_\Z\asll^-$ is spanned by
the unit element, $Z_l$ is the center of $\mc H_l$.

Also note that a linear combination of $\mc A$--sets consisting of a
single snake of length $li-1$ with the sum of all coefficients equal
to zero is necessarily in $U_\Z\asll^-$.

\subsection{The unwinding map}\label{unwind}

For $n\in\Z_{\geq 0}$, let $F_n\subset\mc H_\infty$ be the
$\Z$--submodule spanned by all ${\mc A}$--sets of degree ${\bs d}$,
such that $d_i=0$ if $|i|<n$. Each $F_n$ may be identified with a
quotient of ${\mc H}_\infty$, and we have natural surjective
homomorphisms $F_{n+1} \to F_n$. Let $\wt{\mc H}_\infty$ be the
inverse limit of the projective system $\{ F_n \}_{n\geq 0}$. Then
$\wt{\mc H}_\infty$ is a Hopf algebra, and ${\mc H}_\infty$ is its
Hopf subalgebra.

For each ${\bs d} \in \Z_{\geq 0}^l$ and $p \in \Z$ we
define ${\bs d}(p) \in \Z_{\geq 0}^{\Z}$ by the formulas ${\bs
d}(p)_{i+p} = {\bs d}_i$  ($i=0,\ldots,l-1$), and ${\bs
d}(p)_k = 0$ if $k \not\in \{ p,p+1,\ldots,p+l-1 \}$.

\begin{lem}\label{unwind lem}
There exists a unique embedding of Hopf algebras \bea w:\; \mc
H_l\to\wt{\mc H}_\infty,\qquad f_{\bs d} \mapsto \sum_{k\in\Z}
f_{\bs d(kl)}, \eea
which commutes with the operation of taking direct sum,
i.e., $w(K\oplus L)=w(K)\oplus w(L)$.
\end{lem} 

\begin{proof}
Injectivity and compatibility with comultiplication is clear. We only
have to show the equality $w(K \cdot L)=w(K)w(L)$. It is readily seen
by induction on the number of snakes in the ${\mc A}$--set
$K$. Indeed, suppose $K$ has a snake of length $l$ and head $s$.  Then
the coefficient of some $M$ in product $K \cdot L$ can be computed in
the following way.

Denote by $S_{s,k}$ the snake of length $k$ with the head $s$. First,
we fix the length $a\geq l$ of the snake in $M$ with head $s$, such
that its complement of $S$ is a snake in $L$.  We set $M'(a)$, $K'$,
$L'(a)$ be such that $M'(a)\oplus S_{s,a}=M$, $K'\oplus S_{s,l}=K$ and
$L'(a)\oplus S_{s-l-1,a-l-1}=L$.

Then the coefficient of $M$ in $K \cdot L$ is $\sum_a
k_ag_{K',L'(a)}^{M'(a)}$, where $k_a$ is the number of snakes in $M$
equivalent to $S_{s,a}$.

Now consider a term in $w(M)$. There are $k_a$ ways to represent it in
the form (a term in $w(M'(a))$)$\oplus$(a term of $w(S_{s,a})$). Each
such ${\mc A}$--set will give the contribution $g_{K',L'(a)}^{M'(a)}$
to the coefficient of $M$ in $K \cdot L$, by the induction hypothesis.
\end{proof}

We refer to the map $w$ as the {\em unwinding} map.

\subsection{Hopf algebras $\Z[\wGLl]$ and
$\Z[\wGLln]$}\label{dual hall root}

By Lemma \ref{two dual}, the algebra $\mc H_\infty \simeq U_\Z\sli^-$ is
(restricted) dual to the polynomial algebra $\Z[\wSLi]$ generated by
the entries of a lower triangular matrix $M$. It is clear that the
non-degenerate pairing between $\mc H_\infty$ and $\Z[\wSLi]$ may be
extended to a non-degenerate pairing between $\wt{\mc H}_\infty$ and
$\Z[\wSLi]$. Now we use the unwinding map to find a realization of the
Hopf algebra dual to $\mc H_l$. The algebra $\mc H_l$ is a subalgebra
in $\wt{\mc H}_\infty$, so the dual algebra is a quotient of
$\Z[\wSLi]$.

For $N \in \Z_{>0}$, let $\wGLln$ be the proalgebraic group of
matrices $M = (M_{i,j})_{i,j \in \Z}$ satisfying the conditions
$M_{i,i} = 1$ for all $i$, $M_{i,j} = 0$ whenever $i<j$ or $i>j+N$,
and $M_{i+l,j+l} = M_{i,j}$. We define a product of two such matrices
as their usual product, in which all $(i,j)$ entries with $i>j+N$ are
set to be equal to $0$. It is clear that this operation makes $\wGLln$
into a prounipotent proalgebraic group over $\Z$. Furthermore, we
have natural surjective group homomorphisms $\wt{GL}^{(N+1),-}_l
\twoheadrightarrow \wGLln$. Let $$\wGLl =
\underset{\longleftarrow}{\lim} \; \wGLln.
$$
This is also a prounipotent proalgebraic group of $l$--periodic
matrices $M = (M_{i,j})_{i,j \in \Z}$, such that $M_{i,i} = 1$ for all
$i$, $M_{i,j} = 0$ for all $i<j$, and $M_{i+l,j+l} = M_{i,j}$.

The algebra $\Z[\wGLln]$ of regular functions on $\wGLln$ is by
definition the algebra of commutative polynomials in variables
$M_{i,j}$ $(i,j\in\Z; i>j)$ subject to the relations
$M_{j+l,i+l}=M_{j,i}$. The algebra $\Z[\wGLl]$ of functions on $\wGLl$
is by definition the algebra of polynomials in commuting variables
$M_{i,j}$ $(i,j\in\Z; i>j)$ subject to the relations
$M_{j+l,i+l}=M_{j,i}$.

Both are Hopf algebras, with the comultiplication $\Delta$ and the
antipode $S$ given by the formulas \Ref{delta of M} and \Ref{S of M}.
Furthermore, Hopf algebra $\Z[\wGLl]$ is the inductive limit of its
Hopf subalgebras $\Z[\wGLln]$ ($N \in \Z_{>0}$).

We have the surjective homomorphisms of Hopf algebras
\bea
w^*: \;\Z[\wSLi] &\to& \Z[\wGLl], \qquad M_{i,j}\mapsto
M_{i,j}\qquad (i,j\in\Z).
\eea
Recall the unwinding homomorphism $w$ of Hopf algebras given by Lemma
\ref{unwind lem}.

\begin{prop}    \label{def of pairing for Hl}
The pairing of $\wt{\mc H}_\infty$ and $\Z[\wSLi]$ induces a
non-degenerate pairing of bialgebras $\mc H_l$ and $\Z[\wGLl]$ with
respect to which the maps $w$ and $w^*$ are dual to each other. The
monomial basis in the polynomial ring $\Z[\wGLl]$ is dual to the basis
of ${\mc A}$--sets in $\mc H_l$ under this pairing, so that $\mc H_l$
and $\Z[\wGLl]$ are restricted dual of each other.
\end{prop}

\begin{proof}
The elements $M_{i+l,j+l}-M_{i,j}$ are clearly orthogonal to $\mc
H_l$. Now the lemma follows by dimension counting. Namely, we set the
degree of $M_{i,j}$ to be $i-j$. The pairing clearly respects degrees,
and the dimensions of the corresponding homogeneous components are
equal.
\end{proof}

We call the map $w^*$ the {\em winding} map.

In particular, we obtain an antipode in $\mc H_l$ from the antipode
in $\Z[\wGLl]$.

\begin{remark}

Let ${\mc A}$ denote the quiver $A_\infty$ or $A^{(1)}_l$. Call an $\mc
A$--set {\em $N$--nilpotent} if it contains no snakes of
length greater than $N$. Denote by $H^{(N)}$ the set of all
$N$--nilpotent $\mc A$--sets. The vector space $\mc H^{(N)}$ spanned
by $H^{(N)}$ is a Hopf algebra with multiplication and
comultiplication defined by the formulas 
\bea K\cdot L=\sum_{S\in H^N} g_{K,L}^S
\qquad \Delta(S)=\sum_{K,L\in H^{(N)},\;S=K\oplus L}K\otimes L.  \eea
Moreover, we have an obvious surjective homomorphism of Hopf algebras
\be \pi_N:\;\mc H^{(N+1)}\to \mc H^{(N)}, \ee with the kernel spanned by
all $\mc A$--sets which have a snake of length $N+1$.

The non-degenerate pairing between the Hall algebra ${\mc H}_\infty$
(resp., ${\mc H}_l$) and $\Z[\wSLi]$ (resp., $\Z[\wGLl]$) gives
rise to non-degenerate Hopf algebra pairings
\bea
{\mc H}^{(N)}_\infty \otimes \Z[\wSLin] \to\Z, \qquad
{\mc H}^{(N)}_l \otimes \Z[\wGLln] \to \Z.
\eea
\end{remark}

\subsection{The groups $\wSLl$ and $\wt{\mc Z}_l^-$}    \label{new
sect}

We realize the group $\wGLl$ as a subgroup of loop group of $GL_l$ as
follows. For a ring $R$, the group of $R$--points of $\wGLl$ is the
group of $l \times l$ matrices with coefficients in $R[[u^{-1}]]$ such
that the constant part is a lower triangular matrix with $1$'s on the
diagonal. Denote by $m_{i,j}(t)$ the $(i,j)$ entry of such a matrix,
and by $m_{i,j}(n)$ its $n$th coefficient, where $n\leq 0$. Then under
this identification $m_{i,j}(n)$ corresponds to the entry $M_{i,j+nl}$.

In this realization we define the determinant of $M$, denoted $\on{det}
M$, as the determinant of the corresponding $l \times l$ matrix of
power series in $u^{-1}$. We have
$$
\on{det} M = 1 + a_1 u^{-1} + a_2 u^{-2} + \ldots,
$$
where each $a_i (i>0)$ is a polynomial in the generators $m_{i,j}(n)$
(or equivalently, $M_{s,k}$) of $\wGLl$.

It is straightforward to find the coproduct of the elements $a_i$. If
we introduce the generating function $a(z) = 1 + \sum_{i>0} a_i z^i$,
then
$$
\Delta(a(z)) = a(z) \otimes a(z).
$$
In other words,
\begin{equation}    \label{delta of ai}
\Delta(a_i) = \sum_{j=0}^i a_j \otimes a_{i-j},
\end{equation}
where we set $a_0 = 1$. Therefore the $\Z[a_i]_{i>0}$ is a Hopf
subalgebra of $\Z[\wGLl]$. The spectrum of this subalgebra is a
commutative proalgebraic group, which we denote by $\wt{\mc Z}_l^-$.

We define the group $\wSLl$ as the subgroup of $\wGLl$ determined by
the equations $a_i=0$ ($i>0$), i.e. as the group of matrices of
determinant $1$.  Thus, $\Z[\wSLl]$ is the quotient of $\Z[\wGLl]$ by
the augmentation ideal of $\Z[\wt{\mc Z}_l^-] \subset \Z[\wGLl]$.

\subsection{The Hall algebra ${\mc H}_l$ and the integral form of $U
  \wt{\gw}_l^-$}    \label{on decomp}

We now consider the Lie algebras of the groups introduced in the
previous section.

Let $\wt{\gw}_l^-$ (resp. $\widehat\gw_l^-$) be the Lie subalgebra of
the Lie algebra $\gll[[u^{-1}]]$ (resp. of $\gll[u^{-1}]$) consisting
of all elements $a_0 + a_1 u^{-1} + \ldots$, such that $a_0$ is a
strictly lower triangular $l \times l$ matrix.

Any element of $\wt{\gw}_l^-$ defines an endomorphism of the vector
space of Laurent power series $\C^l((u^{-1}))$. 
Let $\ol{v}_i$ 
($i=0,\ldots,l-1$) be the standard basis of $\C^l$, and consider the
topological basis $\{ v_j \}_{j \in \Z}$ of $\C^l((u^{-1}))$, 
where $v_{i+kl} = \ol{v}_iu^k$ ($i=0,\dots,l-1$, $k\in\Z$). Written
with respect to this basis, each 
element of $\wt{\gw}_l^-$ becomes a strictly lower triangular
matrix. Thus we obtain an embedding of Lie algebras $j_l:
\wt{\gw}_l^- \hookrightarrow \wt\sw_\infty^-$.

For any $A \in \wt{\gw}_l^-$, the matrix
$\exp(j_l(A))$ is a well-defined element of $\wGLl$. Moreover, $\wt\gw_l^-$
is the Lie algebra of the group $\wGLl$, and the exponential map $\exp:
\wt\gw_l^- \to \wGLl$ is an isomorphism.

The Lie algebra $\gll$ is the direct sum of its Lie
subalgebras: $\gll = \sll \oplus {\zz}$, where ${\zz}$ consists
of scalar matrices. Hence we obtain direct sum decompositions
\begin{equation}    \label{dir sum of la}
\wt\gw_l^- = \wt\sw_l^- \oplus \wt\zz_l^-,\qquad \widehat\gw_l^- =
\widehat\sw_l^- \oplus \widehat\zz_l^-,
\end{equation}
where $\wt\sw_l^-$ (resp. $\widehat\sw_l^-$) is a Lie subalgebra of
$\sll[[u^{-1}]]$ (resp. $\sll[u^{-1}]$) defined in the same way as
$\wt\gw_l^-$ (resp. $\widehat\gw_l^-$), and $\wt\zz_l = {\zz} \otimes
u^{-1}\C[[u^{-1}]]$, $\widehat\zz_l = {\zz} \otimes
u^{-1}\C[u^{-1}]$.

Therefore we have the following decomposition of the universal
enveloping algebras: \bean
\label{decom of U} U\agll^- = U\asll^- \otimes U\widehat\zz_l^-. \eean
Here $U\widehat\zz_l^-=\C[z_i]_{i\in\Z_{>0}}$, where
$z_i\in\widehat\zz_l^-$ are scalar matrices $u^{-i}\bs 1$, is the
universal enveloping algebra of $\widehat\zz_l$.
In particular, $U\widehat\zz_l^-$ is the center of $U\agll^-$. 

We define integral forms of $U\asll^-$ and $U\widehat\zz_l^-$. The Lie
algebra $\asll^-$ is generated by $f_i = E_{i+1,i} \otimes 1$
($i=1,\dots,l-1$) and $f_0=E_{1,l} \otimes u^{-1}$. Let $U_\Z \asll^-$
be the $\Z$--subalgebra of $U\asll^-$ generated by $f_i^{(j)}$
($i=0,\ldots,l; j \in \Z_{>0}$). Let $U_\Z \widehat\zz_l^-$ be the
$\Z$--subalgebra of $U\widehat\zz_l^-$ generated by $p_n, n\geq 0$,
defined by formula \eqref{pn}. Then we set
$$
U_\Z \agll^- = U_\Z \asll^- \otimes U_\Z \widehat\zz_l^-.
$$

According to \propref{tensor product}, the algebra $U_\Z\agll^-$
becomes isomorphic $\mc H_l$ after localization at $l$. More
precisely, we have an isomorphism
$$
U_\Z\agll^- \underset{\Z}\otimes \Zl \overset{\sim}\longrightarrow
{\mc H}_l \underset{\Z}\otimes \Zl
$$
sending $f_i^{(j)}$ to $f_i^{(j)}$ ($i=0,\dots,l-1; j \in \Z_{>0}$)
and $p_n$ to $p_n$ ($n\in\Z_{>0}$).

\subsection{Decomposition of $\Z[\wGLl]$}    \label{decomp of gl}

We have the exponential maps $\agll^- \to \wGLl(\C)$, $\asll^- \to
\wSLl(\C)$ and $\widehat\zz_l^- \to \wt{\mc Z}_l^-(\C)$, which are all
isomorphisms. In particular, complex matrices of the form $\exp(A)$
with $A \in \wt\sw_l^-$ are precisely the matrices of determinant
$1$. Now the decomposition \eqref{dir sum of la} gives rise to a
direct product decomposition of groups
$$
\wGLl(\C) = \wSLl(\C) \cdot \wt{\mc Z}_l^-(\C).
$$
Hence we obtain the following decomposition of Hopf algebras over $\C$:
\begin{equation}    \label{decom of CMl}
\C[\wGLl] =\C[\wSLl] \otimes \C[\wt{\mc Z}_l^-].
\end{equation}

There is a similar decomposition over the ring $\Zl=\Z[1/l]$. Indeed,
any element $f(u)$ of $\Zl[[u^{-1}]]$ such that $f(0) = 1$ has a
unique $l$th root $f(u)^{1/l}$ of the same form. Therefore over $\Zl$
any matrix $M$ in $\wGLl$ may be written uniquely as the product
$$
M = M_1 \cdot M_2 = (M (\on{det} M)^{-1/l}) \cdot (\on{det} M)^{1/l}
{\bs 1}
$$ of two matrices, $M_1 \in \wSLl$ and $M_2$, which is a diagonal
matrix. The entries of $M_1$ and $M_2$ are polynomials in the entries
of $M$ with coefficients in $\Zl$, and vice versa. Therefore we
obtain a tensor product decomposition
\begin{equation}    \label{dec of int forms}
\Zl[\wGLl] = \Zl[\wSLl] \underset{\Zl}\otimes
\Zl[\wt{\mc Z}_l^-].
\end{equation}
Under the pairing
$$
({\mc H}_l \underset{\Z}\otimes \Zl) \otimes \Zl[\wGLl] \to
\Zl
$$
induced by the pairing of \propref{def of pairing for Hl}, the
decomposition \eqref{dec of int forms} is dual to the decomposition
\eqref{decom of Hall}. In particular, the orthogonal complement of
$(U_\Z \asll^-)_+$ in $\Zl[\wGLl]$ is equal to $\Zl[\wt{\mc
Z}_l^-]$, and the orthogonal complement of $(U_\Z \widehat\zz_l^-)_+$
in $\Zl[\wGLl]$ is equal to $\Zl[\wSLl]$. The corresponding
pairing $$(U_\Z \asll^- \underset{\Z}\otimes \Zl) \otimes
\Zl[\wSLl] \to \Zl$$ may be described in the same way as in
\lemref{two dual}.

\section{The inductive limit of $\Rep_a
\ur\agln$}\label{ind str root}

In this section we define the inductive limit $\Rep \ur \agli =
\underset{\longrightarrow}{\lim}\;\Rep\,\ur\agln$ and a Hopf algebra
structure on it. We identify this Hopf algebra with the Hopf algebra
$\Z[\wGLl]$ and its Hopf subalgebra $\Rep_{a^*}\urs\agli$ with
$\Z[\wt{Z}^-_l]$. Furthermore, we will identify the quotient of $\Rep
\ur \agli$ by the augmentation ideal of $\Rep_{a^*}\urs\agli$ with
$\Rep_a \fu \agli$ and with $\Z[\wt{SL}^-_l]$. Then we describe the
action of the Hall algebra ${\mc H}_l$ on $\Rep \ur \agli$ in
representation theoretic terms. Finally, we identify the subspace of
$\Rep_a \ur\agli$ (resp., $\Rep_a \fu \agli$) spanned by the
evaluation representations taken at a fixed set of points of
$a\ep^{2\Z}$ with an integrable representation of ${\mc H}_l$ (resp.,
$\wh{\sw}^-_l$).

\subsection{A Hopf algebra structure on $\Rep_a \ur\agli$}\label{ind
limit def root section}

The definition of the Hopf algebras $\Rep\, \ur\agli$ and
$\Rep_a \ur\agli$ is completely parallel to that of $\Rep\,
\uqr\agli$,  $\Rep_a\uqr\agli$; see Sections 
\ref{def of ind section}, \ref{stab section}
and \ref{hopf str section}.

Let ${\mc T}_N^l=\Z[t_{i,n}]_{i=1,\dots,N}^{n=1,\dots,l}$, ${\mathcal
T}_\infty^l=\Z[t_{i,n}]_{i\in\Z_{>0}}^{n=1,\dots,l}$ be commutative rings
of polynomials. We introduce a $\Z_{\geq 0}$--gradation on them by
setting the degree of $t_{i,n}$ to be equal to $i$.  Fix $a\in
\C^{\times}$. By \Ref{iso repa root}, for $N\in\Z_{>0}$, we have
isomorphisms of rings 
\bea
\Rep_a\ur\agln \overset{\sim}\to {\mathcal T}_N^l,\qquad
{}[V_{\om_i}(a\ep^{2n})]\mapsto t_{i,n}. \eea We call the ring
${\mc T}^l_\infty \simeq
\underset{\longrightarrow}{\lim}\;\Rep_a\ur\agln$ the Grothendieck
ring of $\ur\agli$ (restricted to the lattice with base $a$). We also
denote it by $\Rep_a\ur\agli$.

The monomial basis of ${\mathcal T}^l_\infty \simeq
\Rep_a\ur\agli$ corresponds to the basis of tensor products of
the fundamental representations.  We call this basis {\em standard}.

The specialization maps $\wh{\mc S}_N^l$ obviously stabilize, so we have a
map
\bea
\wh{\mc S}^l:\;\Rep_a U_q\agli\to\Rep_a\ur\agli,\qquad
t_{i,j}\mapsto t_{i,j \on{mod} l}.
\eea
By Lemma \ref{spec fund lem}, the specialization map sends 
the PBW basis of $\Rep_a U_q\agli$ to the standard basis of
$\Rep_a\ur\agli$.

Similarly, we set
\be
\Rep\,\ur\gln = \underset{\longrightarrow}{\lim} \; \Rep\,\ur\gln
\simeq \Z[t_i]_{i\in\Z_{>0}},
\ee
where $t_i$ denotes the class of the $i$th fundamental representation
of $\ur\gln$ and the inductive limit of the maps ${\mc S}_N^l$ defines
an isomorphism of Hopf algebras ${\mc S}^l:\;
\Rep\,U_q\gli\to\Rep\,\ur\gli$.

Recall the definition of minimal dominant $M$--tuple from \secref{fin
dim root}. Similarly to Lemma \ref{fund to irrep stable} we obtain the
following 

\begin{lem}\label{can lem}
Let ${\bs Q}$ be a minimal dominant $M$--tuple. Then the sequence of
irreducible representations of $\ur\agln$ corresponding to the
dominant $N$--tuple ${\bs Q}^{(N)}$ ($N\geq M$) stabilizes and
therefore defines an element of $\Rep_a\ur\agli$. The set of all such
elements (for all $M\in\Z_{\geq 0}$) is a basis in $\Rep_a\ur\agli$.
\end{lem}

We call the basis described in Lemma \ref{can lem} {\em canonical}.
Note that the transition matrix from the standard basis to the
canonical basis is triangular with non-negative integer entries and
$1$'s on the diagonal.

We also have another natural basis in $\Rep_a\ur\agli$, given
by specializations $V({\bs Q})^\ep$ of irreducible modules for generic
$q$. The elements of this basis are also labeled by minimal dominant
$M$--tuples, $M\in\Z_{>0}$. By formula \eqref{tri matr}, the
transition matrix from this basis to the canonical basis is triangular
with $1$'s on the diagonal.

Let $N=N_1+N_2$. Then we have the embedding
$\ur{\widehat{\gw}_{N_1}}\otimes \ur{\widehat{\gw}_{N_2}}\to
\ur\agln$. Using this embedding we obtain the comultiplication
structure on $\Rep_a\ur\agli$ as in \secref{hopf str
section}. It is given on the generators by \be \Delta
(t_{i,n})=\sum_{j=0}^i t_{j,n}\otimes t_{i-j,n-j}, \ee where we set
$t_{i,n}=t_{i,n+l}$.

Then we have the following theorem, which identifies $\Rep_a\ur\agli$
together with various structures on it with the Hopf algebra
$\Z[\wGLl]$ (the restricted dual of the Hopf algebra ${\mc H}_l$)
introduced in \secref{unwind}.

\begin{thm}\label{repa to uminus root}
{\em (1)} There exists an antipode $S$, such that
$(\Rep_a\ur\agli,\otimes,\Delta,S)$ is a Hopf
algebra. Furthermore, $\Rep_a\ur\agln$ is a Hopf subalgebra of
$\Rep_a\ur\agli$.

{\em (2)} We have a commutative diagram of Hopf algebras
$$
\begin{CD}
\Z[\wGLl] @>{\sim}>> \Rep_a \ur\agli \\
@AAA @AAA \\ 
\Z[\wGLln] @>{\sim}>> \Rep_a \ur\agln \\
\end{CD}
$$
where the horizontal maps are given by the formula $M_{i+j,j}\mapsto
t_{i,j}$.

{\em (3)} We have a commutative diagram of Hopf algebras
$$
\begin{CD}
\Z[\wGLl] @>{\sim}>> \Rep_a \ur\agli \\
@A{w^*}AA @A{\wh{\mc S}^l}AA \\ 
\Z[\wSLi] @>{\sim}>> \Rep_a \uqr\agli\\
\end{CD}
$$
Furthermore, the standard basis of $\Rep_a\ur\agli$ is
identified with the basis dual to the basis of ${\mc A}$--sets in $\mc
H_l$.
\end{thm}

\subsection{Restriction operators}\label{restrictions section root}

This section is analogous to Sections \ref{def restr section},
\ref{restrictions section}. Since $\Rep_a\ur\agli$ is
(restricted) dual to ${\mc H}_l$, the algebra ${\mc H}_l$ acts on
$\Rep_a\ur\agli$. We will show that the action of the
generators $f_{\bs d}$ of ${\mc H}_l$ is given by the restriction
operators $\res_{\bs d}$ that we now define.

The embedding of algebras
$\ur\widehat{\mathfrak{gl}}_{N-1}\otimes\ur\widehat{\mathfrak{gl}}_{1}\to
\ur\agln$ gives us a linear map \be \Rep_a \ur\agln \to \Rep_a \ur
\wh{\gw}_{N-1} \otimes \Rep_a \ur\widehat{\gw}_1.  \ee The ring
$\ur\widehat{\gw}_1=\Z[X_n]_{n=0,\dots,{l-1}}$, where we set
$t_{1,n}=X_n$, has a monomial basis $\{X_{\bs d} =
X_0^{d_0}X_1^{d_1}\dots X_{l-1}^{d_{l-1}}\}_{d_i\in\Z_{\geq 0}}$
labeled by vectors ${\bs d} \in\Z^{l}_{\geq 0}$. Any representation
$V$ of $\ur\agln$ may be decomposed as
$$
[V] = \sum_{\bs d} [V_{\bs d}] \otimes X_{\bs d}.
$$
We then define the restriction operator $\res^{\ep,N}_{\bs d}$ by the
formula (cf. formula \eqref{resm}):
\bea
\res^{\ep,N}_{\bs d} :\;\Rep_a \ur\widehat{\gw}_N 
\to \Rep_a \ur\widehat{\gw}_{N-1},\qquad {[V]} \mapsto [V_{\bs
d}].
\eea

Introduce the notation ${\mc L}^\ep_N =
\Z[\Lambda_{i,a\ep^{2n}}]_{i=1,\ldots,N}^{n=1,\dots,l}$ and consider
the natural homomorphism $\al^\ep_N: {\mc L}^\ep_N \to {\mc
L}^\ep_{N-1} \otimes {\mc L}^\ep_1$. Given $P \in {\mc L}^\ep_N$, we
write \be \al^\ep_N(P) = \sum_{\bs d} P_{\bs
d} \otimes (\Lambda_{N,a}^{d_0} \La_{N,a\ep^2}^{d_1} \ldots
\La_{N,a\ep^{2l-2}}^{d_{l-1}}). \ee
Define a linear operator $\wt{\res}^{\ep,N}_{\bs d}: {\mc L}_N
\to {\mc L}_{N-1}$ sending $P$ to $P_{\bs d}$. 
Then we have an analogue of formula \eqref{restr and lambda}:
\begin{equation}    \label{restr and lambda root}
\chi_\ep(\res^{\ep,N}_{\bs d} [V]) = \wt{\res}^{\ep,N}_{\bs d}
(\chi_\ep(V)).
\end{equation}
We remark that to ${\bs m} = (m_1 \geq \ldots \geq m_k)$, $0 \leq m_i
\leq l-1$, of formula \eqref{restr and lambda} corresponds ${\bs d} =
(d_0,\ldots,d_{l-1})$, where $d_i$ is the number of times $i$ occurs
in ${\bs m}$.

The maps $\res^{\ep,N}_{\bs d}$ stabilize in the limit $N \to \infty$
and give rise to the operators \be \res^\ep_{\bs
d}:\;\Rep\,\ur\agli\to\Rep\,\ur\agli \qquad ({\bs
d}\in\Z^{l}_{\geq 0}).  \ee

Recall that $f_{\bs d}\in\mc H_l$ is the element corresponding to the
$A^{(1)}_l$--set of degree ${\bs d}$ which has snakes only of length
$0$. By \lemref{all gen}, $\mc H_l$ is generated by these
elements. Using Lemma \ref{f is dual to res} and \thmref{repa to
uminus root} we obtain that the operator $\res_{\bs d}$ is dual to the
operator of right multiplication by $f_{\bs d}$. Hence we obtain the
following analogue of \thmref{action of sl}:

\begin{thm}\label{hall action}
There exists an action of Hall algebra ${\mc H}_l$ on
$\Rep_a\ur\agli$, such that $f_{\bs d}\in{\mc H}_l$ acts as
$\res^\ep_{\bs d}$ ($\bs d\in\Z^l$). In particular, there exists an
action of $U_\Z\sw_l^-$ on $\Rep_a\ur\agli$ such that $f^{(j)}_{i}\in
U_\Z\sw_l^- $ acts as $\res^\ep_{i,\ldots,i}$, where the index $i$
appears $j$ times ($i=0,\dots,l-1$).  Under this action,
$\Rep_a\ur\agli$ is isomorphic to the restricted dual of a free $\mc
H_l$--module of rank one.
\end{thm}

\subsection{The Hopf algebras $\Rep_{a^*}\urs\agli$ and $\Rep_a
  \fu\agli$}\label{tensor product section}

By Lemma \ref{char of frob}, the Frobenius maps
$\on{Fr}^*:\;\Rep\,\urs\agln\to\Rep\,\ur\agln$ stabilize in the limit $N
\to \infty$ and therefore give rise to an embedding
\bean    \label{Frinfty}
\on{Fr}^*:\;\Rep_{a^*}\urs\agli&\to&\Rep_a\ur\agli.
\eean
Lemmas \ref{Fr is homomorphism} and \ref{char of frob} imply:

\begin{lem}    \label{Fr Hopf}
The map \eqref{Frinfty} is a homomorphism of Hopf algebras.
\end{lem}

Given a subalgebra $B$ of an algebra $A$, we denote by $(B)_+$ the ideal
of $A$ generated by the augmentation ideal of $B$.

\lemref{Fr Hopf} implies that the inductive limit of the rings $\Rep\,
\fu\gln$, $\Rep\fu\agln$ with the factorized tensor product structure is
well-defined. We set \bea \Rep\,
\fu\gli&:=&\underset{\longrightarrow}{\lim}\; \Rep\, \fu\gln=\Rep\,
\ur\gli/(\Rep\, \urs\gli)_+,\\ \Rep\,
\fu\agli&:=&\underset{\longrightarrow}{\lim}\; \Rep\, \fu\agln=\Rep\,
\ur\agli/(\Rep\, \urs\agli)_+,\\ \Rep_a
\fu\agli&:=&\underset{\longrightarrow}{\lim}\;\Rep_a\fu\agln = \Rep_a
\ur\agli/(\Rep_{a^*} \urs\agli)_+. \eea Moreover, the rings
$\Rep\fu\gln$, $\Rep\, \fu\agli$ and $\Rep_a \fu\agli$ become Hopf
algebras.

Note that the algebra structures on $\Rep\, \fu\agln$, $\Rep\, \fu\gln$ 
coming from the ordinary tensor
product are not stable in the limit $N \to \infty$. That is why we use
the factorized tensor product structure.

By \thmref{repa to uminus root}, we have an isomorphism of Hopf
algebras $\Rep_a\ur\agli \simeq \Z[\wGLl]$. The Hopf algebra
$\Rep_a\ur\agli$ has the Hopf subalgebra $\Rep_{a^*}\urs\agli =
\Z[{\mb V}_{\omega_i}]_{i>0}$, where ${\mb V}_{\omega_i} =
\on{Fr}^*(V_{\omega_i}(a^*))$, and the Hopf algebra $\Z[\wGLl]$
has the Hopf subalgebra $\Z[\wt{\mc Z}_l] = \Z[a_i]_{i>0}$. We
now show that these two subalgebras are identified under this
isomorphism.

\begin{thm}\label{fr ort sl}
Under the isomorphism $\Rep_a\ur\agli \simeq \Z[\wGLl]$ we have an
isomorphism of the Hopf subalgebras $\Rep_{a^*}\urs\agli$ and
$\Z[\wt{\mc Z}_l]$ sending the generators ${\mb V}_{\omega_i}, i>0$,
to the generators $(-1)^{l+i} a_i, i>0$.
\end{thm}

\begin{proof}
Recall from \secref{decomp of gl} that $\Z[\wt{\mc Z}_l]$ is equal to
the orthogonal complement of the augmentation ideal of $U_\Z \asll^-
\subset {\mc H}_l$ in $\Z[\wGLl]$. This implies that $\C[\wt{\mc
Z}_l]$ is equal to the subspace of $\asll^-$--invariants of
$\C[\wGLl]$, or equivalently, the intersection of kernels of the
operators $f_i$, ($i=1,\ldots,l$). Let us show that each ${\mb
  V}_{\omega_i}$ is also annihilated by the $f_i$'s.

By \thmref{hall action}, the operator $f_i$ acts as the restriction
operator $\res_i$. Lemma \ref{char of frob} and formula \eqref{restr
  and lambda root} immediately imply that 
\be
\res_i(\on{Fr}^*(V))=0
\ee
for any $\urs\agln$ module $V$. Therefore we find that
$\Rep_{a^*}\urs\agli \otimes \C$ is contained inside $\C[\wt{\mc
Z}_l]$. Comparing the graded characters of the two spaces, we find
that $$\Rep_{a^*}\urs\agli \otimes \C = \C[\wt{\mc Z}_l].$$

Both algebras are $\Z_{\geq 0}$--graded Hopf algebras, freely
generated by the generators ${\mb V}_{\omega_i}, i>0$, and $a_i, i>0$,
respectively whose degrees are given by the formulas $\deg {\mb
  V}_{\omega_i} = \deg a_i = i$. Therefore ${\mb V}_{\omega_i}$ is
equal to a non-zero multiple of $a_i$ plus a polynomial in $a_j, j<i$,
of degree $i$. The comultiplication for the $a_i$'s is given by
formula \eqref{delta of ai}, and we find the comultiplication formula
for the ${\mb V}_{\omega_i}$'s from formula \eqref{deltati}
specialized to the case when $q^2=1$ and the
fact that the the map $\on{Fr}^*$ is a homomorphism of Hopf
algebras. The result is
$$
\Delta({\mb V}_{\omega_i}) = \sum_{j=0}^l {\mb V}_{\omega_j} \otimes
      {\mb V}_{\omega_{i-j}}
$$
(where we set ${\mb V}_{\omega_0}=1$), which coincides with the
formula for $\Delta(a_i)$. This readily implies that ${\mb
V}_{\omega_i}$ is equal to a non-zero multiple of $a_i$.

It remains to show that this multiple is equal to $(-1)^{l+i}$. In
order to do that we find the coefficient of one of the monomials in
both ${\mb V}_{\omega_i}$ and in $a_i$. We take this monomial to be
$\ds \prod_{j=0}^{l-1} t_{i,j}$. It enters with coefficient $1$ in the
decomposition of ${\mb V}_{\omega_i}$ into a linear combination of the
products of the fundamental representations $V_{\omega_k}(a\ep^{2m}) =
t_{k,m}$, because the tensor product $\bigotimes_{i=0}^{l-1}
V_{\omega_i}(a\ep^{2j})$ has the same highest weight as ${\mb
  V}_{\omega_i}$.

On the other hand, recall that $t_{k,m}$ corresponds to
$M_{k+m,m}$. In the loop group realization of $\wGLl$ introduced in
\secref{new sect} each element of $\wGLl$ is represented as an $l
\times l$ matrix with entries being formal power series in
$u^{-1}$. In terms of the $t_{k,m}$'s, the $(i,j)$th entry of this
matrix is equal to
\begin{equation}    \label{ij entry}
\sum_s t_{ls+i-j,j} u^{-s},
\end{equation}
where $s\geq 0$, if $i\geq j$ and $s>0$, if $i < j$.

The generator $a_i$ is the $u^{-i}$--coefficient of the determinant of
this matrix. Therefore it contains the monomial $\ds (-1)^{i+l}
\prod_{j=0}^{l-1} t_{i,j}$. This completes the proof.
\end{proof}

Note that in the course of the proof we obtained an explicit formula
expressing the Frobenius pullback ${\mb V}_{\omega_i} =
\on{Fr}^*(V_{\omega_i}(a^*))$ as a linear combination of the products
of the fundamental representations. For example, if $l=2$, then
\be
\on{Fr}^*(V_{\om_{i}}(a^*)) =
t_{i,0}t_{i,1}+\sum_{j=1}^i(-1)^j(t_{i+j,0}t_{i-j,1}+t_{i+j,1}t_{i-j,0}).
\ee

Now recall that have by definition (see \secref{tensor product
  section} and \secref{new sect})
\begin{align*}
\Rep_a\fu\agln &= \Rep_a \ur\agli/(\Rep_{a^*} \urs\agli)_+, \\
\Z[\wSLl] &= \Z[\wGLl]/(\Z[\wt{\mc Z}_l])_+.
\end{align*}
Therefore \thmref{fr ort sl} implies

\begin{thm}    \label{sl ort fr}
The Hopf algebras $\Rep_a\fu\agln$ and $\Z[\wSLl]$ are isomorphic.
\end{thm}

Thus, we obtain the following commutative diagram of Hopf algebras and
homomorphisms between them:
$$
\begin{CD}
\Z[\wt{\mc Z}_l] @>{\sim}>> \Rep_{a^*}\urs\agli \\
@VVV @VVV \\
\Z[\wGLl] @>{\sim}>> \Rep_a \ur \agli \\
@VVV @VVV \\
\Z[\wSLl] @>{\sim}>> \Rep_a \fu \agli .\\
\end{CD}
$$

\subsection{Integrable representations of $\wt\sw_l$}\label{swl}

We recall some facts about the Lie algebra $\asll$ (see \cite{K} for
details). It has generators $e_i, h_i, f_i$ ($i=0,\dots,l-1$) satisfying
the standard relations. In particular, the Lie algebra $\sw_l^-$
introduced in \secref{on decomp} is identified with the Lie subalgebra
of $\asll$ generated by $f_i$ ($i=0,\dots,l-1$). The weight lattice of
$\asll$ is spanned by the fundamental weights $\varpi_i$ ($i=0,\dots,
l-1$). It contains the root lattice spanned by $\al_i = 2 \varpi_i -
\varpi_{i-1} - \varpi_{i+1}$ ($i=0,\dots,l-1$), where we set
$\varpi_{-1}=\varpi_{l-1}$, $\varpi_{l}=\varpi_{0}$.
For weights ${\bs\la},{\bs\mu}$ we
write ${\bs\la} \geq {\bs\mu}$ if ${\bs\la} = {\bs\mu} + \sum_i n_i
\al_i$, where all $n_i 
\in \Z_{\geq 0}$.

Introduce the category ${\mc O}$ of $\asll$--modules. Its objects are
$\asll$--modules $L$ on which the operators $h_i$ ($i=0,\dots,l-1$) act
diagonally, and the weights occurring in $L$ are less than or equal to
some fixed weight ${\bs\la}$. The morphisms in ${\mc O}$ are
$\asll$--homomorphisms. A representation $L$ of $\asll$ is called
integrable if it belongs to ${\mc O}$ and all generators $e_i$ and
$f_i$ act on $L$ locally nilpotently. It is called a highest weight
representation if it is generated by a vector $v$, such that $e_i
\cdot v = 0$ ($i=0,\dots,l-1$). The category of integrable
representations of $\asll$ is semi-simple. Irreducible integrable
representations of $\asll$ are the irreducible highest weight
representations with the highest weights $\bs\nu = \sum_i \nu_i
\varpi_i$, where $\nu_i \in \Z_{\geq 0}$ ($i=0,\dots,l-1$).
Such weights are called dominant integral weights.

The irreducible representation corresponding to a dominant integral
weight $\bs\nu = \sum_i \nu_i \varpi_i$ is generated by a vector
$v_{\bs\nu}$, such that
$$      
e_i v_{\bs\nu} = 0, \qquad h_i v_{\bs\nu} = \nu_i v_{\bs\nu}, \qquad
f_i^{\nu_i+1} v_{\bs\nu} = 0 \qquad (i=0,\dots,l-1).
$$
It is denoted by $L_l(\bs\nu)$. The number $k =
\sum_i \nu_i$ is called the level of $L_l(\bs\nu)$. Since $\sw_l^-$ is
a Lie subalgebra of $\asll$, its universal enveloping algebra
$U\asll^-$ acts on any $\asll$--module. In particular, we
have a canonical surjective homomorphism of $U\asll^-$--modules
\bean
\label{from U to L root} p_{\bs\nu}: U\asll^- \to
L_l(\bs\nu), \qquad g \to g \cdot v_{\bs\nu}. \eean

Fix $n\in\{0,\dots,l-1\}$. Recall from \secref{int section} that
$\bigwedge$ is the vector space with the basis of vectors $|{\bs \la}
\rangle$ labeled by all partitions ${\bs \la}$. Corresponding to $n$,
we have operators $\wt{e}_m$, $\wt{f}_m$ and $\wt{h}_m$ ($m\in\Z$) on
$\bigwedge$ introduced in Section \ref{int section}. The next lemma
follows from \lemref{Hayashi} and from complete reducibility of
integrable $\asll$--modules.

\begin{lem}\label{L and Lambda}
The operators $e_i=\sum_{m\in\Z}\wt{e}_{ml+i}$,
$f_i=\sum_{m\in\Z}\wt{f}_{ml+i}$, $h_i=\sum_{m\in\Z}\wt{h}_{ml+i}$
($i=0,\dots,l-1$) define a 
representation of $\asll$ on $\bigwedge$ denoted by
${\bigwedge}_l(\varpi_n)$. 
Moreover, the irreducible module $L_l(\varpi_n)$ is the direct summand
of ${\bigwedge}_l(\varpi_n)$ equal to its cyclic $\asll$ submodule generated
by the vector $|0\rangle$. 
\end{lem}

\subsection{Evaluation modules and integrable modules}\label{int mod
section root}
Recall the evaluation homomorphism, $ev_{a}: \ur \agln \to \ur
\gln$. The next lemma follows directly from definitions.

\begin{lem}\label{spec and ev}
The following diagram is commutative:
\be
\begin{CD} 
\Rep\, U_q\gli @>{ev_a}>> \Rep_a U_q\agli \\
@V{\mc S_l}VV @VV{\widehat{\mc S}_l}V\\  
\Rep\, \ur\gli @>ev_a>>\Rep_a \ur\agli.
\end{CD}
\ee
\end{lem}
In particular, the $\ep$--characters of evaluation modules are obtained
by specialization of $q$--characters described by Lemma
\ref{evaluation character lemma}.

Fix $n\in\{0,\dots,l-1\}$.
Then $ev_{a\ep^{2n}}$ induces an embedding
\bean\label{gl in agl root} \iota^{(N)}_l(\varpi_n):\;\Rep\,
\ur\gln \to \Rep_a \ur\agln,\qquad [V^{\ep}_{\bs\la}] \mapsto
[V^{\ep}_{\bs\la}(a\ep^{2n})], \eean where $V^\ep_{\bs\la}$
denotes the specialization of the irreducible representation
$V_{\bs\la}$ of $U_q\gln$ of highest weight $\bs\la$. The embeddings
\eqref{gl in agl root} stabilize and hence induce an embedding
\bean\label{gli in agli root}
\iota_l(\varpi_n):\;\Rep\,\ur\gli \to \Rep_a
\ur\agli. \eean

Recall from \thmref{repa to uminus root} and \propref{def of pairing
  for Hl} that we have a pairing of Hopf algebras
$$
\langle,\rangle: \mc H_l \otimes \Rep_a\ur\agli \to \Z.
$$
In particular, $\Rep_a\ur\agli$ is an $\mc H_l$--module.
The following lemma is obtained in the same way as Lemma \ref{iotan}.

\begin{lem}    \label{iotan root}
The $\Z$--submodule $\iota_l(\varpi_n)(\Rep\, \ur\gli)\subset\Rep_a
\ur\agli$ is a right comodule of $\Rep_a \ur\agli$ and an $\mc
H_l$--submodule of $\Rep_a\ur\agli$.
\end{lem}

We also have a pairing \bean\label{pairing of level 1 root}
\langle,\rangle_n:\;\bigwedge{}_l(\varpi_n)\otimes
(\iota_l(\varpi_n)(\Rep\,\ur\gli) \otimes \C) \to \C, \eean such that
the bases $\{|\bs\la\rangle\}$ and $\{[V^\ep_{\bs\la}(a\ep^{2n})]\}$
are dual.

The algebra ${\mc H_\infty}$ acts on $\bigwedge_l(\varpi_n)$, and this
action extends to its completion $\wt{\mc H}_\infty$.  Therefore, $\mc
H_l$ acts on $\bigwedge_l(\varpi_n)$ via the unwinding map $w:\;\mc
H_l\to \wt{\mc H}_\infty$. Moreover, the restriction of this action to
$U_\Z\wh{\mathfrak{sl}_l^-}$, coincides with the action of
$U_\Z\wh{\mathfrak{sl}}_l^-\subset U \wh{\mathfrak{sl}}_l$ described
in Section \ref{swl}.

Combining Proposition \ref{level one} and Lemma \ref{spec and ev} we
obtain the following

\begin{prop}    \label{level one root}
The $\mc H_l$--module $\iota_l(\varpi_n)(\Rep\, \ur\gli) \otimes \C$
is restricted dual to $\bigwedge_l(\varpi_n)$ with respect to the
pairing \eqref{pairing of level 1 root}. Furthermore, the embedding
$$\iota_l(\varpi_n):\;\Rep\, \ur\gli \otimes \C \to \Rep_a \ur\agli
\otimes \C$$ is dual to the surjective $\mc H_l \otimes
\C$--homomorphism $$\wt p_{\varpi_n}: \mc H_l \otimes \C \to
\bigwedge{}_l(\varpi_n)$$ sending $1$ to $v_{\varpi_n}$.
\end{prop}

\propref{level one root} implies that as a $\Rep_a
U_q\agli$--comodule, $\iota_l(\varpi_n)(\Rep\, \ur\gli)$ is an
integral form of the ${\mc H}_l \otimes \C$--module which is the
restricted dual to $\bigwedge_l(\varpi_n)$. Moreover, this integral
form is preserved by ${\mc H}_l$.

By Theorem \ref{sl ort fr}, the pairing $\langle,\rangle$ gives rise
to the pairing \bean\label{pairing of level 1 root factor}
\langle,\rangle:\; U_\Z\wh{\mathfrak{sl}}_l^-\otimes \Rep_a\fu\agli \to
\Z.  \eean Moreover, since the Frobenius and
evaluation maps commute, the map $\iota_l(\varpi_n)$ gives rise to the
map \be \iota_l(\varpi_n):\; \Rep\fu\gli\to\Rep_a\fu\agli.  \ee
Alternatively, this maps sends $[V_{\bs\la}^\ep]$ with
$\la_i-\la_{i+1}<l$ ($i\in\Z_{>0}$) to $[V_{\bs\la}^\ep(a\ep^{2n})]$.

Combining Proposition \ref{level one root} and Theorem \ref{sl ort
fr}, we obtain the following

\begin{prop}    \label{level one root factor}
The $U\wh{\mathfrak{sl}}_l^-$--module $\iota_l(\varpi)(\Rep\, \fu\gli)
\otimes \C$ is restricted dual to the irreducible integrable module
$L_l(\varpi_n)$ with respect to the pairing \eqref{pairing of level 1
root factor}. Furthermore, the embedding $$\iota_l(\varpi):\;\Rep\,
\fu\gli \otimes \C \hookrightarrow \Rep_a \fu\agli \otimes \C$$ is dual
to the surjective $U\wh{\mathfrak{sl}}_l^-$--homomorphism
$$p_{\varpi_n}:U\wh{\mathfrak{sl}}_l^-\to L_l(\varpi_n).$$
\end{prop}

Now we describe the integrable modules of higher levels.
Given a dominant integral weight $\bs\nu = \sum_i \nu_i \varpi_i$, set
$\bigwedge_l(\bs\nu) = \otimes_i (\bigwedge_l({\varpi_i})^{\otimes
\nu_i})$. This 
is an integrable $\wh{\sw}_l$--module of level $k= \sum_i
\nu_i$. Denote the tensor product of the highest weight vectors of
$\bigwedge({\varpi_i})$ by $|0\rangle_{\bs\nu}$.

Similar to Lemma \ref{fock as cyclic submodule}, we obtain the
following

\begin{lem}\label{fock as cyclic submodule root}
The irreducible module $L_l(\bs\nu)$ is a direct summand in
$\bigwedge_l(\bs\nu)$, which is equal to the cyclic $U\wh{\mathfrak
s\mathfrak l}_l^-$--submodule of $\bigwedge_l(\bs\nu)$ generated by
$|0\rangle_{\bs\nu}$.
\end{lem}

Denote by $\wt{p}_{\bs\nu}$ the homomorphism $U\mathfrak s\mathfrak
l_\infty^- \to \bigwedge(\bs\nu)$ obtained by composing $p_{\bs\nu}$ given by
formula \eqref{from U to L root} and embedding
$L_l(\bs\nu)\to\bigwedge_l(\bs\nu)$.

Write $\bs\nu$ in the form $\sum_{j=1}^k \varpi_{s_j}$. Using
\propref{level one root}, we identify the restricted dual
$\bigwedge_l^*(\bs\nu)$ with $(\Rep\, \ur\gli)^{\otimes k}$ as $\mc
H_l$--modules. Denote by $\langle,\rangle_{\bs\nu}$ the pairing
$\bigwedge_l(\bs\nu) \otimes (\Rep\, \ur\gli)^{\otimes k} \to \C$.

The map \bean\label{many gl to agl root} \iota_l(\bs\nu):\;
(\Rep\,\ur\gli)^{\otimes k}&\to& \Rep_a \ur\agli,\\
{[V^\ep_{\bs\la_{1}}]\otimes\dots\otimes [V^\ep_{\bs\la_{k}}]}
&\mapsto& [V^\ep_{\bs\la_{1}}(a\ep^{2s_1})\otimes\dots\otimes
V^\ep_{\bs\la_{k}}(a\ep^{2s_k})]\notag \eean is a homomorphism of
$\mc H_l$--modules.

We have a diagram of $\mc H_l \otimes \C$--modules and $\mc H_l
\otimes \C$--homomorphisms
$$\begin{CD} (\Rep\, \ur\gli)^{\otimes k} \otimes \C
@>{\iota_l(\bs\nu)}>> \Rep_a \ur\agli \otimes \C \\
@V{\langle,\rangle_{\bs\nu}}VV @VV{\langle,\rangle}V\\
\bigwedge_l(\bs\nu)^* @>\wt{p}_{\bs\nu}^*>> (\mc H_l \otimes \C)^*
\end{CD}$$

Similar to Proposition \ref{many gl to agl lemma}, we obtain the following

\begin{prop}\label{many gl to agl root lemma}
The above diagram is commutative. The image of the map
$\iota_l(\bs\nu)$ is an ${\mc H}_l \otimes \C$--submodule of $\Rep_a
\ur\agli \otimes \C$ isomorphic to the restricted dual of
$\bigwedge_l(\bs\nu)$.
\end{prop}

Finally, we factorize the diagram above to obtain
a diagram of $U\wh{\mathfrak{sl}}_l^-$--modules and
$U\wh{\mathfrak{sl}}_l^-$--homomorphisms
$$\begin{CD} (\Rep\, \ur\gli)^{\otimes k} \otimes \C
@>{\wt\iota_l({\bs\nu})}>> \Rep_a \fu\agli \otimes \C \\
@V{\langle,\rangle_{\bs\nu}}VV @VV{\langle,\rangle}V\\
\bigwedge(\bs\nu)^* @>\wt{p}_{\bs\nu}^*>> (U\wh{\mathfrak{sl}}_l^-)^*
\end{CD}$$
Similar to Proposition \ref{many gl to agl lemma}, we obtain the following

\begin{prop}\label{many gl to agl root factor lemma}
The above diagram is commutative. The image of the map
$\wt\iota_l(\bs\nu)$ is a $U\wh\sw_l^-$--submodule of
$\Rep_a \ur\agli \otimes \C$ isomorphic to the restricted dual of
$L_l(\bs\nu)$.
\end{prop}

Propositions \ref{many gl to agl root lemma} and \ref{many gl to agl
root factor lemma} imply that the image of $(\Rep\, U_q\gli)^{\otimes
n}$ under $\iota(\bs\nu)$ (resp., $\wt{\iota}(\bs\nu)$) is an integral
form of the restricted dual to the ${\mc H}_l \otimes \C$--module
(resp., $U\wh{\sw}_l^-$--module) $\bigwedge_l(\bs\nu)$ (resp.,
$L_l(\bs\nu)$), which is preserved by ${\mc H}_l$ (resp.,
$U_\Z\wh\sw_l^-$).

We denote the image of the map $\wt\iota_l(\bs\nu)$ by
$\Rep_a^{\bs\nu} \ur\gli$.

\subsection{Combinatorial identity}
A dominant $\bs\nu$--admissible (see
Section \ref{comb section}) $N$--tuple of polynomials ${\bs Q}$ is
called {\em $l$--acyclic} if $Q_i(u)/Q_{i+1}(u\ep^2)$ ($i=1,\dots,
N$) is an $l$--acyclic polynomial (we set $Q_{N+1}(u)=1$).

\begin{conj}\label{can conj root}
The classes of irreducible representations $V({\bs Q})$ which have 
$l$--acyclic $\bs\nu$--admissible 
Drinfeld polynomials ${\bs Q}$, form a basis in
$\Rep_a^{\bs\nu} \ur\gli$.
\end{conj}

Clearly, Conjecture \ref{can conj root} holds for $k=1$. 
 
We end this section by describing a combinatorial corollary of
Proposition \ref{many gl to agl root factor lemma} and Conjecture
\ref{can conj root}.

Fix ${\bs s}=(s_1\geq\dots\geq s_k=0$), such that $s_i\in\Z$. Let
$\bs\la$ be an $\bs s$--diagram (see Section \ref{comb section}). We
replace each box in $\bs\la$ by another box by the rule: a box with
coordinates $(i,j)$ is replaced by the box with coordinates $(i,j \mod
l)$. The result is a union of boxes (with multiplicities) in a
vertical strip of width $l$, which we call an {\em $l$--folded
${\bs s}$--diagram}.

Here is an example of folding for $l=2$.

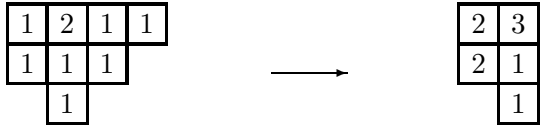
\begin{picture}(0,23)(0,18)
\setlength{\unitlength}{1mm}

\put(55,26){\vector(1,0){10}}
\thicklines
\put(20,30){\framebox(5,5){1}}
\put(25.3,30){\framebox(5,5){2}}
\put(30.6,30){\framebox(5,5){1}}
\put(35.9,30){\framebox(5,5){1}}

\put(20,24.7){\framebox(5,5){1}}
\put(25.3,24.7){\framebox(5,5){1}}
\put(30.6,24.7){\framebox(5,5){1}}

\put(25.3,19.42){\framebox(5,5){1}}

\put(80,30){\framebox(5,5){2}}
\put(85.3,30){\framebox(5,5){3}}

\put(80,24.7){\framebox(5,5){2}}
\put(85.3,24.7){\framebox(5,5){1}}

\put(85.3,19.42){\framebox(5,5){1}}
\end{picture}

Here the number inside each box shows its multiplicity.

An {\em acyclic $l$--folded ${\bs s}$--diagram} is an 
$l$--folded ${\bs s}$--diagram $\bs\la$ with the property:  for any $i\in\Z$,
such that the $i$th row of $\la$ is non-empty, there exists
$j\in\{1,\dots,l\}$ such that multiplicities of boxes with coordinates
$(i,j)$ and $(i+1,j)$ are the same.

In our example the $2$--folded diagram above is acyclic: for the first and
the third rows $j=1$ and for the second row $j=2$.

Let $N_n^{(l)}({\bs s})$ be the number of distinct acyclic
$l$--folded ${\bs s}$--diagrams with
$n$ boxes and let the series 
$\chi_{\bs s}^{(l)}(\xi)=\sum_{i=0}^\infty N_i^{(l)}({\bs s}) \xi^i$
be the corresponding generating function.

By Proposition 
\ref{many gl to agl root factor lemma}, we obtain the following
\begin{lem}
If Conjecture \ref{can conj root} is true then the function $\chi^{(l)}_{\bs
s}(\xi)$ coincides with the 
character of the irreducible integrable $\widehat{\sw}_l$ module of weight
$\sum_{i=1}^k\varpi_{s_i}$ in the principal gradation.
\end{lem}
This character may be found, e.g., in \cite{K}, formula (10.10.1). 
It would be interesting to find a combinatorial proof of this fact.

\section{Summary}\label{summary section}

The general picture described in this paper is summarized in Figure
\ref{figure 2}.

\begin{figure}
\setlength{\unitlength}{1mm}
\begin{picture}(40,116)(50,-57)    \label{pic2}
\put(50,45){\vector(-3,-2){20}} 
\put(63,49){\vector(1,0){50}}
\put(36,27){\vector(1,0){50}}
\put(11,26){$\Rep_a U_q \widehat{\gw}_\infty$}
\put(38,48){$\Rep\, U_q {\gw}_\infty$}
\put(115,48){$\bigwedge(\varpi_0)^*_\Z$}
\put(88,26){$\Z[\wSLi]$}
\put(116,45){\vector(-3,-2){20}} 

\put(85,50){${\sim}$} 
\put(61,28){${\sim}$} 
\put(75,45){${[V_{\bs\la}]\mapsto |\bs\la\rangle^*}$} 
\put(39,34){${ev_a^*}$}

\put(25,23){\vector(0,-1){32}} 
\put(118,13){\line(0,1){32}}
\put(120,13){\line(0,1){32}}
\put(52,24){\vector(0,-1){11}}
\put(52,44){\line(0,-1){15}}
\put(92,23){\vector(0,-1){32}}

\put(20,5){$\wh{\mc S}_l$}
\put(47,19){${\mc S}_l$}
\put(87,-2){$w^*$}
\put(114,27){}

\put(50,5){\vector(-3,-2){20}}
\put(63,9){\line(1,0){26}}
\put(95,9){\vector(1,0){18}}
\put(36,-13){\vector(1,0){50}}
\put(11,-14){$\Rep_a \ur\agli$}
\put(38,8){$\Rep\, \ur\gli$}
\put(115,8){$\bigwedge_l(\varpi_0)^*_\Z$}
\put(88,-14){$\Z[\wGLl]$}
\put(116,5){\vector(-3,-2){20}} 

\put(80,10){${\sim}$} 
\put(61,-12){${\sim}$} 
\put(70,5){$[V^\ep_{\bs\la}]\mapsto |\bs\la\rangle^*$} 
\put(32,-2){${ev_a^*}$}

\put(25,-17){\vector(0,-1){32}} 
\put(119,5){\vector(0,-1){32}}
\put(52,-16){\line(0,-1){11}}
\put(52,-11){\line(0,1){15}}
\put(92,-17){\vector(0,-1){32}}

\put(22,-35){}
\put(49,-21){}

\put(50,-35){\vector(-3,-2){20}}
\put(68,-31){\line(1,0){21}}
\put(95,-31){\vector(1,0){18}}
\put(40,-53){\vector(1,0){46}}
\put(15,-54){$\Rep_a \fu\agli$}
\put(43,-32){$\Rep\, \fu\gli$}
\put(117,-32){$L_l(\varpi_0)^*_\Z$}
\put(88,-54){$\Z[\wSLl]$}
\put(116,-35){\vector(-3,-2){20}} 

\put(80,-30){${\sim}$} 
\put(61,-52){${\sim}$} 
\put(80,-35){${}$} 
\put(37,-46){$ev_a^*$} 
\put(105,-46){$$}

\end{picture}
\caption{}\label{figure 2}
\end{figure}
\bigskip

We recall that $\Rep_a U_q\agli$ and $\Rep U_q\gli$ are defined in
Section \ref{ind limit def generic section}; $\Rep_a \ur\agli$ and
$\ur\gli$ in Section \ref{ind limit def root section} (see also
Section \ref{fin dim root}); $\Rep_a \fu\agli$ and $\fu\gli$ in
Section \ref{tensor product section} (see also Section \ref{fin
modules section}); the rings of functions $\Z[\wSLi]$, $\Z[\wGLl]$ and
$\Z[wSLl]$ in \secref{Minfty}, \secref{dual hall root} and \secref{new
sect}, respectively; the Fock spaces $\bigwedge(\varpi_0)$ and
$\bigwedge_l(\varpi_0)$ in Sections \ref{int section} and \ref{swl};
the module $L_l(\varpi_0)$ in Section \ref{swl}. Here we consider the
integral forms of their duals. The integral forms of
$\bigwedge(\varpi_0)^*$ and $\bigwedge_l(\varpi_0)^*$ are generated
over $\Z$ by the basis dual to $\{ |\bs\la\rangle \}$ and the integral
form of $L_l(\varpi_0)^*$ is the $\Z$--span of the projections of the
elements of this dual basis in $\Z[\wSLl]$.

The evaluation maps $ev_a$ are defined in Sections \ref{evaluation
section} and \ref{evaluation section root}, 
the specialization maps $\wh{\mc S}_l$ in
Section \ref{ind limit def root section} (see also \ref{spec
section}), the winding map $w^*$ in Section \ref{dual hall root}.

\begin{thm}\label{summary thm}
All squares of the diagram in Figure \ref{figure 2} are commutative.
\end{thm}

\begin{proof}
The upper floor is commutative by Lemma \ref{Hayashi} and Proposition
\ref{level one}. 
The middle floor is commutative by Proposition \ref{level one root}.
The bottom floor is commutative by Proposition \ref{level one root
factor}.

The upper left face is commutative because the 
evaluation homomorphism for $\ur\agln$ is by definition 
induced by the evaluation homomorphism for $\uqr\agln$. 
The bottom left face is commutative because the 
homomorphism $ev^*_a$ for $\Rep_a\fu\agln$ is by definition 
induced by the homomorphism $ev^*_a$ for $\Rep_a\ur\agln$.

The back upper face is commutative by definitions of all participating
maps. The back lower face is commutative because all participating
maps are maps of $U_\Z\wh{\mathfrak{sl}}_l^-$ modules, and uniqueness
of such maps.

The front upper face is commutative by Theorem \ref{repa to uminus
root}, part (2). The front bottom face is commutative by Theorem
\ref{fr ort sl}.

The upper right face is commutative by definitions of all
participating maps. The lower right face is commutative by Lemma
\ref{L and Lambda}.
\end{proof}

The diagram in Figure \ref{figure 2} 
corresponds to the vacuum module of level 1.
The same commutative diagrams with obvious changes hold for all
irreducible integrable $\wh{\mathfrak{sl}}_l$--modules (see Propositions
\ref{many gl to agl lemma}, \ref{many gl to agl root factor lemma}).

\end{document}